\documentclass[final]{siamart190516}
\usepackage{booktabs}
\usepackage[mathscr]{euscript}
\usepackage{color}
\usepackage{optidef}
\usepackage[algo2e,ruled,vlined]{algorithm2e}
\usepackage{float}
\usepackage{graphicx}
\usepackage{tikz,pgfplots}
\usepackage{upgreek}
\usepackage{multirow}
\usepackage{yfonts}
\usepackage{mathtools}
\usepackage[normalem]{ulem}
\usepackage{latexsym}
\usepackage{url}
\usepackage{amsmath,amssymb,amsbsy,amsfonts}
\usepackage{enumitem}
\usepackage{mathrsfs}
\usepackage[export]{adjustbox}

\definecolor{utorange}{rgb}{0.8,0.33,0.}
\definecolor{themec}{RGB}{51,108,121}
\definecolor{darkred}{rgb}{.6,.1,.1}
\definecolor{darkblue}{rgb}{.1,.1,.9}
\definecolor{greenback}{rgb}{.19,.94,.13}
\definecolor{orange}{rgb}{.76,.39,.13}
\definecolor{grass}{rgb}{.19,.64,.13}
\definecolor{sierp}{RGB}{209,28,209}
\definecolor{bgorange}{rgb}{1.,.95,.78}
\definecolor{grassgreen}{RGB}{92,135,39}
\definecolor{thinbox}{rgb}{.7,.8,1.}

\renewcommand{\vec}[1]{{\mathchoice
                     {\mbox{\boldmath$\displaystyle{#1}$}}
                     {\mbox{\boldmath$\textstyle{#1}$}}
                     {\mbox{\boldmath$\scriptstyle{#1}$}}
                     {\mbox{\boldmath$\scriptscriptstyle{#1}$}}}}

\newcommand{\eps}{\varepsilon}

\newcommand\restr[2]{{
  \left.\kern-\nulldelimiterspace 
  {#1}\vphantom{\big|} \right|_{#2}}}

\DeclareMathOperator*{\argmin}{argmin} 

\newcommand{\R}{\mathbb{R}}

\newcommand{\B}{\mathcal{B}}

\newcommand{\J}{\mathcal{J}}
\newcommand{\loss}{\mathcal{L}}

\newcommand{\Reg}{\mathcal{R}}

\newcommand{\iparspace}{\mathcal{M}}
\newcommand{\ipar}{m}
\newcommand{\iparpr}{m_{\text{pr}}}
\newcommand{\iparmap}{m_{\scriptscriptstyle\text{MAP}}}

\newcommand{\dpar}{\vec{m}}

\newcommand{\istatespace}{\mathcal{U}}
\newcommand{\istate}{u}
\newcommand{\dstate}{\vec{u}}

\newcommand{\iadj}{p}
\newcommand{\dadj}{\vec{p}}

\newcommand{\obsspace}{\mathcal{D}}

\newcommand{\noise}{\bs e}
\newcommand{\data}{ \bs d }

\newcommand{\ff}{\mathcal{F}}

\newcommand{\ipzspace}{\mathcal{Z}}
\newcommand{\ipz}{z}

\newcommand{\dpz}{\vec{z}}

\newcommand{\imul}{y}

\newcommand{\gbf}[1]{\boldsymbol{#1}}
\newcommand{\bs}[1]{\ensuremath{\boldsymbol{#1}}}

\renewcommand{\vec}[1]{\gbf{#1}}

\usepackage{xspace}

\newcommand{\LI}{\mathcal{L}}
\renewcommand{\H}{\mathcal{H}}

\def\addresswustl{Electrical \& Systems Engineering,
  Washington University in St Louis, St. Louis, MO, USA}
\usepackage{setspace}

\headers{Consensus ADMM for Inverse Problems Governed by Multiple PDEs}{L. Lozenski and U. Villa}

\title{Consensus ADMM for Inverse Problems Governed by Multiple PDE Models\thanks{
\funding{This work was funded by the National Science Foundation under Grant No ACI-1550593.}}}

\author{Luke Lozenski\thanks{\addresswustl}
\and Umberto Villa\thanks{\addresswustl 
  (\email{uvilla@wustl.edu}).}}

\pgfplotsset{compat = 1.16}
\begin{document}
\maketitle
\begin{abstract}
The Alternating Direction Method of Multipliers (ADMM) provides a natural way of
solving inverse problems with multiple partial differential equations (PDE) forward
models and nonsmooth regularization. ADMM allows splitting these large-scale 
inverse problems into smaller, simpler sub-problems, for which computationally
efficient solvers are available. In particular, we apply large-scale second-order
optimization methods to solve the fully-decoupled Tikhonov regularized inverse
problems stemming from each PDE forward model. We use fast proximal methods to
handle the nonsmooth regularization term. In this work, we discuss several
adaptations (such as the choice of the consensus norm) needed to maintain
consistency with the underlining infinite-dimensional problem. We present two imaging applications inspired
by electrical impedance tomography and quantitative photoacoustic tomography to demonstrate the proposed method's effectiveness. 
\end{abstract}



\section{Introduction}
Partial differential equations (PDEs) are used in various fields to model complex systems adhering to physical principles. However, PDEs often depend on unknown or uncertain parameters that cannot be measured directly \cite{Reyes15}, such cases can be framed as a type of inverse problem. An inverse problem has the goal of estimating a parameter field given a set of possibly noisy data. This data is related to the parameter via a \emph{parameter-to-observable map}. We will focus on this map being the solution to a PDE model. The inverse problem can then be formulated as an infinite-dimensional optimization problem with the objective being to minimize some distance between the measured data and solution to the PDE model with a regularization term on the parameter.

The main issue that separates an inverse problem from other types of optimization problems is the concept of ill-posedness \cite{Bui-ThanhGhattasMartinEtal13}. Ill-posedness means that solutions to these PDE relationships can face high-sensitivity to noise or small perturbations in the data creating larger perturbations in solutions or non-uniqueness; two distinct solutions could arise from one state. To account for noise, we then add a regularization term that enforces desirable results. Similarly, the proper regularization function can be used to fix the problem of non-uniqueness so that the objective is minimized around the correct solution. The regularization is decided by choosing what properties to enforce and can thus be a broad class of functions. 

 The infinite-dimensional nature of inverse problems means that we have to utilize a discretization process for numerical solutions at some point in the solution process.  To solve this issue we implement an \emph{optimize-then-discretize}  approach. This means we first derive the optimality conditions and derivatives in a continuous setting then implement a discretization process for numerical results \cite{Reyes15}. In this work, we use the finite element method (FEM) for discretization. FEM is the process in which a finite mesh approximates a continuous domain.  Functions over the domain are approximated by the span of a set of basis functions over the mesh; usually, piecewise polynomial \cite{ErnGuermond04,OdenReddy76,StrangFix88}. This numerical solution can be made arbitrarily accurate with a sufficiently fine mesh and a wider array of basis functions. However, this increased accuracy creates the trade-off of the problem becoming very large in scale.
 
In solving inverse problems of this form, we will need to apply regularization to account for noise and ill-posedness of the problem. In many cases, it is advantageous to apply nonsmooth regularization to reinforce desired behaviors and prior expectations. One such example of a useful nonsmooth regularization is total variation. Total variation has similar behaviors to regularization on the $L^2$ norm of the gradient but places a higher penalty on smaller values and a lower penalty on larger values of the gradient. Total variation allows the reconstructed parameter to vary more but still be primarily piecewise constant\cite{RudinOsherFatemi92}. However, using nonsmooth regularization is problematic because the methods that best solve PDE constrained optimization problems, INCG, require well-defined derivatives.

One method of solving problems of this form involves the Alternating Direction Method of Multipliers(ADMM).  The ADMM is a proximal point algorithm that is well suited for distributed convex optimization. The method was first formulated in the 1970s with works such as \cite{Rockafellar76} and splits large scale problems into multiple convex subproblems. Its potential for a variety of large scale problems, including machine learning, has been widely detailed in such works as \cite{BoydParikhChuetal10}, which also provides a rigorous analysis of the ADMM.  The authors in \cite{WangYinZeng18} further demonstrated that ADMM is well suited in a general optimization setting and is globally convergent for nonsmooth and nonconvex problems. With the development of ADMM focusing on large-scale problems and nonsmooth regularization, it is natural to have applications to image reconstruction and inverse problems governed by PDE forward models. The first application of ADMM for optimization involving PDEs was demonstrated in \cite{GabayMercier76} and applied to solve several continuum mechanics problems. In \cite{WangYangYinetal07, AlfonsoBioucas-DiasFigueiredo11} ADMM is applied to classical image reconstruction problems with a linear imaging operator, such as blurring and denoising operators, with total variation regularization. In \cite{FungRuthotto19} this methodology is further developed to explore inverse problems with PDE constraints with linear forward models and a version of ADMM with a weighted consensus relationship to increase progress made in early iterations. Similarly, \cite{ZandSiahkoohiMalcolmetal20} implements ADMM for solving an inverse problem related to seismology and demonstrates the compatibility of ADMM with nonsmooth regularization terms, including $L^1 $ regularization and second-order total variation.

This paper presents an application of the alternating direction method of multipliers (ADMM) for solving inverse problems governed by multiple PDE forward problems and nonsmooth regularization. We demonstrate ADMM's natural ability to handle multiple PDE forward problems and nonsmooth regularization functionals by splitting this large-scale problem into subproblems for which efficient solution methods are available. The specific contribution of our work is the following. We demonstrate the effectiveness of ADMM for solving inverse problems governed by multiple PDE models. In particular, we will demonstrate this framework with PDE relationships leading to nonlinear parameter to observable maps. Our framework will also utilize a consensus ADMM equipped with an arbitrary inner product for equality enforcement. We will also demonstrate how using different inner products and norms can lead to numerical stability for solving variational problems. Throughout, we will quantitatively demonstrate how ADMM reduces computational costs for large-scale problems compared to traditional methods while still achieving sufficient accuracy.  

The remainder of the paper is structured as follows. Section \ref{sec:ADMM} provides a brief overview of ADMM. In particular, we recall the scaled formulation of the ADMM algorithms, adaptive weights, and robust stopping criteria. In Section \ref{sec:PDE_CIP}, we provide a theoretical basis for PDE-constrained inverse problems and how the scaled-ADMM can naturally be applied to such problems. In Section \ref{sec:iNCG} we provide a general outline of deterministic inverse problems and a general optimization method, inexact Newton Conjugate Gradient, for inverse problems. In Section\ref{sec:PDE_CIP}, we will also introduce two algorithms for applying ADMM to large inverse problems. In Section \ref{sec:eit}, we look at a model problem related to electrical impedance tomography and perform numerical studies to show the novelty and merit in using ADMM for various problem sizes with multiple PDE models. In these numerical studies, we also show the effect of implementing a modified consensus norm, solving for local inversion parameters inexactly. In Section \ref{sec:qpact} we demonstrate the application of ADMM to an inverse problem found in quantitative photoacoustic tomography.

\section{The Alternating Direction Method of Multipliers(ADMM)}\label{sec:ADMM}

The \newline ADMM  is an algorithm for solving large scale constrained optimization problems whose objective involves the sum of a well behaved twice differential term and another term that may be nonsmooth. The ADMM splits this large problem into separate smaller subproblems, one of which will involve the nonsmooth objective term. These subproblems are solved iteratively, and equality between their solutions is reinforced with a consensus term.

In this section, we recall the various formulation of ADMM and some commonly used heuristic to accelerate ADMM convergence and robust stopping criteria.
The ADMM is part of a class of algorithms, called proximal point algorithms \cite{Rockafellar76}, which require very little to guarantee convergence.

We consider the following minimization problem
\begin{equation}\label{eq:minimization_problem}
    \min_{\ipar \in \iparspace} \J(\ipar) := \frac{1}{q}\sum_{i=1}^q \loss_i(\ipar) + \Reg(\ipar),
\end{equation}
where $\ipar$ is the sought after parameter in possibly infinite dimensional Hilbert space $\iparspace$, and $q>0$. The functionals $\loss_i: \iparspace \mapsto \mathbb{R}$ are assumed to be smooth (twice differentiable) and expensive to evaluate. The functional $\mathcal{R}: \iparspace \mapsto \mathbb{R}$ is assumed convex and non-smooth. Furthermore, we assume that computationally efficient methods are available to solve proximal problems stemming from $\mathcal{R}$.

For ease of notation, we assume $q=1$ in the presentation of the ADMM algorithm below. The general case $q\geq 1$ is presented in Section \ref{sec:PDE_CIP} in the contest of inverse problems governed by partial differential equation forward models.

The ADMM introduces a new variable $z \in \iparspace$ and changes the optimization problem \ref{eq:minimization_problem} to the form \eqref{eqn:general_opt}.

\begin{equation}\label{eqn:general_opt}
\begin{array}{cc}
    \displaystyle{\min_{\ipar,\ipz\in\iparspace,}}  & \loss(\ipar)+\Reg(\ipz),    \\
     s.t. & \ipar - \ipz = 0 
\end{array}
\end{equation}

For a general equality constrained optimization problem given by \eqref{eqn:general_opt} we first form the augmented Lagrangian for some $\rho > 0$ by  \eqref{eqn:aug_lag}. 
\begin{equation}\label{eqn:aug_lag}
    L_\rho (\ipar,\ipz,\imul) = \loss(\ipar) + \Reg(\ipz) + \langle \imul ,\ipar - \ipz \rangle + \frac{\rho}{2}||\ipar - \ipz||^2
\end{equation}

Where $\imul \in \iparspace$ is a Lagrange multiplier for the constraint $\ipar = \ipz$. Here $\langle \cdot, \cdot \rangle$ denotes the inner product and $||\cdot||$ is the norm induced by the inner product. 

The augmented Lagrangian has a few noteworthy properties. First, if the Lagrangian has a unique stationary point, then the augmented Lagrangian will have the same unique stationary point \cite{BirosGhattas03}. This invariance on the stationary point means optimizing using the augmented Lagrangian will result in the same outcome as the regular Lagrangian, and a choice of $\rho$ does not affect the outcome. 

Second, the added quadratic term for penalizing the constraint ensures convergence of the ADMM iterative process \cite{WangYinZeng18}. The augmented Lagrangian is also more desirable than other penalization or barrier functions because it will create well-behaved steps in the iterative process \cite{BirosGhattas03}.

Third, the lack of requirements on $\rho$ means that one's choice of $\rho$ can tuned to accelerate the convergence of the ADMM. A larger $\rho$ will make the variables $\ipar$ and $\ipz$ more accurately agree at every iteration. A smaller $\rho$ will make the $\ipar$ and $\ipz$ quickly approximate the optimal values of each subproblem but lack complete agreement. Choosing the optimal $\rho$ is then a trade-off between these two qualities. Algorithm 1 summarizes the ADMM.\\
\begin{algorithm2e}[H]
\SetAlgoLined
 Begin with starting points $(\ipar^0, \ipz^0, \imul^0)$ \\
 \While{While convergence criterion is not met}{
  $\ipar^{k+1} = \argmin_{\ipar} L_\rho (\ipar, \ipz^k, \imul^k)$\\
  $\ipz^{k+1} = \argmin_{\ipz} L_\rho (\ipar^{k+1}, \ipz, \imul^k)$\\
  $\imul^{k+1} = \imul^k + \rho(\ipar^{k+1} - \ipz^{k+1})$
 }
 \caption{The General ADMM}
\end{algorithm2e}
If we scale $\imul$ by a factor of $\frac{1}{\rho}$ and do a least squares completion then we can transform the ADMM to the scaled ADMM given by Algorithm 2\cite{BoydParikhChuetal10}. This form is useful because the problem is now only in terms of the objective functions and quadratic term. Furthermore, in the scaled form of ADMM the update of $\imul$ is independent of $\rho$.\\ 
\begin{algorithm2e}[H]
\SetAlgoLined
 Begin with starting points $(\ipar^0, \ipz^0, \imul^0)$ \\
 \While{While convergence criterion is not met}{
  $\ipar^{k+1} = \argmin_{\ipar} \loss(\ipar)+\rho ||\ipar-\ipz^k+\imul^k||^2$\\
  $\ipz^{k+1} = \argmin_{\ipz} \Reg(\ipz)+ \rho ||\ipar^{k+1}-\ipz+\imul^k||^2$\\
  $\imul^{k+1} = \imul^k + \ipar^{k+1} - \ipz^{k+1}$\\
  
 }
 \caption{The Scaled ADMM}
\end{algorithm2e}

If the Lagrangian $L_0$ has a unique saddle point at $(\ipar^*,\ipz^*,\imul^*)$ and $\loss$ and $\Reg$ are proper closed and convex functions then as $k\rightarrow \infty$ then $(\ipar^k,\ipz^k,\imul^k)\rightarrow (\ipar^*,\ipz^*,\imul^*)$ \cite{GabayMercier76}. That being said, ADMM will only reach moderate accuracy in a few iterations and requires many following iterations for high-precision convergence\cite{BoydParikhChuetal10, FungRuthotto19}. Luckily in many applications, including those considered here, this is sufficient. The ADMM is so useful for large-scale problems because it splits these problems into multiple sub-problems updated iteratively. This splitting means we only require the resources to solve these smaller problems, reducing memory requirements, problem complexity, and processing power. The ADMM is designed to be entirely parallelizable and only requires communication to the consensus from each instance variable once per update. Therefore it requires significantly less communication time compared with other parallelized solutions \cite{ZandSiahkoohiMalcolmetal20}.

\subsection{Adaptive $\rho$ choice}
Since the saddle point of \eqref{eqn:aug_lag} is independent of $\rho$, one can implement an adaptive choice of $\rho$ depending on the disagreement between the parameters and the rate at which they change.

To make this choice, two types of residuals are used as a measure of convergence following the framework outline in \cite{Wohlberg17}. We denote the \emph{primal residual} and the \emph{dual residual} at the $k$-th iteration respectively as $$r^k = \ipar^k - \ipz^k, \  s^k = \rho(\ipz^{k} - \ipz^{k-1}).$$
The primal residual serves as a measure for the agreement between $\ipar$ and $\ipz$. The primal residual being zero implies that the system has reached primal feasibility.
Meanwhile, the dual residual serves as a measure for the settling of $\ipz$ and it being zero indicates the system has reached dual feasibility\cite{Wohlberg17}. The primal residual being much greater than the dual residual means too much weight is given to the objective function over the agreement of $\ipar$ and $\ipz$. This issue is fixed by increasing the value of $\rho$. Similarly, if the dual residual is much greater than the primal residual, there is not enough weight on the consensus; thus, $\rho$ should be decreased.

To speed the convergence of this process the following heuristic can be used. For a fixed $\rho^k$ following the calculations of $\ipar^k, \ipz^k, r^k$ and $s^k$ we can the make the following choice of $\rho^{k+1}$ by 
\begin{equation}\label{eqn:rho_update}
     \rho^{k+1} = \begin{cases} \tau \rho^k & ||r^k|| > \mu || s^k|| \\ \rho^k/\tau & ||s^k|| > \mu ||r^k || \\ \rho^k & otherwise \end{cases} \end{equation}
     
Where $\mu,\tau >1$ are tunable parameters. Commonly these are chosen to by $\mu = 10$ and $\tau = 2$.

\subsection{Stopping criterion} Using the residuals we can also 
implement a gradient free stopping criterion adapted from 
\cite{Wohlberg17, BoydParikhChuetal10}. To do this we simply choose
a small absolute tolerance $\varepsilon_a > 0 $ and a small relative 
tolerance $\varepsilon_r >0 $. We then stop the iterative process at the 
first iteration when the following criterion is met. 
$$||r^k|| \leq \varepsilon_{abs} + \varepsilon_{rel}||\ipar^k||
 \textnormal{ and } ||s^k|| \leq\varepsilon_{abs} + \varepsilon_{rel}||\ipz^k||$$

 It should be noted that $||\cdot||$ is the norm for $\iparspace$.

\section{Inexact Newton conjugate gradient method for PDE-constrained inverse problems}\label{sec:iNCG}

This section presents a brief outline of the formulation of deterministic inverse problems and a specific method for their solution. We begin by outlining the notation that we will use relating to inverse problems and their formalism. Then we explain how to apply the Inexact Newton Conjugate Gradient(INCG) method to solve problems of this form. 

\subsection{Deterministic inverse problems} \label{det_inv_probs}

An inverse problem has the goal of reconstructing for a parameter $\ipar$ given a measurement $\data \in \obsspace$. Mathematically the forward model of an inverse problem can be expressed as \eqref{eqn:general_forward_model}.

\begin{equation}\label{eqn:general_forward_model}
    \data = \ff(\ipar) + \noise,
\end{equation}

 where $\ff:\iparspace \rightarrow \obsspace$ is a map from the parameter space to the observation space, and $\noise$ is random noise present in measurements. 
 The forward map $\ff$ is often a relationship that is not directly available; for example, we will be focusing on implicit PDE relationships. For direct evaluations we can then introduce a state variable $\istate\in \istatespace $ satisfying some relationship $r(\ipar,\istate) = 0$. The goal of this transformation is to create an explicit relationship  $\data = \B(\istate)$, where $\B: \istatespace \rightarrow \obsspace$ is an observation operator.

 Given $\data$, our goal is to reconstruct for $\ipar$ while adjusting for the presence of the noise $\noise$. This can be characterized as minimizing the cost functional 
 
 \begin{equation}\label{eqn:general_cost_f}
     \J(\ipar) := \LI(\ipar)  + \Reg(\ipar)\text{ where } \LI(\ipar) = \frac{1}{2} \| \ff(\ipar) - \data\|^2
 \end{equation}
Here $\Reg$ is a regularization function and $\LI$ is a data fidelity term that is minimized when $\ff(\ipar) = \data$. Substituting in $\ff(\ipar) = \B(\istate)$, 
the new goal of the inverse problem is to solve the minimization problem in

\begin{equation}\label{eqn:inverse_prob_opt}
\begin{array}{cc}
    \displaystyle{\min_{\ipar\in\iparspace, \istate \in \istatespace}}  & \J(\ipar) = \frac{1}{2} \| \B(\istate) - \data\|^2+ \Reg(\ipar)   \\
     s.t. &  r(\ipar,\istate) = 0
\end{array}
\end{equation}
With this problem now formulated as a constrained optimization problem we can express the Lagrangian in \eqref{eqn:variational_Lagrangian} in terms of $\ipar$, $\istate$ and an adjoint variable $\iadj \in \iparspace$. 

\begin{equation}\label{eqn:variational_Lagrangian}
    \mathscr{L}(\ipar, \istate, \iadj) = \frac{1}{2}||\B(\istate) - \data||^2 + \langle \iadj, r(\ipar,\istate) \rangle.
\end{equation}

With this Lagrangian formalism expressed, we can derive the infinite-dimensional analogs of the gradient and Hessian. Letting subscripts denote Gâteaux derivative, we can denote the gradient of the cost functional \eqref{eqn:general_cost_f} in a a direction $\Tilde{\ipar} \in \iparspace$ at a point $\ipar_0$ as

\begin{equation}\label{eqn:variational_gradient}
    (\mathcal{G}(\ipar_0),\Tilde{\ipar}) = (\Reg_\ipar(\ipar_0),\Tilde{\ipar}) + \langle \iadj_0,r_\ipar(\istate_0,\ipar_0)[\Tilde{\ipar}]\rangle.
\end{equation}

 Above $\istate_0$ is the solution of the forward problem

\begin{equation}\langle \label{eqn:forward_problem} \Tilde{\iadj},r(\istate_0,\ipar_0)\rangle = 0 , \  \forall \Tilde{\iadj}, \end{equation}
which is obtained by requiring variations of \eqref{eqn:variational_Lagrangian} to vainish for all directions $\Tilde{\iadj}.$ $\iadj_0$ is the solution to the adjoint problem

\begin{equation}\label{eqn:adjoint_problem}\langle\iadj_0, r_\istate(\istate_0,\ipar_0)[\Tilde{\istate}]\rangle + \langle \istate(\B(\istate_0) - \data, \B \Tilde{\istate}\rangle = 0, \ \forall \Tilde{\istate},\end{equation}

which is obtained by requiring variations of \eqref{eqn:variational_Lagrangian} to vanish for all directions $\Tilde{\istate}\in \istatespace$.

Similarly, to derive the Hessian action, we consider the second-order Lagrangian

\begin{equation}\label{eqn:second_order_hessian}\begin{array}{ll}
   \mathscr{L}^{\H}(\ipar,\istate, \iadj;\hat{\ipar},\hat \istate, \hat \iadj)  & = (\mathcal{G}(\ipar), \hat \ipar)  \\
     & +  \langle \hat{\iadj}, r(\istate,\ipar) \rangle \\
     & + \langle \iadj, r_\istate (\istate, \ipar)[\hat{\istate}]\rangle + \langle \B(\hat{\istate}), \B(\istate) - d\rangle .
\end{array}\end{equation}

The Hessian in a direction $\hat{m}\in \iparspace$ evaluated at $\ipar = \ipar_0$ is then the Gâteaux derivative of $\mathscr{L}^\H$ with respect to $\ipar$ and given as 

\begin{equation}\label{eqn:variational_Hessian} \begin{array}{lll} (\Tilde{\ipar},\H(\ipar_0)\hat{\ipar}) & = (\Tilde{\ipar}, \Reg_{\ipar \ipar}(\ipar_0)[\hat \ipar]) + (\iadj_0, r_{\ipar \ipar}(\istate_0, \ipar_0)[\Tilde{\ipar}, \hat{\ipar}]) & \\ 
 & + \langle \hat \iadj, r_\ipar(\istate_0, \ipar_0)[\Tilde{\ipar}]\rangle + \langle \iadj_0, r_{\istate \ipar})(\istate_0, \ipar_0)[\hat{\istate}, \Tilde{\ipar}]\rangle, & \forall \Tilde{\ipar} \in \iparspace
\end{array}\end{equation}

Where 
incremental state $\hat{\istate}$ and incremental adjoint $\hat{\iadj}$ solve the so-called incremental forward and incremental
adjoint problems, which are obtained by setting to zero variations of \eqref{eqn:second_order_hessian} with respect to $\iadj$ and $\istate$ respectively. 
Optimality conditions still hold in the infinite dimensional setting. This means that the cost functional in \eqref{eqn:general_cost_f} and \eqref{eqn:inverse_prob_opt} is minimized at a point when the gradient \eqref{eqn:variational_gradient} is identically zero for every $\Tilde{\ipar}$ and the Hessian is positive definite.  
\subsection{Inexact Newton Conjugate Gradient(INCG) for solving inverse problems}

With the infinite-dimensional derivatives derived in Section \ref{det_inv_probs} we can apply traditional minimization algorithms, including gradient descent and Newton descent. Here we will outline the application of Inexact Newton Conjugate Gradient(INCG) for solving an inverse problem. Now we proceed with discretizing the problem for a numerical solution. It is then necessary to note that the gradient, Hessian, and optimality conditions are evaluated as discretized statements of variational problems instead of first discretizing the problem and then treating it as a finite-dimensional optimization problem. The INCG algorithm is shown in Algorithm \ref{alg:iNCG}.

\begin{algorithm2e}[t]
\SetAlgoLined
Start with $i = 0$. \\
Given $\dpar_0$ solve the forward problem \eqref{eqn:forward_problem} to obtain $\dstate_0$.\\
Given $\dpar_0,\dstate_0$ compute the cost functional $\J_0$ using \eqref{eqn:general_cost_f}.\\
 
 \While{$i < $max\_iter }{
 Given $\dpar_i,\dstate_i$ solve the adjoint problem \eqref{eqn:adjoint_problem} to obtain $\dadj_i$\\
 
Given $\dpar_i, \dstate_i, \dadj_i$ compute the gradient $\vec{g}_i$ using \eqref{eqn:variational_gradient}.\\
\If{$||\vec{g}_i|| \leq \tau$}{
    \textbf{break}
    }
Given $\dpar_i, \dstate_i,\dadj_i$ define a linear operator $\vec{H}_i$ implementing the Hessian action \eqref{eqn:variational_Hessian}.\\
Using Conjugate gradients, find a search direction $\hat{\dpar_i}$ such that 
\begin{equation}\label{eqn:eis_walk}
    ||\vec{H}_i\hat{\dpar_i} + \vec{g}_i|| \leq \eta_i||\vec{g}_i||, \ \eta_i = \left(\frac{||\vec{g}_i||}{||\vec{g}_0||}\right)^{1/2} \end{equation} \\
Set $j = 0, \alpha^{(0)}=1$
 \While{$j < $ max\_backtracking\_iter}{
  Set $\dpar^{(j)} = \dpar_i + \alpha^{(j)}\hat{\dpar_i}$\\
  Given $\dpar^{(j)}$ solve the forward problem \eqref{eqn:forward_problem} to obtain $\dstate^{(j)}$ \\ 
  
  Given $\dpar^{(j)}$ and $\dstate^{(j)}$ copute the cost $\J^{(j)}$ using \eqref{eqn:general_cost_f}\\
  \If{$\J^{(j)} < \J_i + \alpha^{(j)} c_{armijo}\vec{g}_i^T\hat{\dpar}_i$}{
    $\dpar_{i+1} \leftarrow \dpar^{(j)}, \ \J_{i+1}\leftarrow \J^{(j)}$ \\
    \textbf{break}}
    $\alpha^{(j+1)} \leftarrow \alpha^{(j)}/2, j \leftarrow j+1$ \\
    
  }
  $i \leftarrow i +1$
 }
\caption{The Inexact Newton Conjugate Gradient algorithm for solving inverse problems}
\label{alg:iNCG}
\end{algorithm2e}

Here \eqref{eqn:eis_walk} is known as the Eisenstat-Walker condition and results in desirably fast local convergence \cite{EisenstatWalker96}. This condition leads to superlinear convergence of Algorithm 3 while at the same time drastically reducing the number of necessary iterations to solve the Newton system.

\section{Application of ADMM to the solution of inverse problems governed by multiple PDE constraints} \label{sec:PDE_CIP}

In this section, we consider the minimization problem \eqref{eq:minimization_problem}, in the context of an infinite dimensional inverse problem with PDE forward problems. Then, $\ipar \in \iparspace$ belongs to some Sobolev space defined on a domain $\Omega \subset \mathbb{R}^d$ ($d=1, 2, 3$), the functionals $\mathcal{L}_i(\ipar)$ represent the smooth data fidelity terms, those evaluation involve the solution of the PDE, and $\mathcal{R}(\ipar)$ the regularization functional.
Specifically, we consider the following form of the data-fidelity term  
\begin{equation}
\mathcal{L}(\ipar) = \frac{1}{q}\sum_{i=1}^q\mathcal{L}_i(\ipar) :=  \frac{1}{2q}\sum_{i=1}^q\left\| \ff_i(\ipar) -\data_i \right\|^2,
\end{equation}
where $\data_i \in \obsspace$ ($i=1,\ldots,q$) represent the data, and $\ff_i: \iparspace \mapsto \obsspace$ is the parameter to observable map. For the applications we focus on $\ff_i$ will be the composition of a PDE solver and an observation operator. \cite{VillaPetraGhattas20}.

\subsection{ADMM with consensus equations}

In Section \ref{sec:ADMM}, we introduce ADMM for a single data fidelity term $\loss$(i.e. for $q=1$). Here we can generalize this to deal with multiple PDE-based forward models in a special version of the global consensus problem \cite{BoydParikhChuetal10}. If we are given $q$ data sets $\{\data_i\}_{i=1}^{q}$ and their corresponding parameter to observable maps $\{\ff_i:\iparspace \mapsto \obsspace \}_{i=1}^{q}$, our goal would normally be to solve for $\iparmap$ as in \eqref{eqn:inv_obj_multiple}.
\begin{equation}\label{eqn:inv_obj_multiple}
\iparmap = \argmin_{\ipar} \frac{1}{2q} \sum_{i=1}^q||\ff_i(\ipar) - \data_i||^2 + \Reg(\ipar) \end{equation} We can instead split this parameter for each model and data set and apply an equality constraint between its multiple instances.
\begin{equation}\label{eqn:multiple_equality}
\ipar_1 = \hdots = \ipar_q = \ipz.
\end{equation} 
This also requires the introduction of $q$ different multipliers $\{\imul_i\}_{i=1}^q \subset \iparspace$ for each of these equality relationships. $\imul_i$ will act as the Lagrange multiplier for the equality relationship $\ipar_i = \ipz$.
Using this, we can form the scaled augmented Lagrangian over all these variables and change our goal to solving \eqref{eqn:multiple_pde_constrained_opt}.

\begin{equation}\label{eqn:multiple_pde_constrained_opt}
\begin{array}{cc}
    \displaystyle{\min_{\ipar_i \in 
    \iparspace, \ipz \in \ipzspace}}  & 
    \frac{1}{2q}\displaystyle{\sum_{i=1}^q}||\ff_i(\ipar_i) - \data_i||^2 
    +\Reg(\ipz)   \\
     s.t. & \ipar_i - \ipz = 0 \ for \ i=1,\hdots, q 
\end{array}
\end{equation}

If we scale $\rho$ by a factor of $\frac{1}{q}$ then this problem will result in the augmented Lagrangian given by \eqref{eqn:ADMMobj_mult}.
\begin{equation}\label{eqn:ADMMobj_mult}
\begin{array}{ccl}
    L_\rho(\{\ipar_i\}_{i=1}^q, z, \{\imul_i\}_{i=1}^q)   & =  & \frac{1}{2q} \sum_{i=1}^q||\ff_i(\ipar_i) - \data_i||^2 \\
     & + & \Reg(\ipz)  + \frac{1}{q}\sum_{i=1}^q\langle\imul_i, \ipar_i - \ipz \rangle\\
     & + &\frac{\rho}{2q}\sum_{i=1}^q ||\ipar_i - \ipz ||^2
\end{array}
\end{equation}.

With this form the optimal argument for each $\ipar_i$ will be independent of all other parameters $\ipar_1,\hdots,\ipar_{i-1},\ipar_{i+1},\hdots, \ipar_q$. This now means that at each step of the scaled ADMM each $\ipar_i^k$ can be updated in the simplified form given by \eqref{eqn:ADMMobj_split}. \begin{equation}\label{eqn:ADMMobj_split}
    \ipar^{k+1}_i = \argmin_{\ipar_i}\frac{1}{2q} ||\ff_i(\ipar_i) - \data_i||^2 + \frac{\rho^k}{2q} ||\ipar_i - \ipz^k + \imul_i^k||^2.
\end{equation}
This greatly reduces the computational complexity of each update opposed to only having one parameter variable for every model and simplifies the regularization to always being a Tikhonov regularization for the update. Similarly consensus variable $\ipz$ will be updated as
\begin{equation}\label{eqn:consensus_update}
z^{k+1} = \argmin{z} \Reg(\ipz) + \frac{\rho}{2q}\sum_{i=1}^q ||\ipar_i - \ipz + \imul_i||^2.
\end{equation}

This update is entirely free of the terms related to the PDE, which means we can solve it with a broader class of optimization methods. Thus the update process can be described by Algorithm \ref{alg:scaled_ADMM_invProb}.

\begin{algorithm2e}[ht]
\SetAlgoLined
Let $q$ be the number of PDE relationships \\
 Begin with starting points $(\{\ipar_i^0\}_{i=1}^q, z^0, y^0)$  \\
 \While{While convergence criterion is not met, $k=1,\hdots$}{
 \For{$i = 1,\hdots,q$}{
  $m_i^{k+1}$ is updated as in \eqref{eqn:ADMMobj_split} \\
  }
  $z^{k+1}$ is updated as \eqref{eqn:consensus_update} \\
  \For{$i = 1,\hdots,q$ }{
  $\imul_i^{k+1} = \imul_i^k+ (\ipar_i^{k+1} - \ipz^{k+1})$  \\
  }
  Update $\rho^{k+1}$ following  \eqref{eqn:rho_update}
 }
\caption{The Scaled ADMM for parameter inversion with multiple PDE's}
\label{alg:scaled_ADMM_invProb}
\end{algorithm2e}

At this point, we observe that 

$$\frac{1}{q}\sum_{i=1}^q||\ipar_i + \imul_i- \ipz||^2 = ||\ipz||^2 - 2\langle \ipz, \frac{1}{q}\sum_{i=1}^q \ipar_i + \imul_i\rangle + \frac{1}{q}\sum_{i=1}^q ||\ipar_i + \imul_i||^2=$$

$$||\ipz||^2 - 2\langle \ipz, \bar{\ipar} + \bar{\imul}\rangle + ||\bar{\ipar} + \bar{\imul}||^2 +\frac{1}{q}\sum_{i=1}^q ||\ipar_i + \imul_i||^2 - ||\bar{\ipar} + \bar{\imul}||^2 = $$

$$ || \bar{\ipar}+\bar{\imul}-\ipz||^2+ \frac{1}{q}\sum_{i=1}^q ||\ipar_i + \imul_i||^2 - ||\bar{\ipar} + \bar{\imul}||^2 \textnormal{, where } \bar{\ipar} = \frac{1}{q}\sum_{i=1}^q \ipar_i, \bar{\imul} = \frac{1}{q}\sum_{i=1}^q \imul_i$$

 The term $\frac{1}{q}\sum_{i=1}^q ||\ipar_i + \imul_i||^2 - ||\bar{\ipar} + \bar{\imul}||^2$ is constant in $\ipz$, which means that the update for $\ipz$ is equivalent to
\begin{equation}\label{eqn:mean_based_consensus}
    z^{k+1} = \argmin_{z} \Reg(\ipz) + \frac{\rho}{2} ||\bar{\ipar}- \ipz + \bar{\imul}||^2.
\end{equation}
We can then implement a mean based approach that simplifies the optimization process for updating $\ipz$. This mean based approach is shown in Algorithm \ref{alg:scaled_ADMM_mean_based}.

\begin{algorithm}[t]
\SetAlgoLined
Let $q$ be the number of PDE relationships \\
 Begin with starting points $(\{\ipar_i^0\}_{i=1}^q, z^0, y^0)$  \\
 \While{While convergence criterion is not met, $k=1,\hdots$}{
 \For{$i = 1,\hdots,q$}{
  $\ipar_i^{k+1}$ is updated as in\eqref{eqn:ADMMobj_split}  \\
  }
  Set $\bar{\ipar} = \frac{1}{q}\sum_{i=1}^q\ipar^{k+1}_i$ and, $\bar{\imul} = \frac{1}{q}\sum_{i=1}^q\imul^{k+1}_i$\\ 
  $\ipz^{k+1}$ is updated as in \eqref{eqn:mean_based_consensus} \\
  \For{$i = 1,\hdots,q$ }{
  $\imul_i^{k+1} = \imul_i^k+ \frac{1}{q}(\ipar_i^{k+1} - \ipz^{k+1})$  \\
  }
  Update $\rho^{k+1}$ following  \eqref{eqn:rho_update}
 }
\caption{The Mean based Scaled ADMM for parameter inversion with multiple PDE's }
\label{alg:scaled_ADMM_mean_based}
\end{algorithm}
This process can further be expanded by splitting the parameter variables by spatially dependent subregions and implementing variable asynchronous weights for consensus update as demonstrated in \cite{FungRuthotto19}.

\section{Numerical studies: electrical impedance tomography problem}\label{sec:eit}

To demonstrate the effectiveness of ADMM, we consider a model problem based on electrical impedance tomography\cite{CheneyIsaacsonNewell99, BorceaGrayZhang03}. With this model problem, we will consider four different experiments. In the first, we will consider using the $H^1$ norm for ADMM consensus compared to the $L^2$ norm. In the second, we will demonstrate the effectiveness of using inexact updates to accelerate the global solution. In our third experiment, we will consider the computational cost of the ADMM compared with the monolithic approach on a discrete mesh at multiple refinements with a fixed number of PDE models. Our fourth experiment will analyze the computational cost of the ADMM compared with the monolithic approach with a varying number of PDE models on a mesh of fixed size. Here the monolithic approach references. The term monolithic refers to solving a single large problem without breaking it into smaller sub-problems. The monolithic approach then implemented a traditional INCG descent method, found in Algorithm 3, to directly optimize the cost functional of the inverse problem. With each ADMM inversion, we performed the same inversion using the monolithic approach to compare accuracy, solution time, and computational cost.

To calculate the instance parameter updates and the monolithic approach, we applied the INCG solver found in hIPPYlib, an extensible software framework for large-scale inverse problems governed by PDEs \cite{VillaPetraGhattas20}.
The update of $\dpz$ in Algorithm 5 was calculated using the PETScTAOSolver built into Fenics\cite{LoggMardalWells12}, a comprehensive library designed for numerical solutions for PDEs.
\subsection{Formulating electrical impedance tomography in the continuous setting}
Electrical impedance tomography is an imaging modality that relies on inputting an electrical current to a portion of the domain boundary and measuring the resulting electric potential on the rest of the domain's boundary. The electric potential is dependent on the conductivity of the material throughout the domain. 
In this example, we consider a compact domain $\Omega \in \R^2$ representing the object to be imaged of and let $\iparspace := H^1(\Omega)$ be the Sobolev space of square-integrable functions with square-integrable gradients.
The data fidelity terms  $\mathcal{L}_i$ in \eqref{eq:minimization_problem} have the form:
\begin{equation*}
    \mathcal{L}_i(\ipar) = \frac{1}{2}\int_{\Gamma_{i} } (\istate_i - \data_i)^2 d\mathbf{s},
\end{equation*}

\noindent where $\Gamma_i \subset \partial\Omega$ is portion of the boundary where the state variable (electric potential) $\istate_i$ is measured. The potential $\istate_i$ solves the electrostatic Maxwell equation:

\begin{equation}
    \begin{cases}
    -\nabla \cdot e^m \nabla u_i = 0 & x \in \Omega \\ 
    
    \frac{\partial}{\partial \eta}u_i= g_i & x \in \Gamma^i_N \\ 
    
    u_i = 0 & x \in \Gamma_D^i
    \end{cases}\label{eqn:eit_pde}
\end{equation}

Here $\sigma:= e^\ipar$ is the conductivity of the domain, and $u_i$ is the electric potential resulting from introducing the current $g_i$. $\Gamma_N^i$ denotes the Neumann boundary corresponding to the current injected, and $\Gamma_D^i$ is the Dirichlet boundary corresponding to the electrical ground. $\partial \Omega = \Gamma_N^i \cup \Gamma_D^i$
Suppose then that we perform $q$ measurements with $q$ different currents, resulting in $d_i$.

Our goal is then to find a minimize \eqref{eqn:eit_obj} satisfying \eqref{eqn:eit_pde}.

\begin{equation}\label{eqn:eit_obj}
    \frac{1}{2 q}\sum_{i=1}^q \mathcal{L}_i(\ipar) + \mathcal{R}(\ipar), 
\end{equation}
where $\mathcal{R}(\ipar)$ is a combination of Total Variation and $L^2(\Omega)$ regularization defined as 
\begin{equation}
\mathcal{R}(\ipar) = \alpha_{TV}\int_{\Omega}|\nabla (\ipar - \iparpr) |_\eps d\boldsymbol{x} + \frac{\alpha_{TK}}{2}\int_{\Omega} (\ipar - \iparpr)^2d\boldsymbol{x},
\end{equation}
where $\iparpr \in \iparspace$ is a reference value for the inversion parameter, $\alpha_{TV}, \alpha_{TK} > 0$ are the regularization parameters. Finally, $$|\nabla (\ipar - \iparpr) |_\eps = \sqrt{ (\nabla (\ipar - \iparpr))^T(\nabla (\ipar - \iparpr) + \eps}$$ is a smooth approximation to make the TV functional differentiable. The parameter $\eps >0$ controls the smoothness of the functional.

\subsection{Discretization}
The unit disc was selected as our domain of interest $\Omega$. This continuous domain was then discretized with a uniform mesh with triangular elements. On this mesh, we chose our set of basis functions to be continuous piecewise linear finite polynomials for both $\iparspace$ and $\istatespace$.  On our coarsest mesh, we then had 8044 degrees of freedom on the parameter and state variables. However, for our experiment in \ref{ref}, we will perform multiple mesh refinements. This will then result in parameter and state variables with  $8044, 31816, 71280,$ and  $126428$ degrees of freedom.

\subsection{Ground truth, synthetic data}
\begin{figure}[tbh]
    \centering
    \resizebox{\textwidth}{!}{
    \includegraphics[scale = 1]{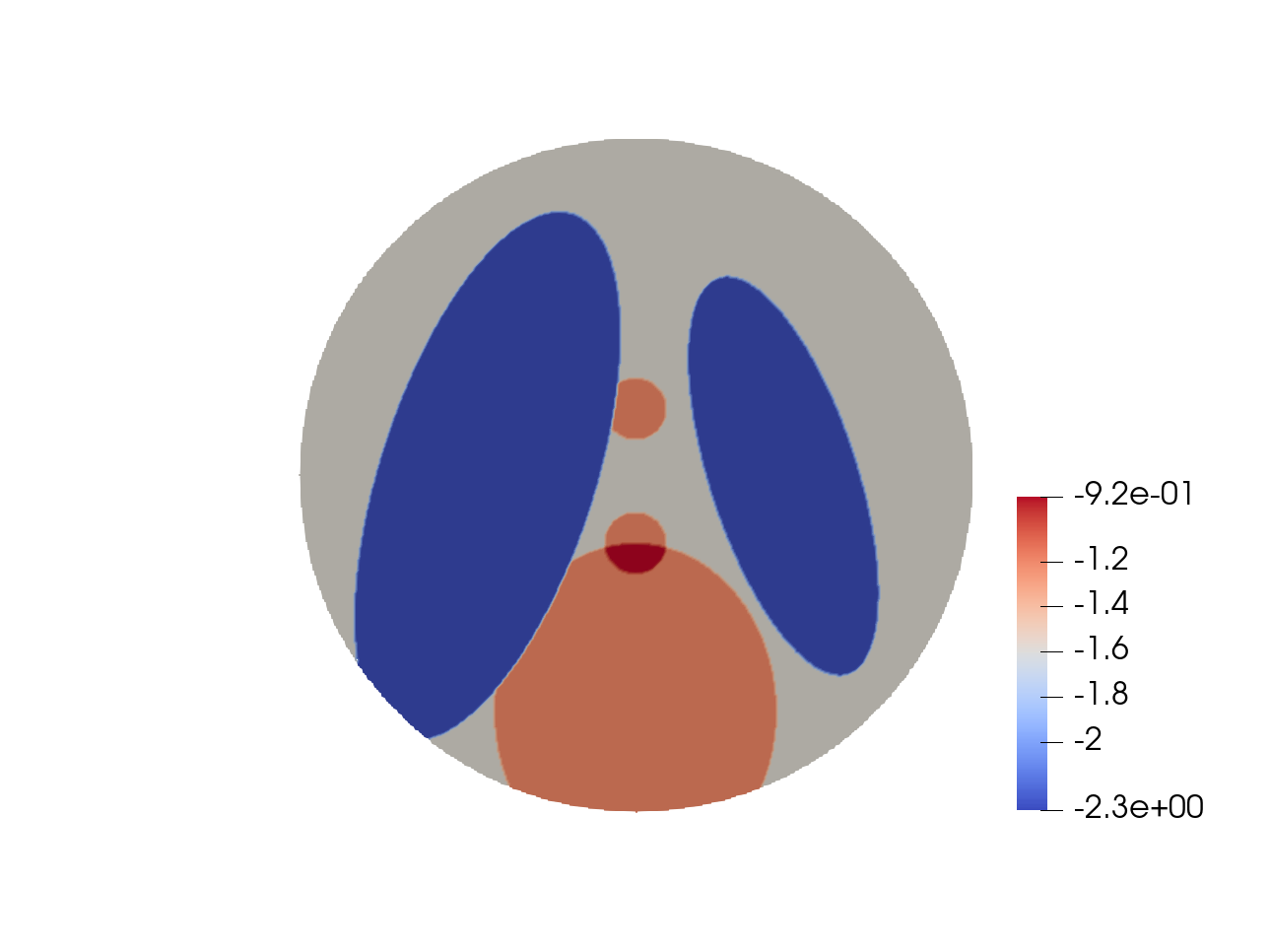}
    \includegraphics[ trim={0 0 10cm 0},clip]{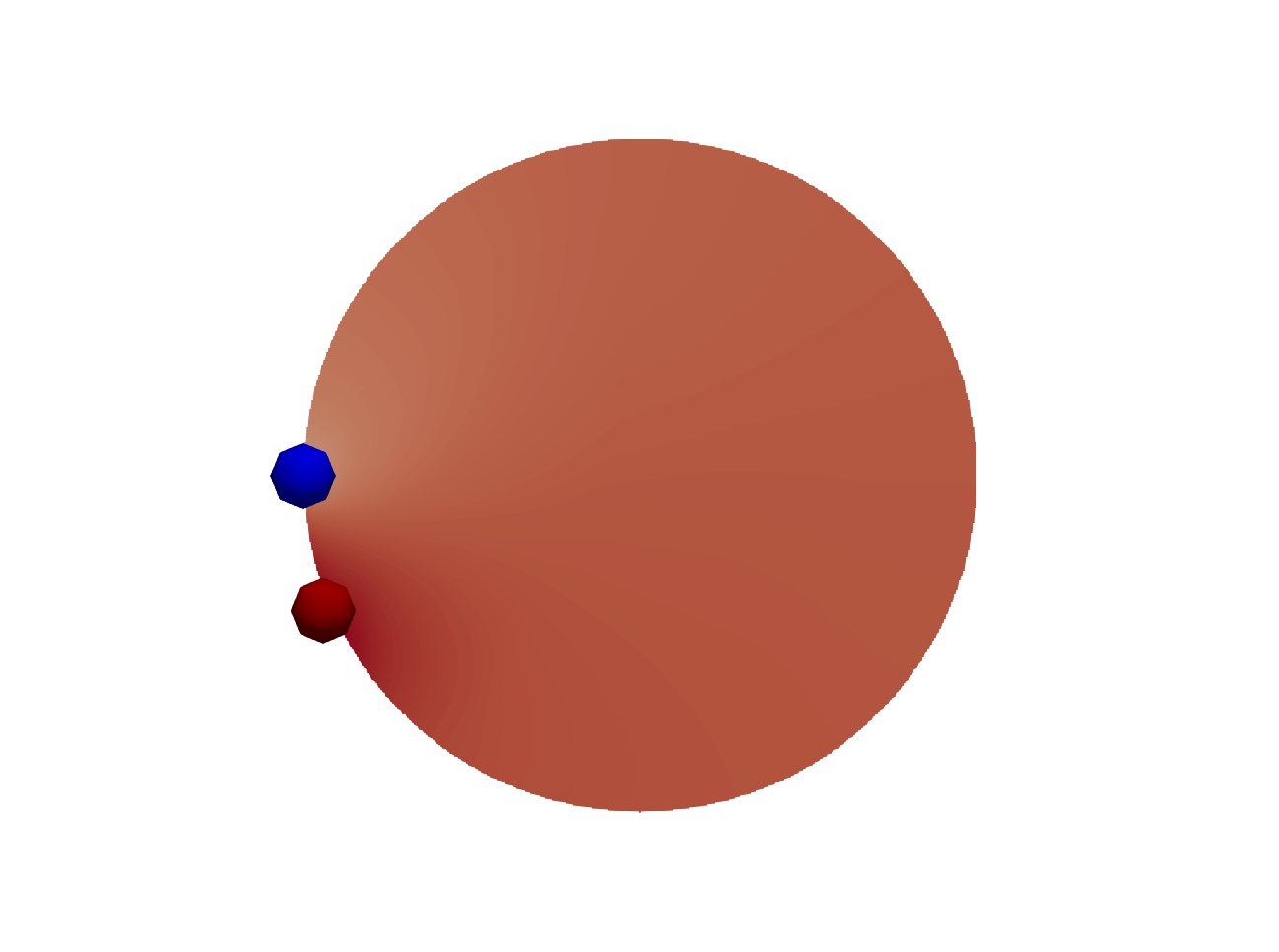}
    \includegraphics[ trim={5cm 0 5cm 0},clip]{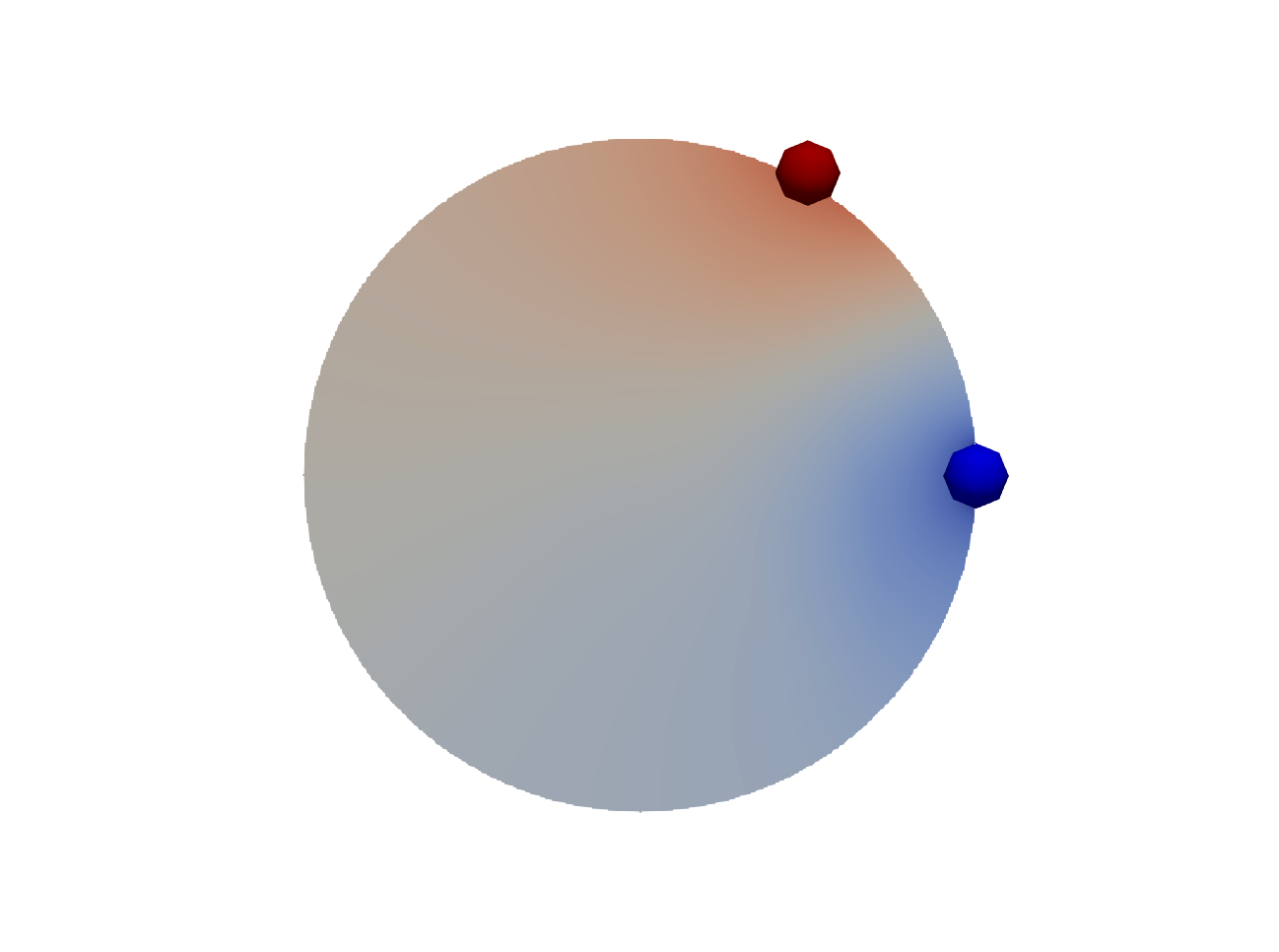}
    \includegraphics[ trim={5cm 0 0  0},clip]{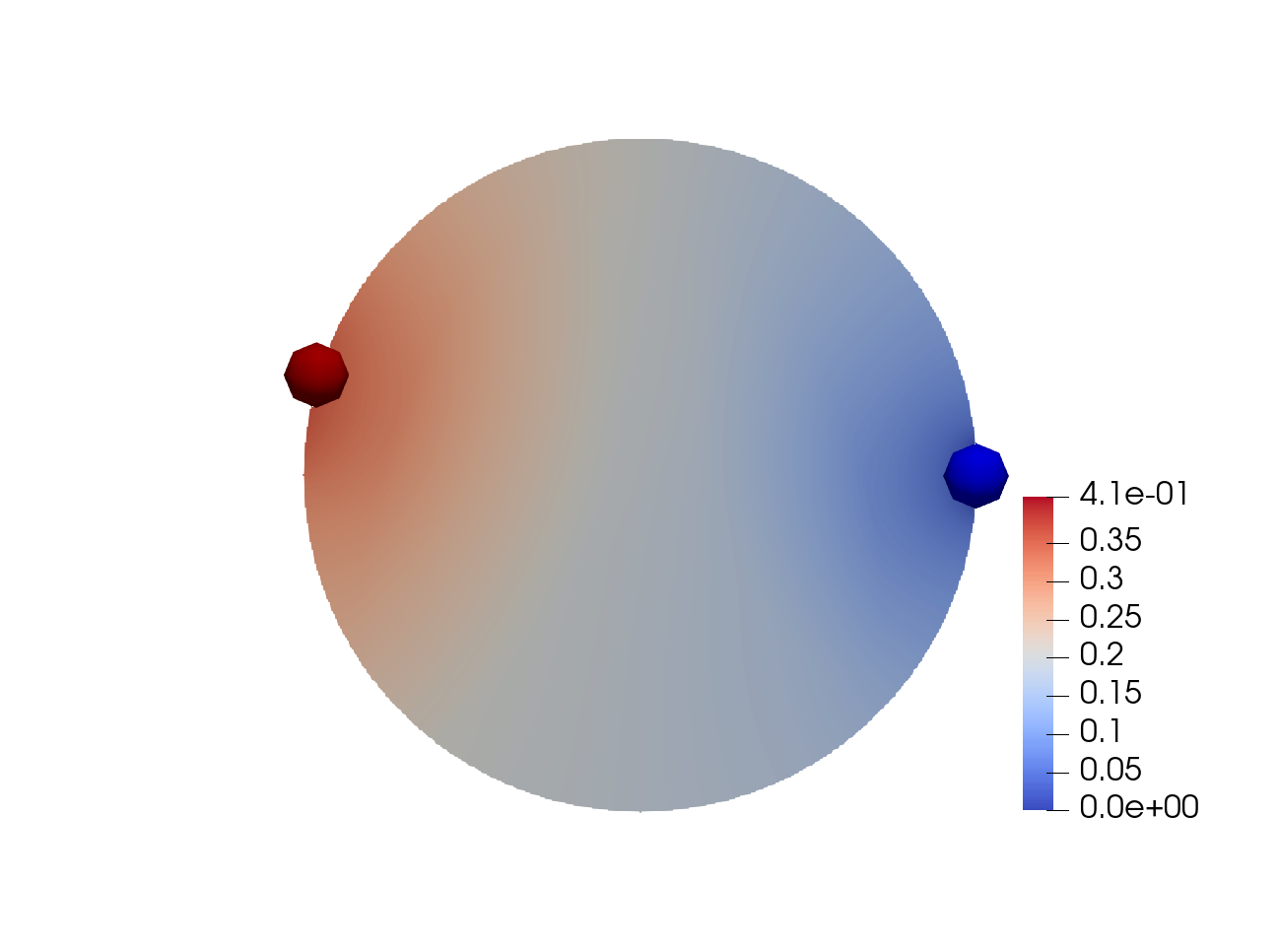}}
    \caption{True parameter (left) and true states  corresponding to sources 1,11,16 (right) for EIT problem with $q=16$  }
    \label{fig:eit_true_param_states}
\end{figure}

The true parameter used was a modified Shepp-Logan Phantom on the unit circle displayed in Figure \ref{fig:eit_true_param_states} (left).

Next we let the incident current $g_i$ be given by \eqref{eqn:def_current}

\begin{equation}\label{eqn:def_current}
    g_i(\theta) = \gamma\exp(-\beta(\theta - \theta_i)^2)
\end{equation}
where $\theta$ is the angle a point on $\partial \Omega$, $\theta_i$ dictate the position of the electrical source are dispersed evenly along the boundary, and $\gamma,\beta >0$ are constants dictating the amplitude and decay of the source. We will use $\gamma = 0.1$ and $\beta =10$. The Dirichlet boundary $\Gamma^D_i$, which acted as electrical ground, was chosen to be a single point on the boundary. The rest of the boundary was considered to be the Neumann boundary, $\Gamma_N^i = \partial \Omega \setminus \Gamma_D^i$.

The true states for models $1,11,16$ for the $q=16$ case are displayed in Figure \ref{fig:eit_true_param_states} (right).  The electrical source is highlighted with a red sphere, and the electrical ground is highlighted with a blue sphere.

The states were then perturbed by random noise with a standard deviation of $0.01*||u||_\infty,$ or one hundredth of the maximum state. To account for this noise we performed reconstructions with regularization decided by $\alpha_{TV} = 0.1$, $\varepsilon = 10^{-4}$ and $\alpha_{TK} = 0.01$. 

Here we will also use the adaptive $\rho$ scheme with $\mu = 2$ and $\tau = 3 $. The termination constants are set to $\varepsilon_a= 10^{-5}$ and $\varepsilon_r =2\cdot 10^{-2}.$

\subsection{Reconstruction with $L^2$ norm and $H^1$ norms}
In this section we compare reconstruction performed with the norms given by the $L^2(\Omega)$ and $H^1(\Omega)$ norms. $L^2(\Omega)$ is the space of square intolerable functions with a norm given by the square root of this integral. The Sobolev space$H^1(\Omega)$ is the space of functions such that the norm given in \eqref{eqn:H1_norm} is bounded.
\begin{equation}\label{eqn:H1_norm}
    ||\ipar||_{H^1(\Omega)} = \left(\int_\Omega 
    ||\ipar||^2 +  ||\nabla \ipar||^2 dx \right)^{1/2}
\end{equation}
For these reconstructions, we fix $q=16$. The ADMM solver utilized a maximum of 10 global iterations and used an INCG solver with a maximum of 10 iterations to find a more precise estimate of individual parameters at each global iteration. The consensus variable was updated using the PETScTAOSolver implementing a Newton, trust-region method with an absolute tolerance on the gradient of $10^{-12}$, relative gradient tolerance of $10^{-9}$, and at most 10 iterations. 

\begin{figure}[tbh]
    \centering
    \includegraphics[scale = 0.15]{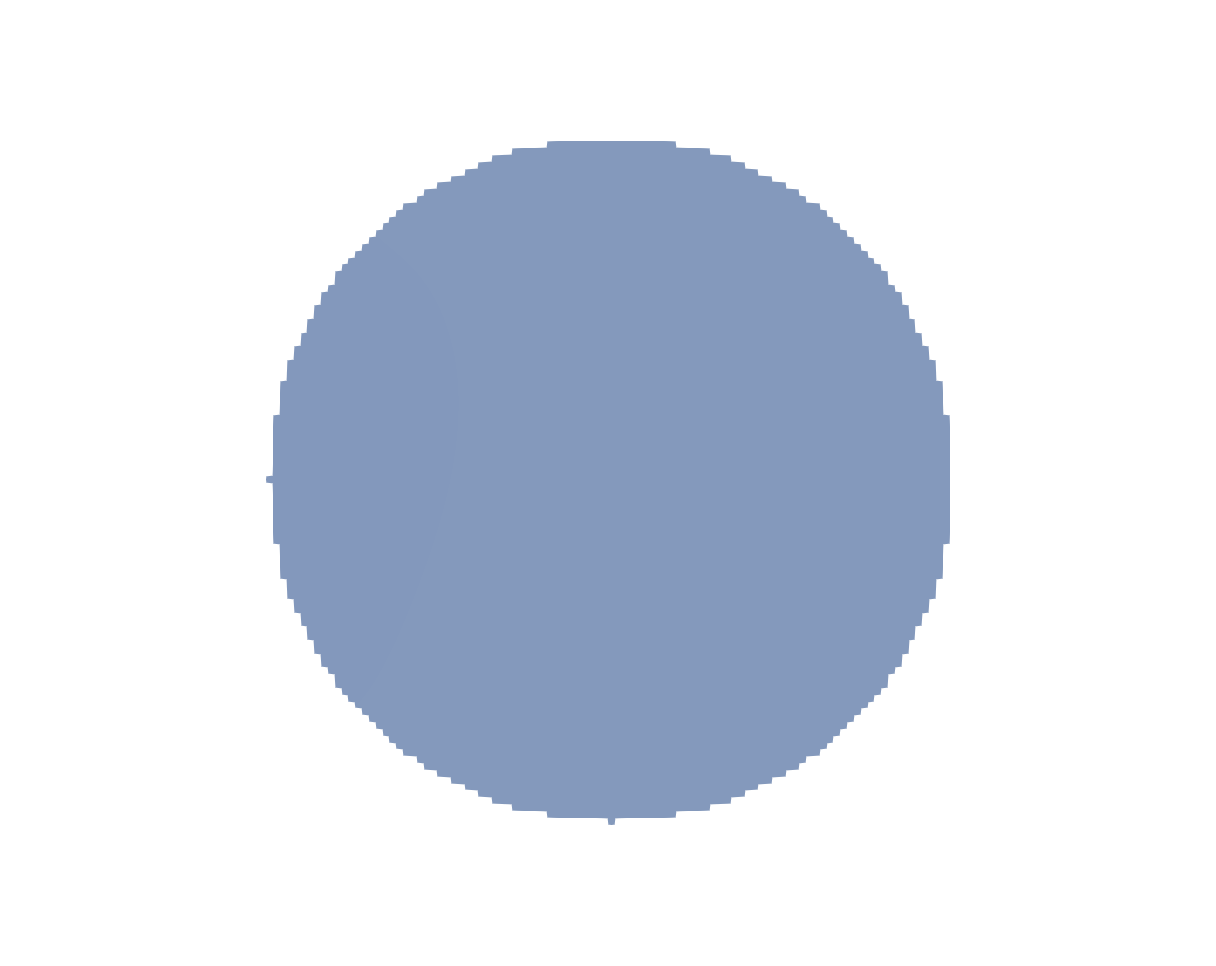}\includegraphics[scale = 0.15]{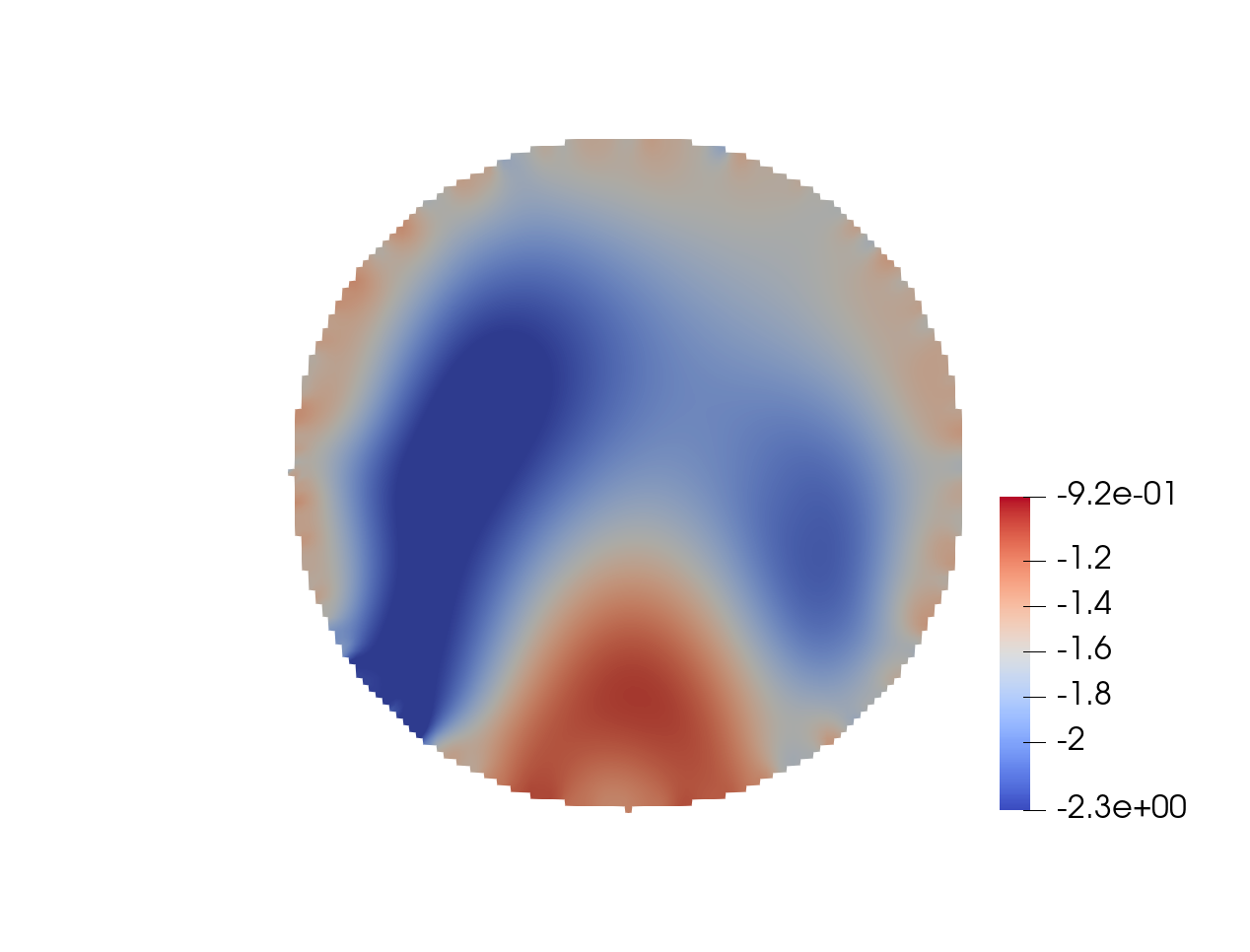}
    \caption{Inverted consensus for EIT problem with ADMM solver using $L^2$ and $H^1$ norms}
    \label{fig:l2_ex_con}
\end{figure}

The ADMM solver with the $L^2$ consensus began with $\rho^0 = 1000$ whereas the ADMM solver with the $H^1$ consensus began with $\rho^0 = 0.1$. The higher value of $\rho^0$ for the $L^2$ consensus is needed to give approximately the same starting consensus weight to both solution methods. The solvers had a global tolerance of $10^{-3}$ and a relative tolerance of $10^{-2}$. The $L^2$ consensus solver terminated in $7$ iterations and the $H^1$ solver terminated in 10 iterations. The final consensus for the $L^2$ solution is pictured on the left of Figure \ref{fig:l2_ex_con} and the final consensus for the $H^1$ solution is on the right.

\begin{figure}[tbh]
    \centering
\resizebox{\textwidth}{!}{\includegraphics[scale = 0.35]{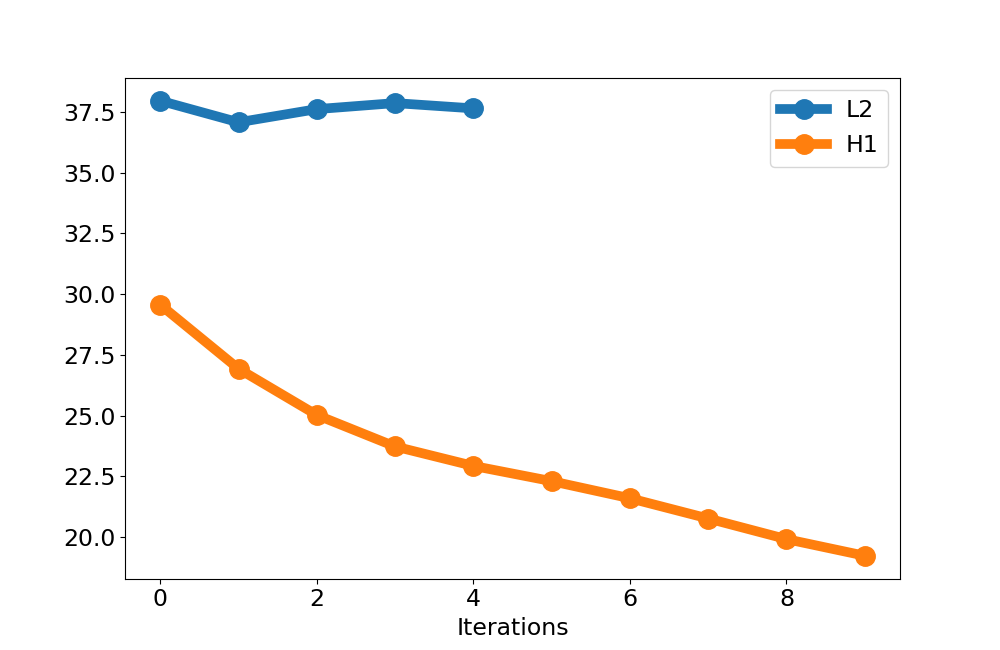}\includegraphics[scale = 0.35]{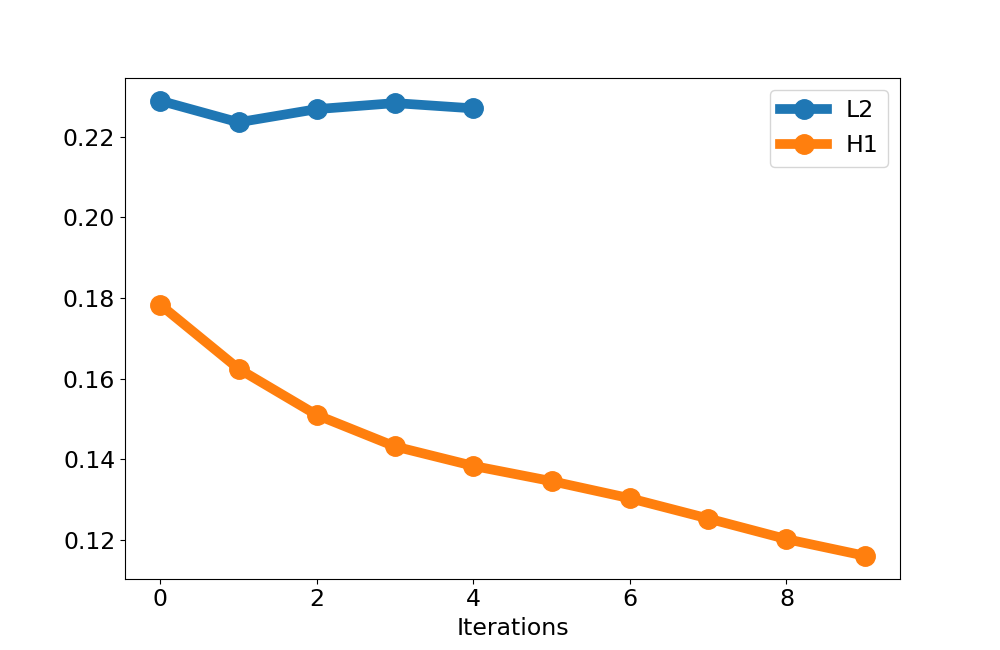}}
    \caption{Error and relative error using $L^2$ and $H^1$ consensus norms}
    \label{fig:L2_H1_error}
\end{figure}

 Qualitatively, we can not discern the original image's traits from the $L^2$ reconstruction, which is approximately constant. Compare this to the $H^1$ reconstruction, which is much more faithful to the ground truth. Quantitatively, this $L^2$ reconstruction had a final relative error of $.2279$ where we define the error in \eqref{eqn:error} and was not reduced across global iterations. The $H^1$ reconstruction had a clear and gradual reduction in error across iterations and had a final relative error of $0.1552$.  Figure \ref{fig:L2_H1_error} displays the error and relative error resulting from the different consensus norms.

 \begin{equation}\label{eqn:error}
\textnormal{Relative error } = \frac{||m_{true} - \iparmap||_{L^2(\Omega)}}{|| m_{true}||_{L^2(\Omega)}}
\end{equation}

\begin{figure}[tbh]
    \centering
    \resizebox{\textwidth}{!}{\includegraphics[scale = 0.35]{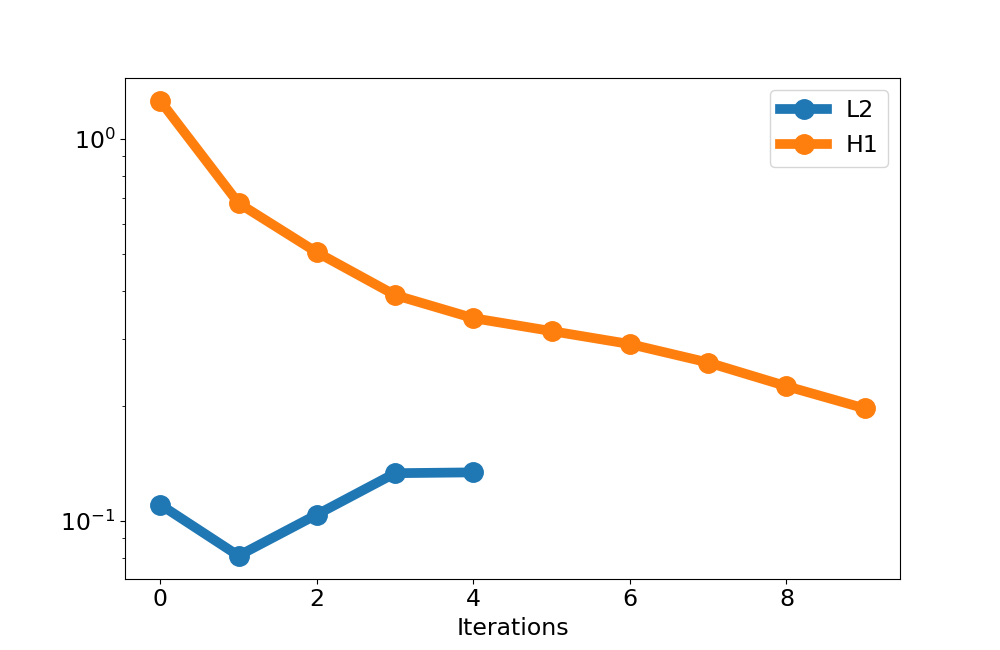}\includegraphics[scale = 0.35]{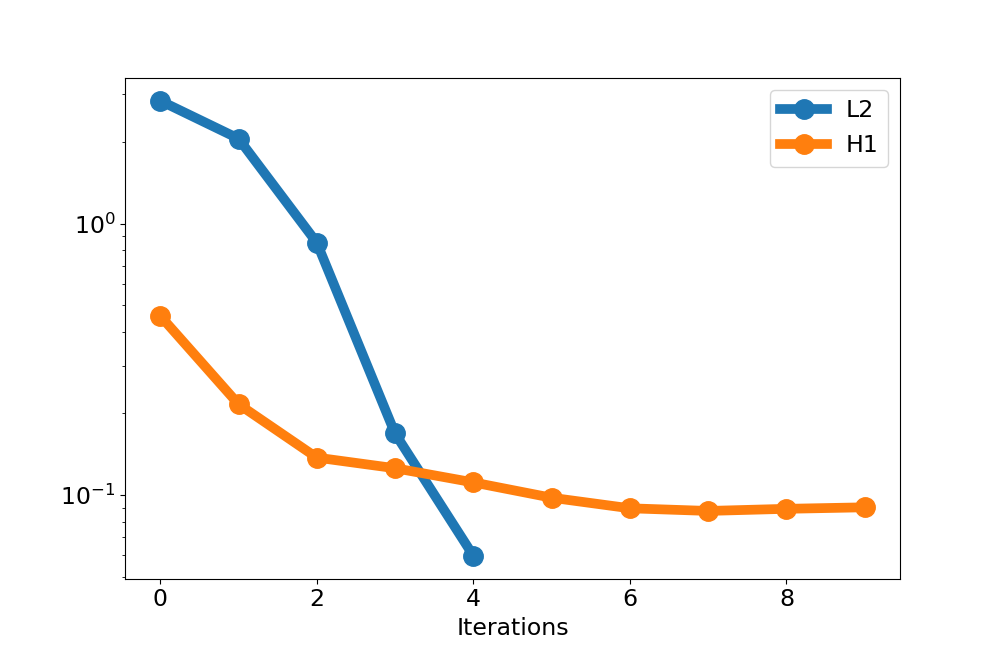}}
    \caption{Primal and dual residuals using $H^1$ and $L^2$ consensus norms}
    \label{fig:L2_H1_residuals}
\end{figure}

We can further analyze the effect of the different consensus norms
by looking at the primal and dual residuals' behavior for ADMM displayed in Figure \ref{fig:L2_H1_residuals}.

Using the $H^1$ norm creates a stable and clear decrease in the residuals compared to the $L^2$ norm. This improvement in the inversion performance and stability demonstrates that the $H^1$ norm is better suited for applying the ADMM to variational problems of this form compared to using the $L^2$ norm.
The Table \ref{Tab:l2_h1_table} summarizes the accuracy of the the $L^2$ and $H^1$ models via the relative error and state misfit. Table \ref{Tab:l2_h1_table} also outlines the computational cost of each method via the solution time and forward, adjoint, and incremental solves.

\begin{table}[tbh]
\caption{Comparison between use of $L^2$ and $H^1$ consensus norm  in ADMM}
\label{Tab:l2_h1_table}
\resizebox{\textwidth}{!}{\begin{tabular}{ l|llllllll} 
 Norm & Iterations & Solution time & Relative Error  & State misfit & Forward solves & Adjoint solves & Incremental Solves \\
 \hline
 $L^2$  &5 & 1m 57s & 0.2249  & 3.249 & 1014 & 799  & 4444 \\
$H^1$  &10& 6m 59s & 0.1160  & 1.228 &  1946 & 1583 & 33968
\end{tabular}}
\end{table}

\subsection{Reconstruction with $H^1$ norm using inexact parameters solutions}

In the previous section, we demonstrated the performance of ADMM by using the $H^1$ norm for consensus reinforcement. We can further improve upon this solution's performance by solving the inverse problems associated with each PDE model inexactly. This will increase the number of global iteration, but each iteration becomes progressively less expensive.

.

\begin{figure}[tbh]
    \centering
    \includegraphics[scale = 0.25]{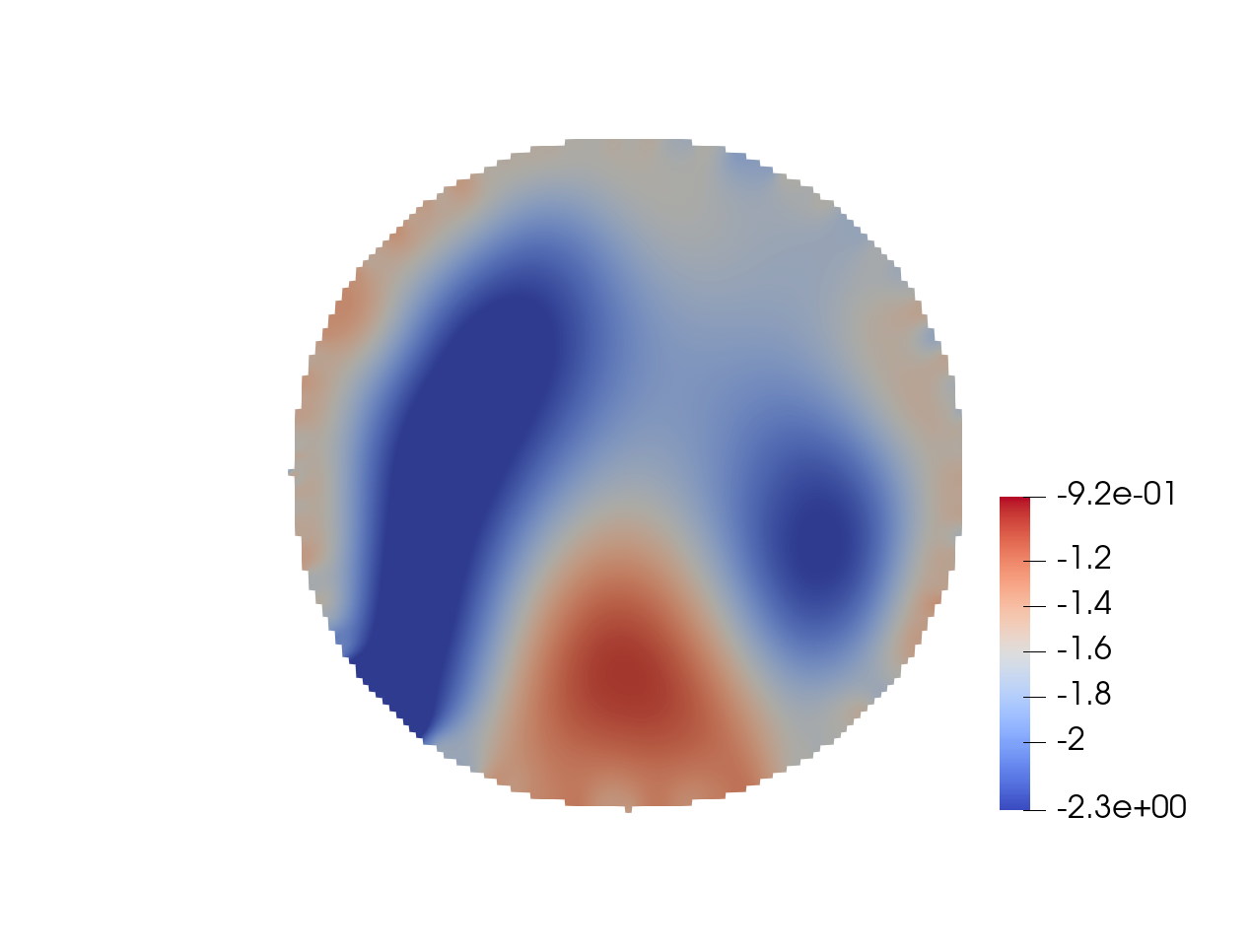}
    \caption{Inverted consensus for EIT problem with ADMM solver using $H^1$ norm}
    \label{fig:in_con}
\end{figure}

With this idea, we proceed with an inversion for $q=16$. The ADMM solver utilized a maximum of 40 global iterations and used an INCG solver with a maximum of 3 iterations for an inexact estimate of the individual parameters at each global iteration. The ADMM solver began with $\rho^0 = 0.1$ and had a global absolute tolerance of $10^{-5}$ and a global relative tolerance of $2\cdot10^{-2}$. 
This solution terminated at 20 iterations when the tolerances were reached and resulted in the final consensus shown in Figure \ref{fig:in_con}. 

\begin{figure}[tbh]
    \centering
    \resizebox{\textwidth}{!}{\includegraphics[scale = 0.35]{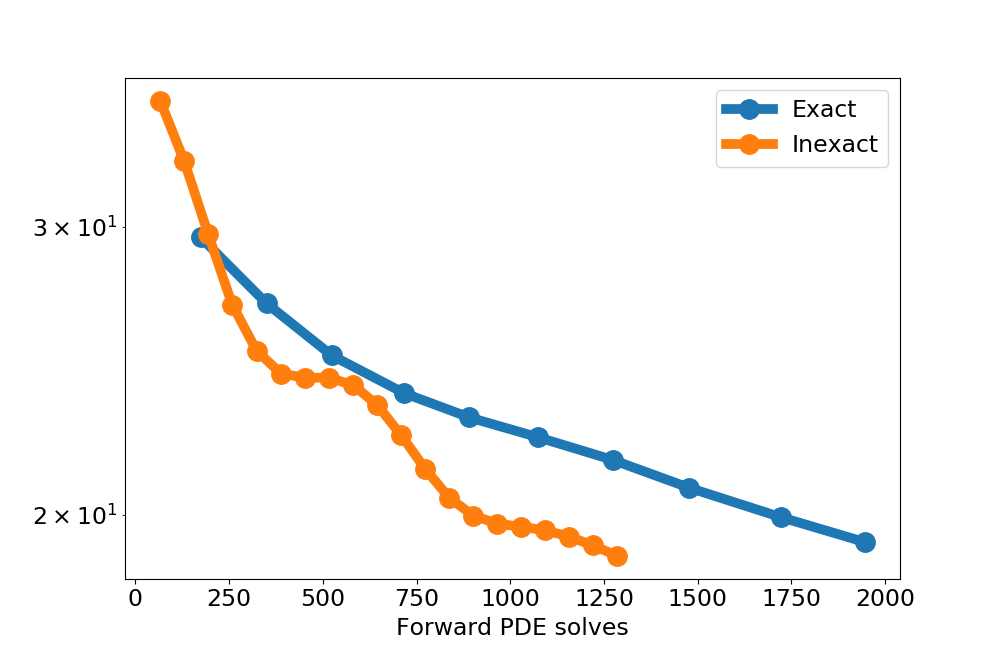}\includegraphics[scale = 0.35]{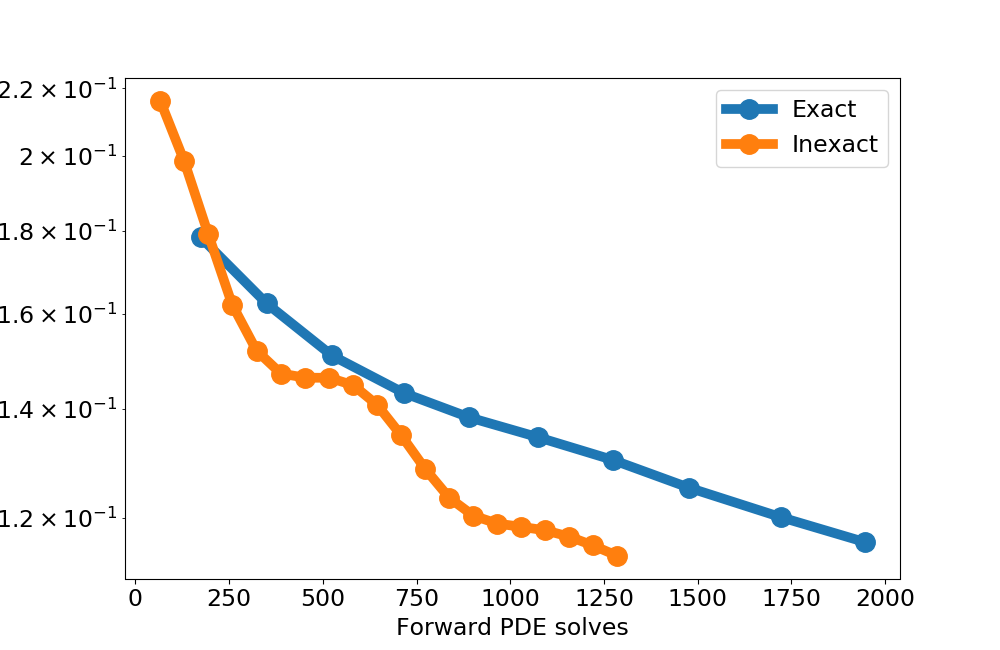}}
    \caption{The error and relative error of the $H^1$ consensus across global iterations}
    \label{fig:in_err}
\end{figure}

\begin{figure}[tbh]
    \centering
    \resizebox{\textwidth}{!}{\includegraphics[scale = 0.35]{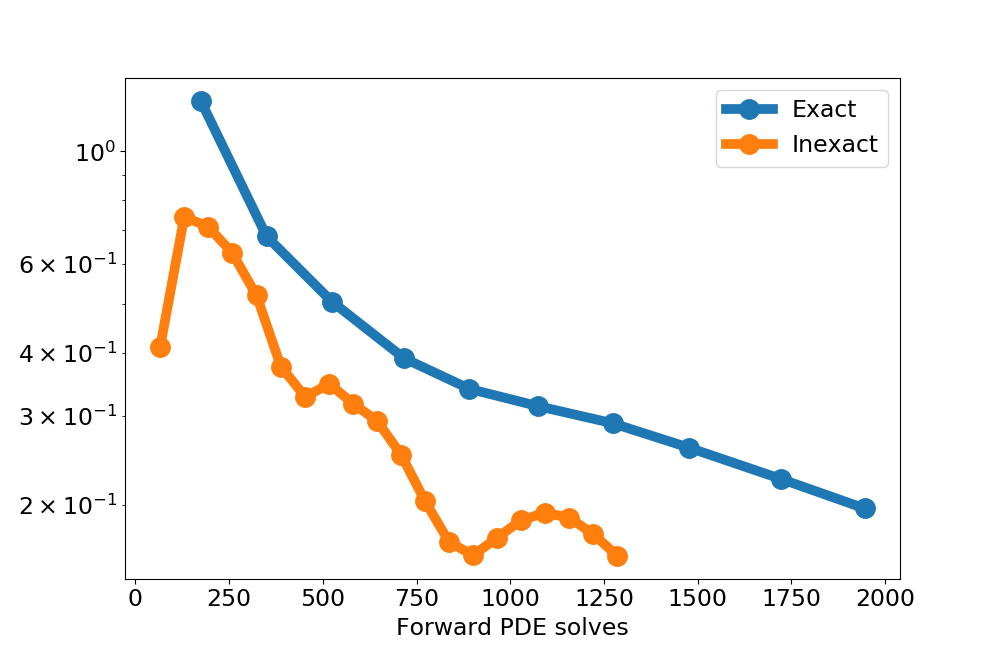}
    \includegraphics[scale = 0.35]{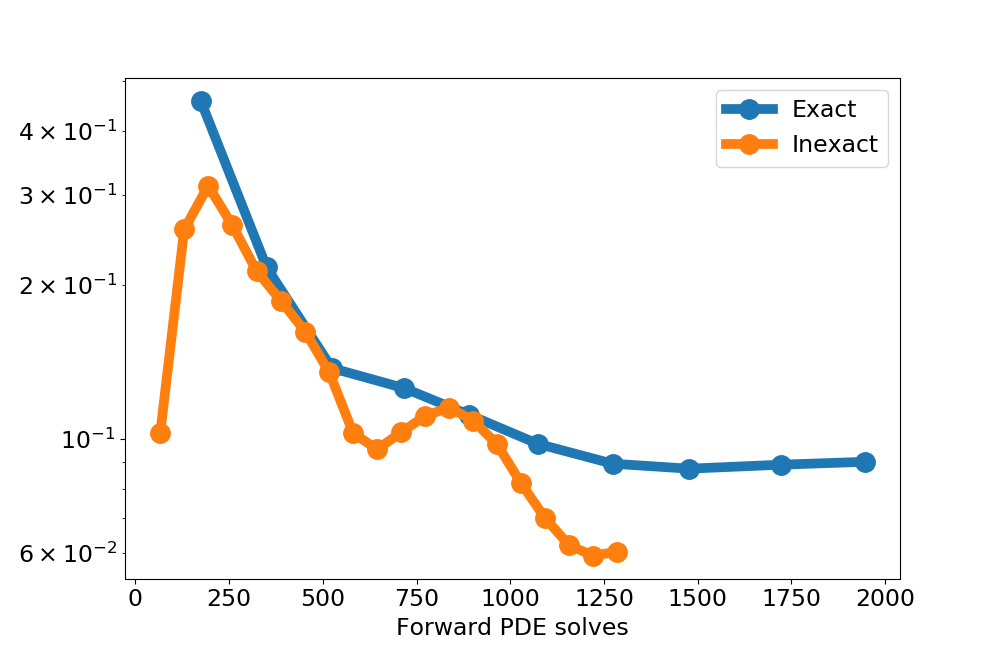}}
    \caption{Primal and dual residuals using $H^1$ consensus}
    \label{fig:in_residuals}
\end{figure}

This reconstruction had even better performance than the instance with exact solutions for each global iteration parameter. The final relative was $0.1160$ and was reduced with each global iteration as seen in Figure \ref{fig:in_err}. Similarly, the primal and dual residuals exhibited very nice behaviors with constant reductions, as seen in Figure \ref{fig:in_residuals}. The plots below have both the Inexact and Exact results plotted against the number of forward PDE solves at that stage.

Utilizing inexact solutions for the parameters also lead to a much faster inversion that met convergence requirements. The exact solution took 2 minutes and 59 seconds and terminated. The inexact solution took 2 minutes and 30 seconds and terminated when the convergence requirements were met. The table \ref{Tab:ex_in_table}
 summarizes the performance of using inexact parameter solver over exact parameter solves.
 
\begin{table}[tbh]
\caption{Comparison between exact and inexact solution of subproblems in ADMM}
\label{Tab:ex_in_table}
\resizebox{\textwidth}{!}{\begin{tabular}{ l|llllllll} 
 Parameter solves &Iterations & Solution time & Relative Error  & State misfit & Forward solves & Adjoint solves & Incremental Solves \\
 \hline
Exact    &10& 6m 59s & 0.1160  & 1.228 &  1946 & 1583 & 33968 \\
Inexact  &20& 3m 2s & 0.1138  & 1.259 &  1285 & 960 & 6616

\end{tabular}}
\end{table}

\subsection{Comparing the ADMM and monolithic scalability with respect to problem size}\label{ref}
We compared the effectiveness of the ADMM with the monolithic approach by performing multiple reconstructions using both approaches on meshes at various levels of refinement. Using a finer mesh meant a larger number of degrees of freedom for both the inversion parameters and state variables. With this, we performed inversions at four different levels of refinement with $8044, 31816, 71280,$ and  $126428$ degrees of freedom. We then fixed the number of PDE models to $q=16$ for these inversions. For both approaches, we also ran their solution in parallel across 8 processes.

For the ADMM method, we utilized a solver with a global absolute tolerance of $10^{-5}$ and a global relative tolerance of $2\cdot 10^{-2}$. For each global iteration, the local inversion parameter was solved using an INCG solver with 3 iterations, a relative tolerance of $10^{-6}$, and an absolute tolerance of $10^{-9}$. This INCG solver utilized a maximum of 100 conjugate gradient evaluations at each iteration. We also used the $H^1$ norm for consensus reinforcement and began the adaptive $\rho$ scheme with $\rho^0 = 0.1$. The consensus variable was updated using the PETScTAOSolver implementing a Newton, trust-region method with an absolute tolerance on the gradient of $10^{-12}$, relative gradient tolerance of $10^{-9}$, and at most 10 iterations. 

We implemented an INCG solver for the monolithic solution with a relative tolerance of $10^{-6}$ and an absolute tolerance of $10^{-2}$. This INCG solver utilized a maximum of 100 conjugate gradient evaluations for each iteration and had a maximum number of 75 iterations.

\subsubsection{Reconstructions over varying levels of mesh refinement}

\begin{figure}[tbh]
    \centering
    \resizebox{\textwidth}{!}{\includegraphics[scale = 0.4]{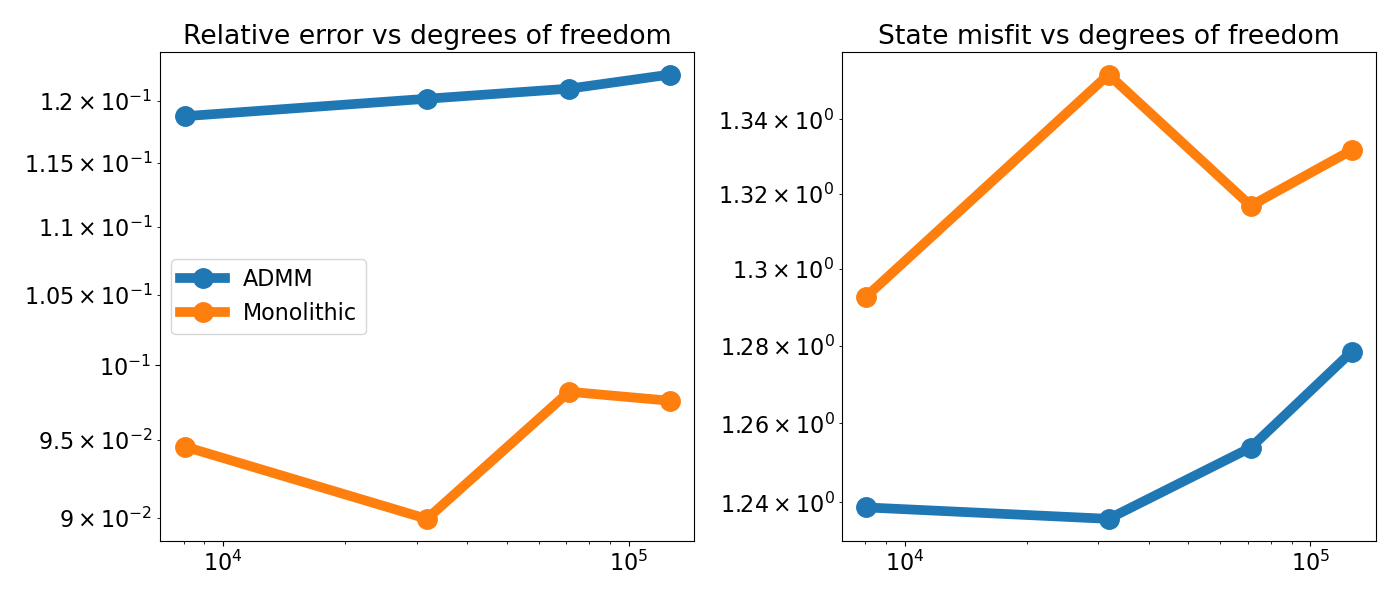}}
    \caption{Relative error and state misfit for ADMM and monolithic approaches vs number of degrees of freedom}
    \label{fig:acc_ref}
\end{figure}

Using the ADMM and monolithic approaches, we achieved the following reconstructions at each level of refinement. The top row contains the ADMM inversions corresponding to 8044, 31816, 71280, and 126428 degrees of freedom. The bottom row contains the monolithic inversions corresponding to 8044, 31816, 71280, and 126428 degrees of freedom. The relative error and state misfits are plotted against the number of degrees of freedom in Figure \ref{fig:acc_ref}

\begin{figure}[tbh]
    \centering
    \resizebox{\textwidth}{!}{\includegraphics[scale = 0.1, trim={6cm 5cm 9cm 5cm},clip, valign=t]{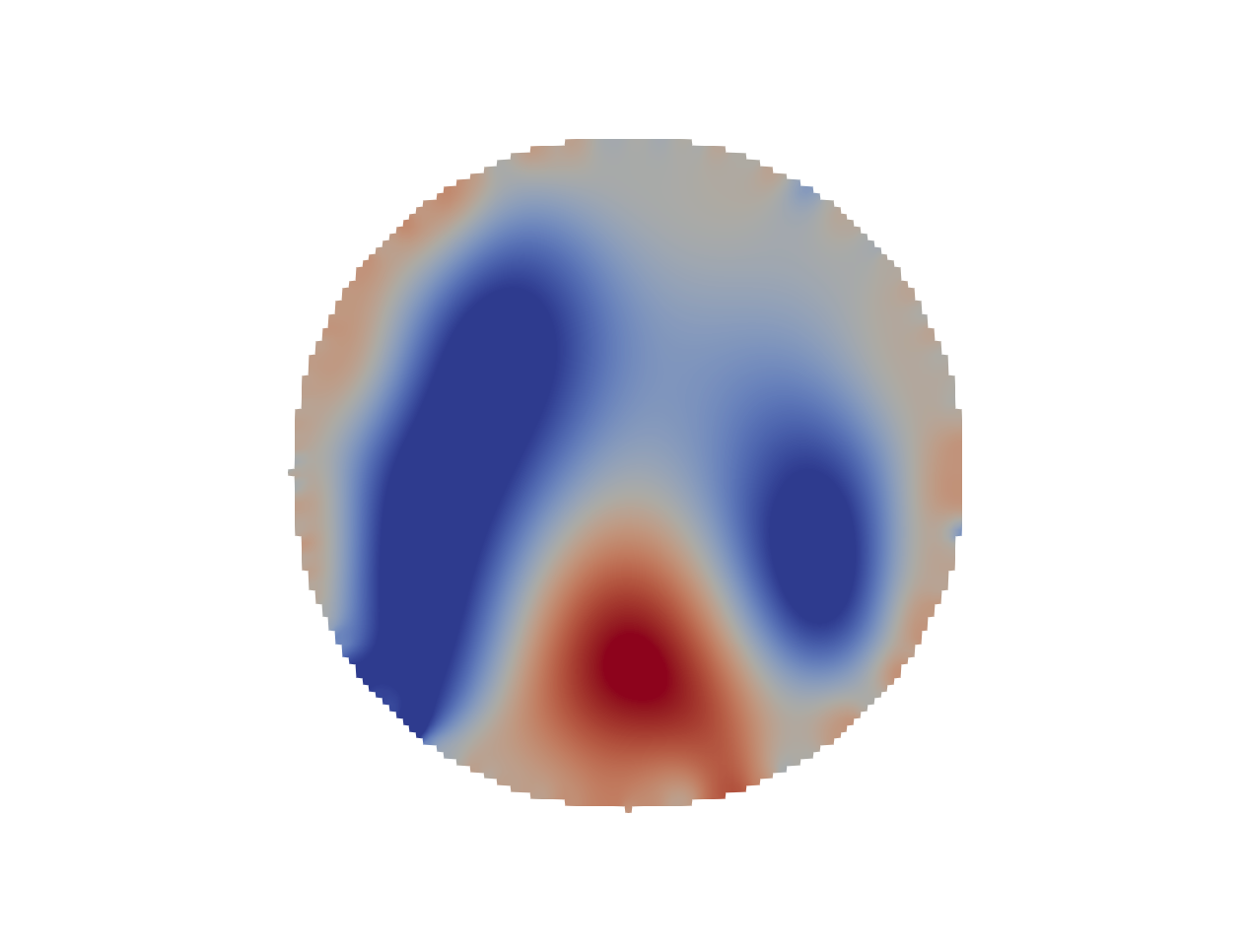}\includegraphics[scale = 0.1, trim={6cm 5cm 9cm 5cm},clip, valign=t]{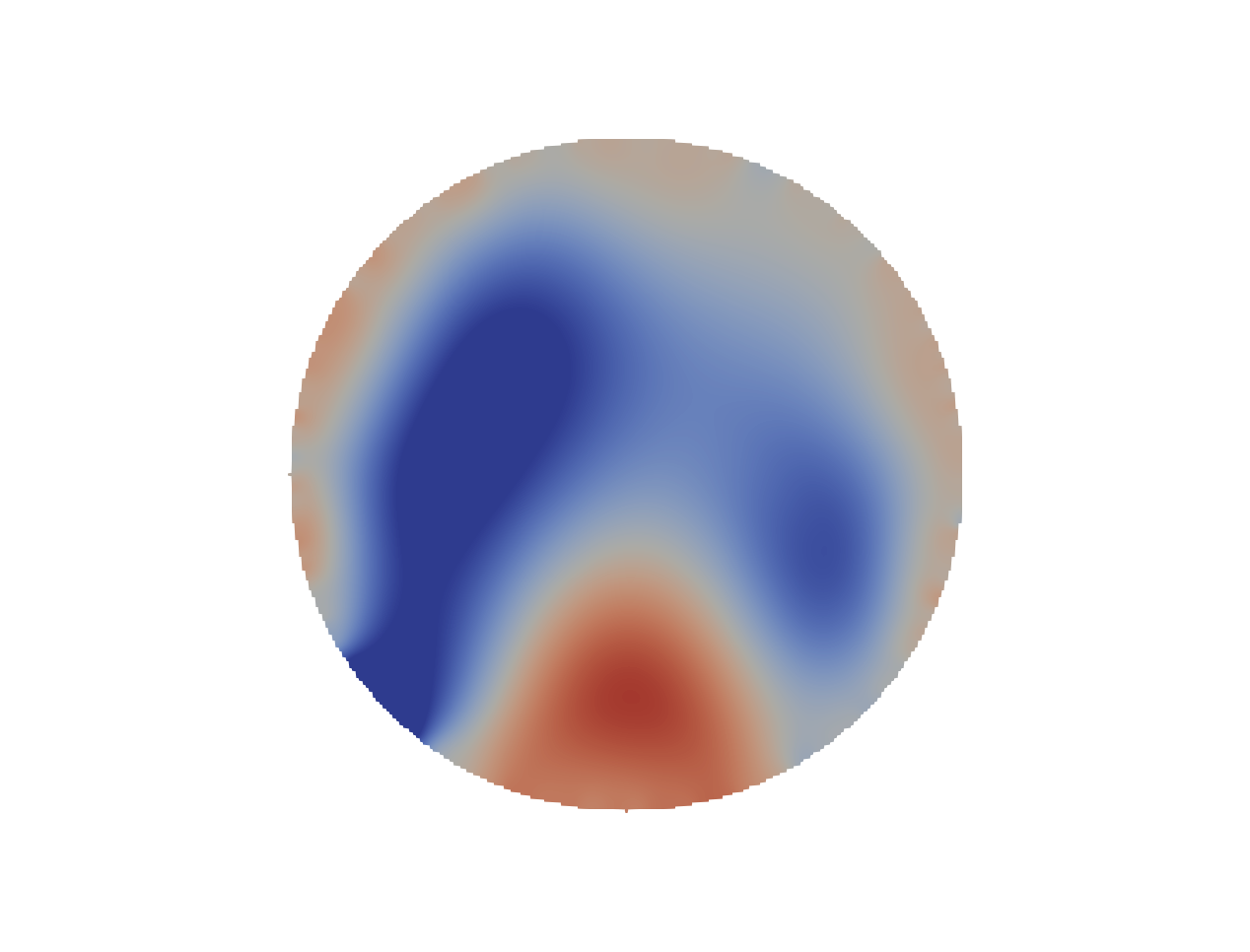}\includegraphics[scale = 0.1, trim={6cm 5cm 9cm 5cm},clip, valign=t]{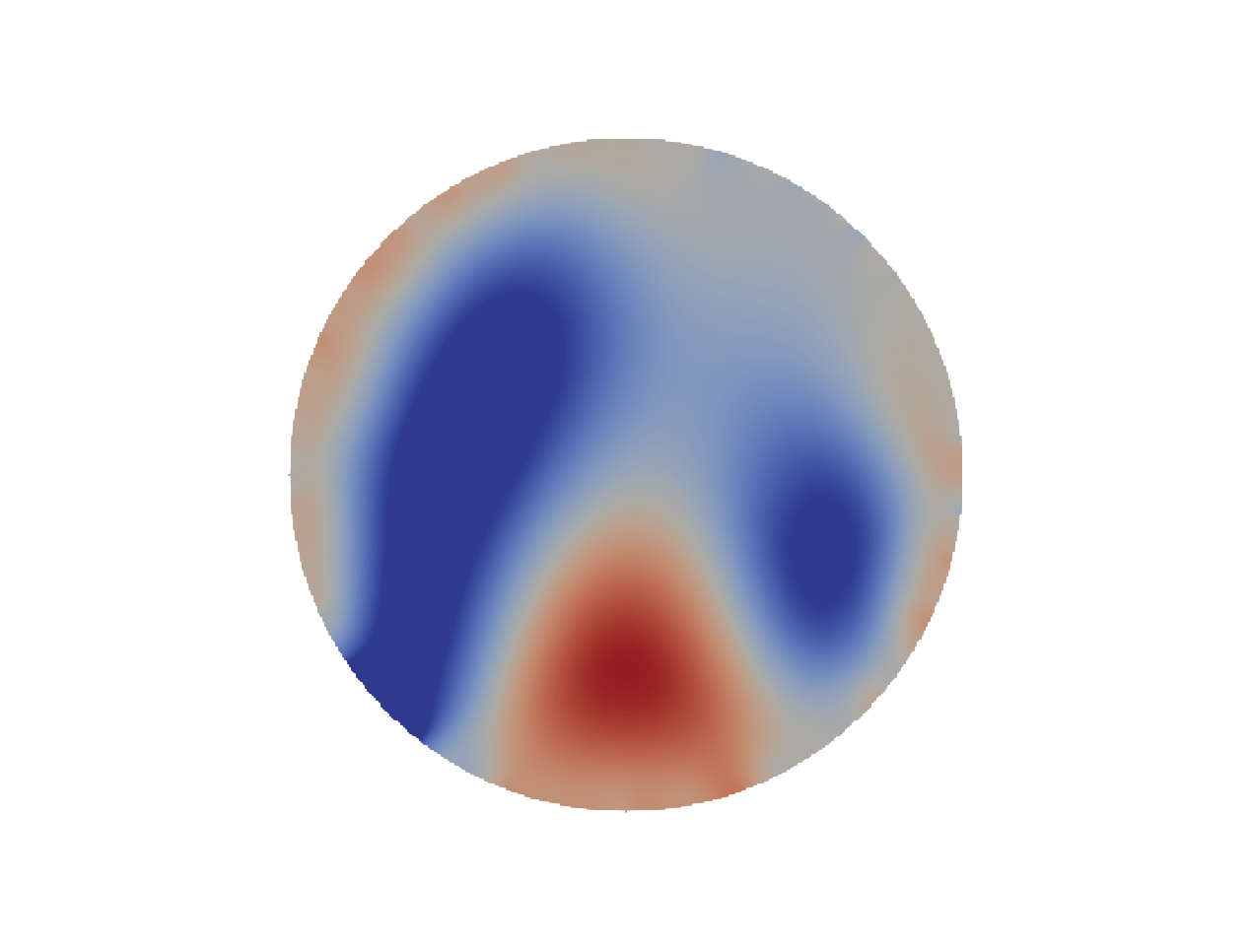}\includegraphics[scale = 0.1, trim={6cm 5cm 9cm 5cm},clip, valign=t]{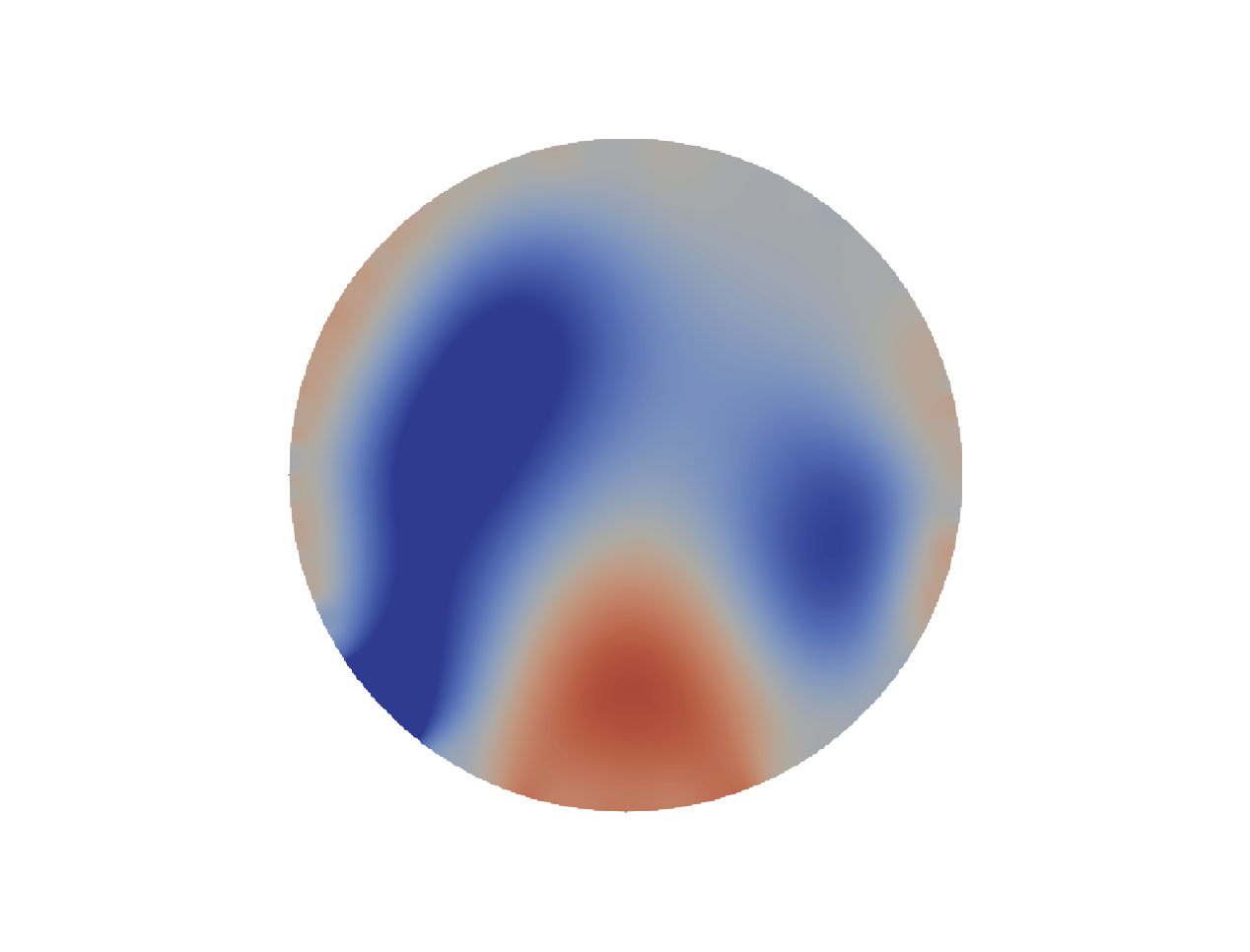}\includegraphics[scale = 0.20, trim={35cm 0 0 15cm},clip, valign=t]{eit_figures/ref4.png}}
    
    \resizebox{\textwidth}{!}{\includegraphics[scale = 0.1, trim={6cm 5cm 9cm 5cm},clip, valign=t]{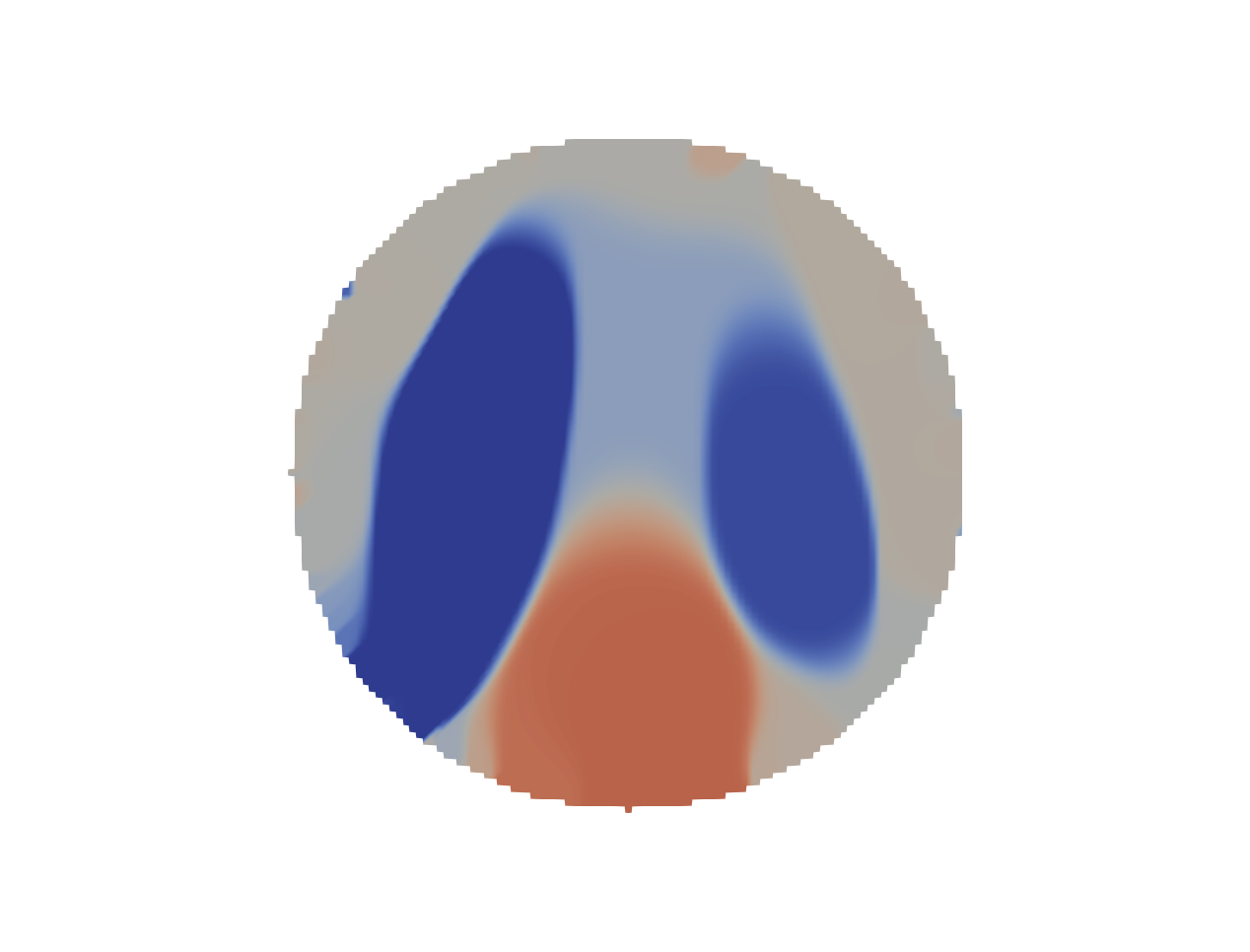}\includegraphics[scale = 0.1, trim={6cm 5cm 9cm 5cm},clip, valign=t]{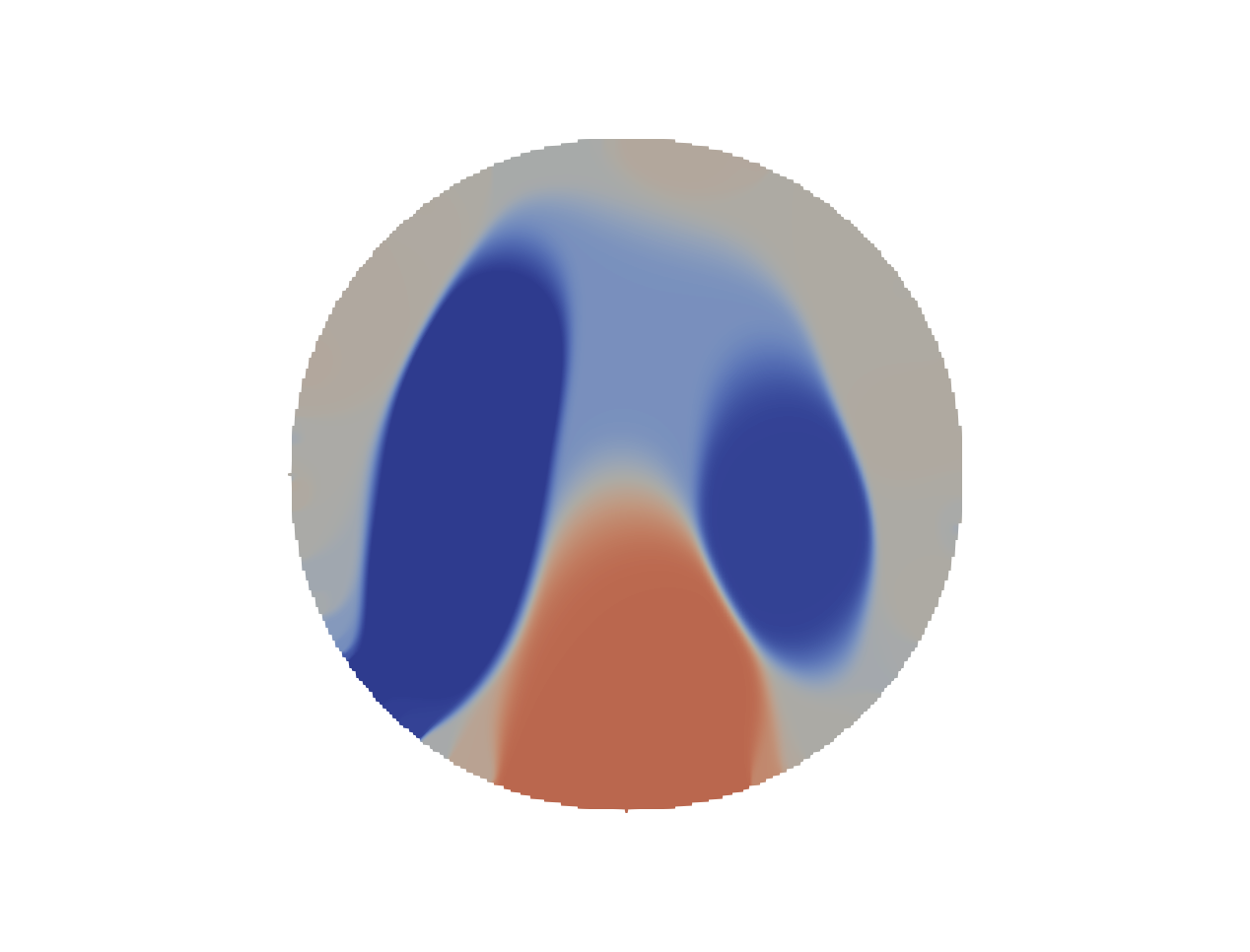}\includegraphics[scale = 0.1, trim={6cm 5cm 9cm 5cm},clip, valign=t]{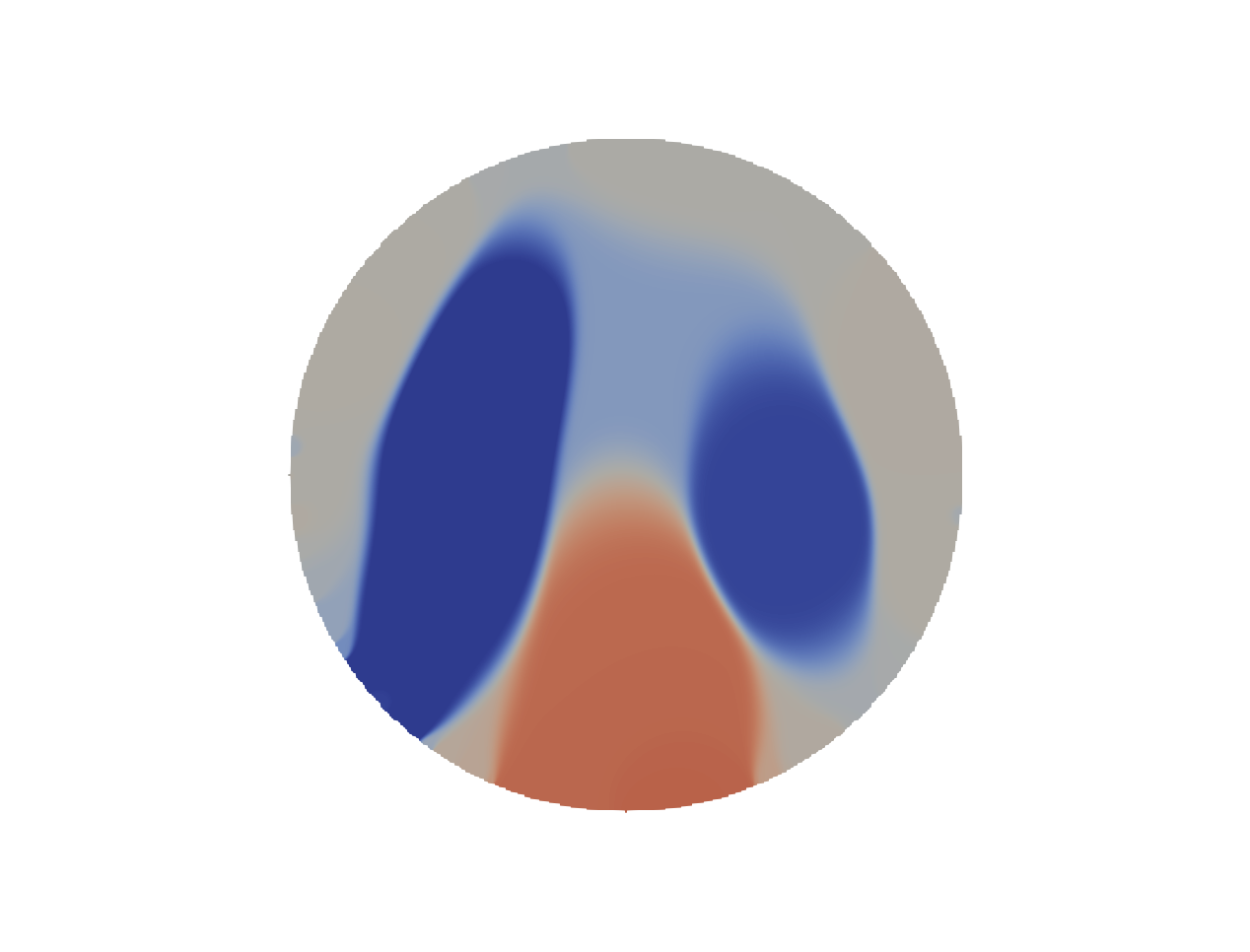}\includegraphics[scale = 0.1, trim={6cm 0 9cm 5cm},clip, valign=t]{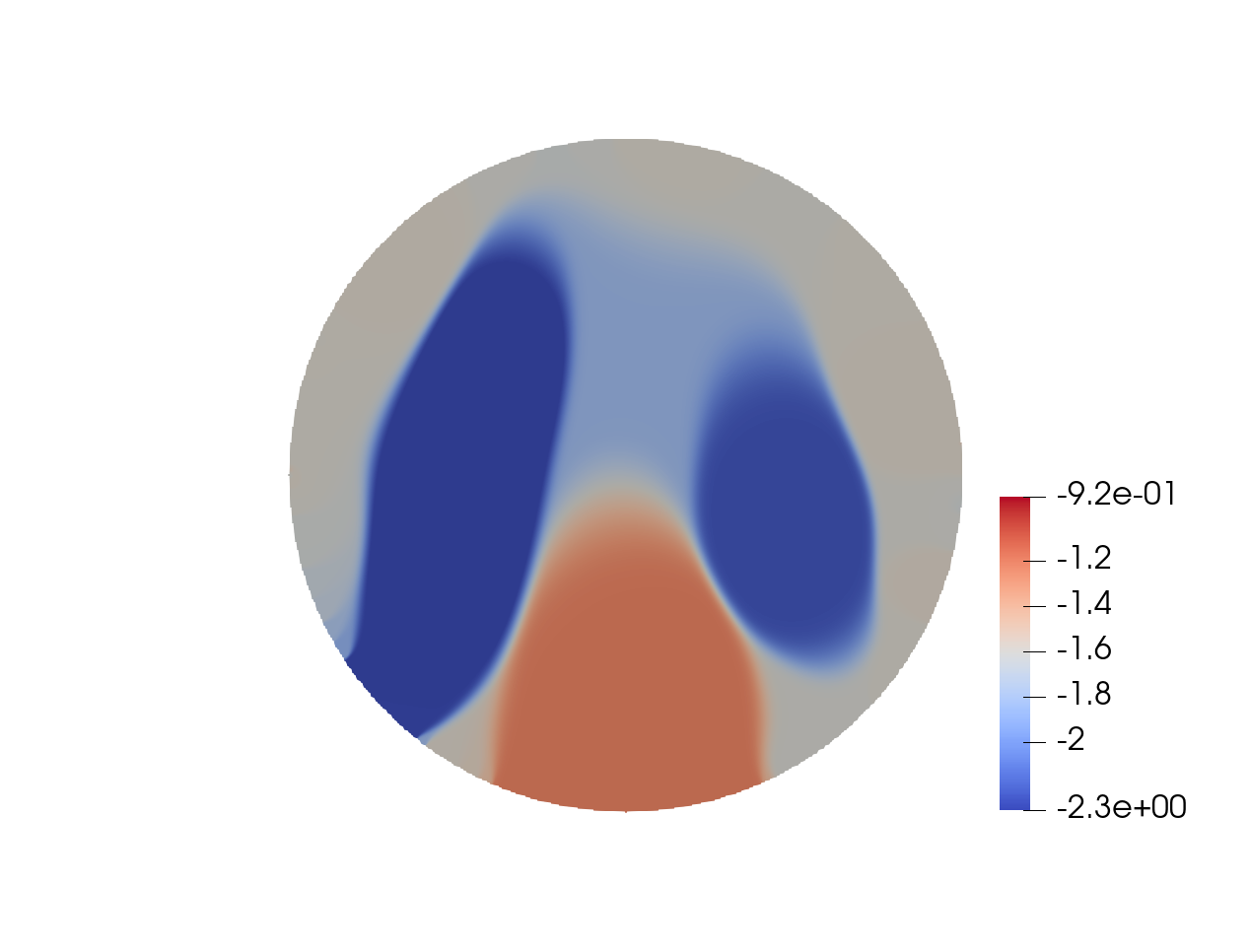}\includegraphics[scale = 0.20, trim={35cm 0 0 15cm},clip, valign=t]{eit_figures/mon_ref4.png}}
    \caption{ADMM and monolithic reconstructions over multiple mesh refinements}
    \label{fig:ref}
\end{figure}

\subsubsection{Computational scalability with respect problem size}

\begin{figure}[tbh]
    \centering
    \includegraphics[scale = 0.3]{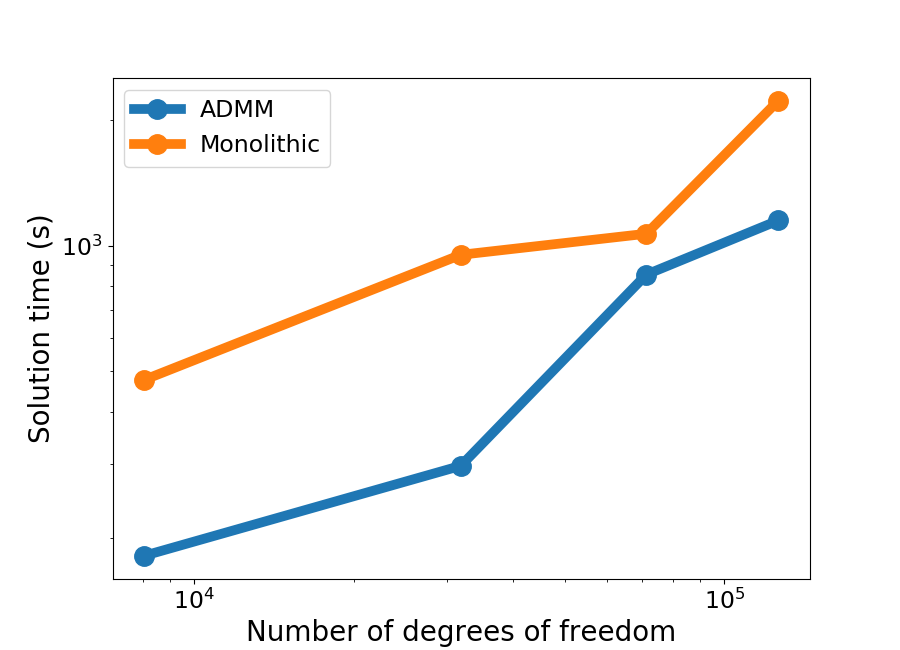}
    \caption{Solution time(s) vs number of degrees of freedom}
    \label{fig:ref_time}
\end{figure}

\begin{figure}[tbh]
    \centering
    \resizebox{\textwidth}{!}{\includegraphics[scale = 0.4, trim={3cm 0 3.5cm 0},clip]{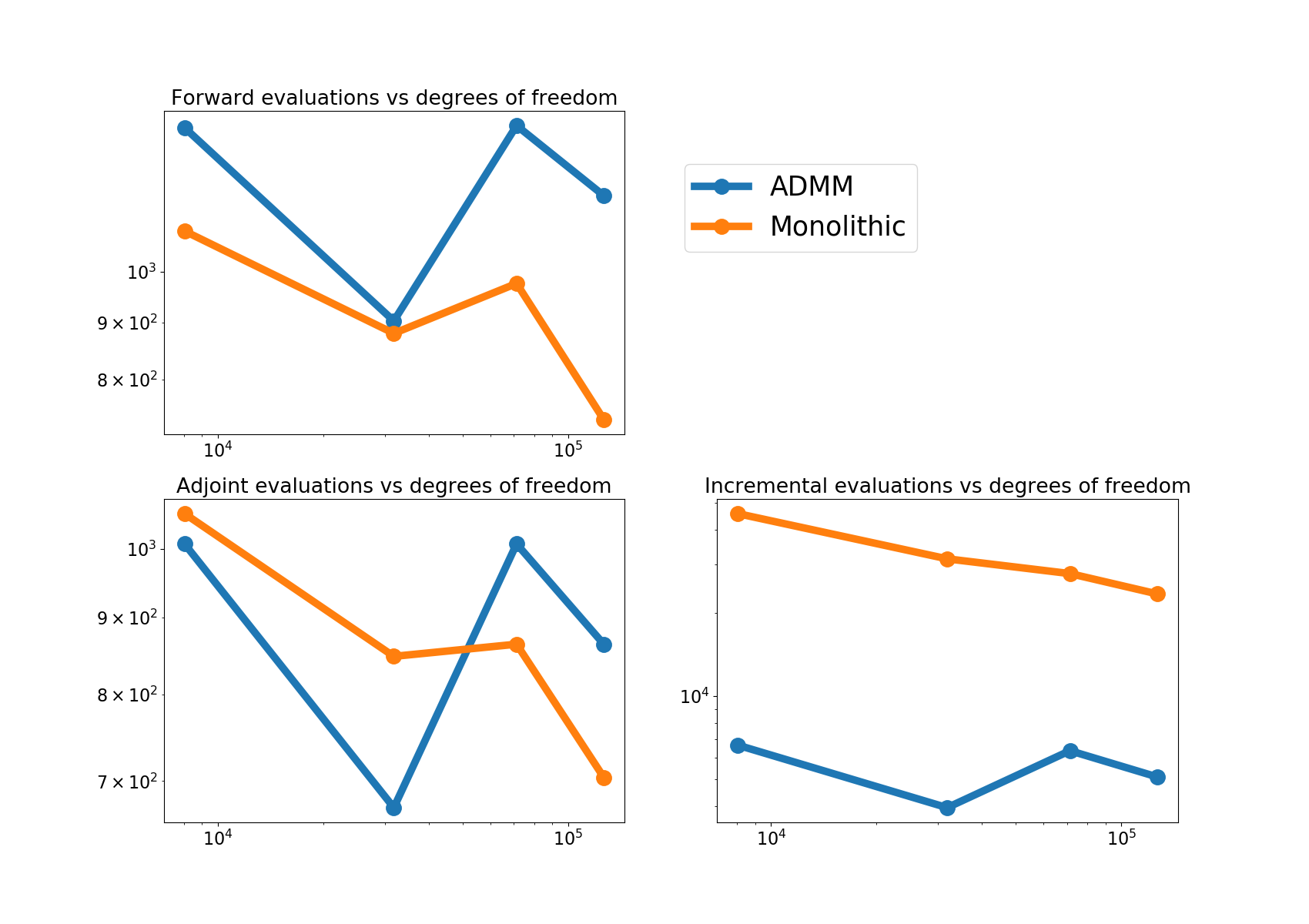}}
    \caption{Number of forward solves, adjoint solves, and incremental solves for ADMM and monolithic approaches vs number of degrees of freedom}
    \label{fig:comp_ref}
\end{figure}

 The time required for each reconstruction is plotted against the number of degrees of freedom in Figure \ref{fig:ref_time}. The computational cost was assessed in three different ways: tracking how many times the forward problem was solved, how many times the adjoint problem was solved, and how many times the incremental problems were solved. The total number of forward solves, adjoint solves, and incremental solves is plotted against the number of degrees of freedom in Figure \ref{fig:comp_ref}. 

The ADMM approach achieved satisfactory accuracy compared to the monolithic approach at varying problem sizes with similar amounts of forward and adjoint solves.   However, the ADMM solutions required much fewer incremental evaluations. Having fewer incremental evaluations then reduced the amount of time for the solution time. These results are summarized in Table \ref{Tab:ref_table}.

\begin{table}[tbh]
\caption{Execution time comparison between ADMM and monolithic approach as a function of the number of inversion parameters}
\label{Tab:ref_table}

\resizebox{\textwidth}{!}{\begin{tabular}{ l|llll } 
 \textnormal{Degrees of freedom}& 8044 & 31816  & 71280 &  126428 \\
 \hline
 \textnormal{Solution time} & 3m 1s/7m 32s & 4m 57s/15m 50s &14m 12s/17m 47s & 19m 8s /37m 0s\\ 
 
 \textnormal{Relative error} & 0.1104/ 0.0864 & 0.1206/0.0954 & 0.1194/0.0952& 0.1222/0.0985\\ 
 \textnormal{State misfit} & 1.259/1.249 & 1.277/1.342 &1.294/1.317 & 1.37/1.314\\ 
 \textnormal{Forward solves} & 1348/1088 & 903/880 &1354/976 & 1171/736\\
 \textnormal{Adjoint solves} & 1008/1065 & 672/848 & 1008/864& 864/704\\ 
 \textnormal{Incremental solves} & 6649/45616 & 3955/31328  & 6351/27712& 5108/23408\\

\end{tabular}}
\end{table}

\subsection{Comparing the ADMM and monolithic scalability with respect to number of forward models}\label{prob}
We compared the effectiveness of the ADMM with the monolithic approach by performing multiple reconstructions using both approaches with a varying number of PDE models on a mesh with a fixed number of degrees of freedom. The number of PDE models varied according to $q=8,16,32,64$.  The solutions for both approaches were found in parallel across 8 processes. The solver settings were the same as in subsection \ref{ref}.

\subsubsection{Reconstructions over with multiple PDE models}

\begin{figure}[tbh]
    \centering
    \resizebox{\textwidth}{!}{\includegraphics[scale = 0.4]{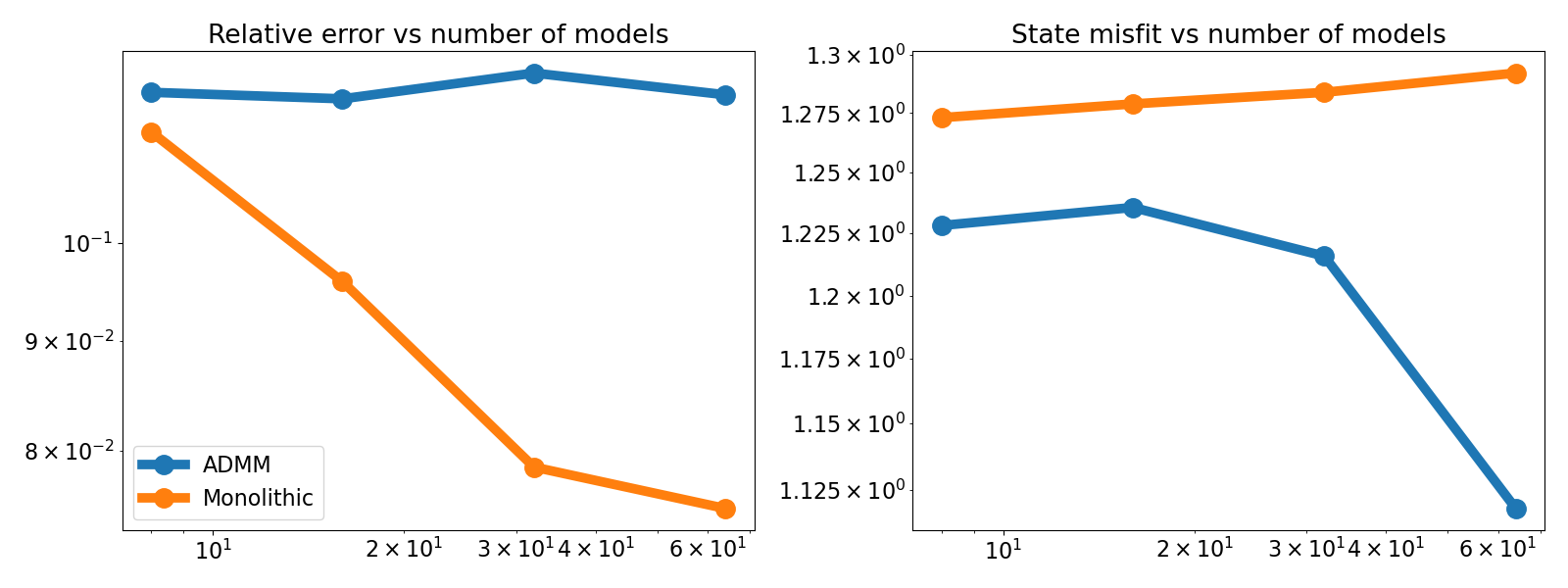}}
    \caption{Relative error and state misfit for ADMM and monolithic approaches vs number of PDE models}
    \label{fig:acc_prob}
\end{figure}

Using the \newline ADMM and monolithic approaches, we achieved the following reconstructions at each level of refinement. The top row contains the ADMM inversions corresponding to 8, 16, 32, and 64 PDE models. The bottom row contains the monolithic inversions corresponding to 8, 16, 32, and 64 PDE models.
The relative error and state misfits are plotted against the number of PDE models in Figure \ref{fig:acc_prob}

\begin{figure}[tbh]
    \resizebox{\textwidth}{!}{\includegraphics[scale = 0.1, trim={6cm 5cm 9cm 5cm},clip, valign=t]{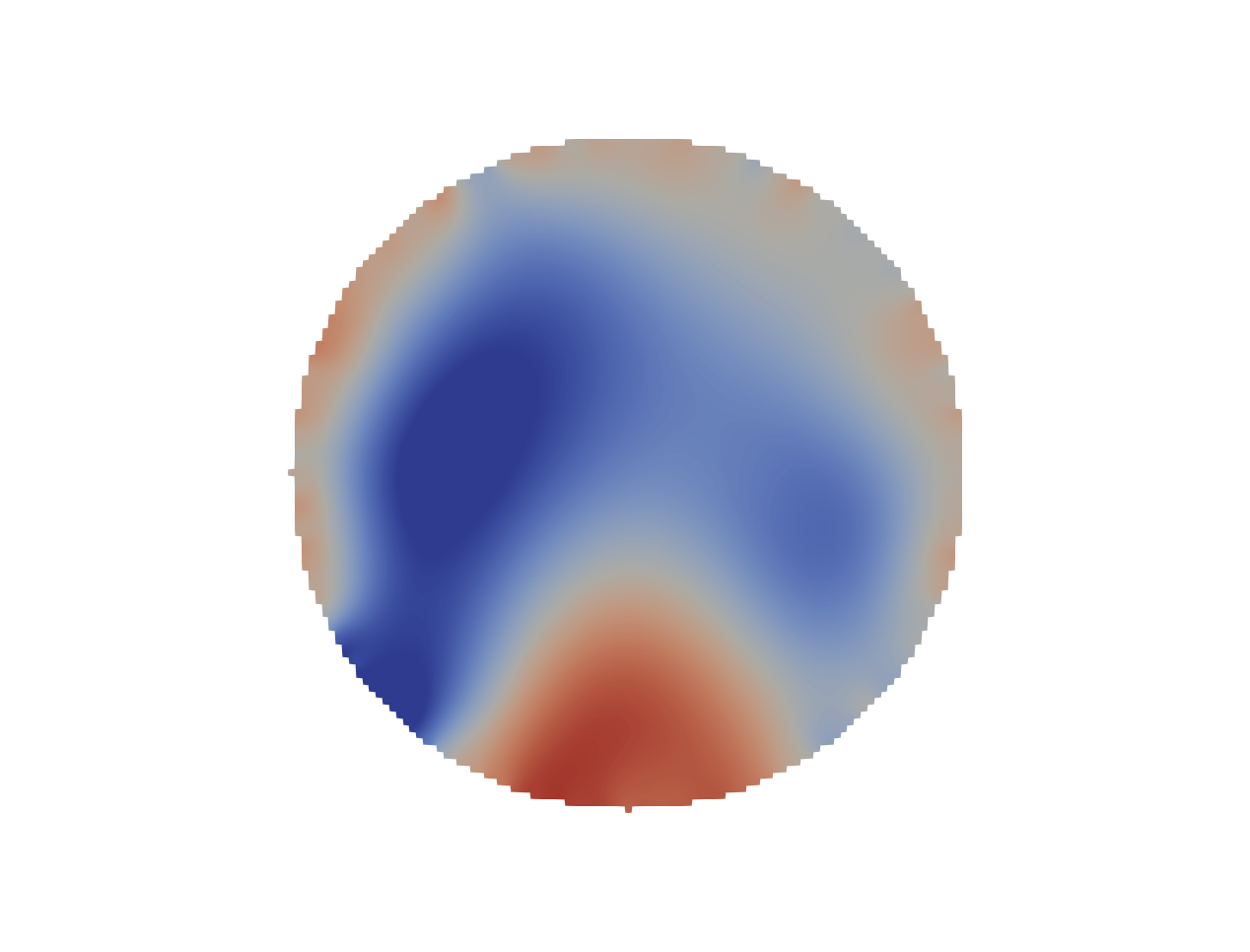}\includegraphics[scale = 0.1, trim={6cm 5cm 9cm 5cm},clip, valign=t]{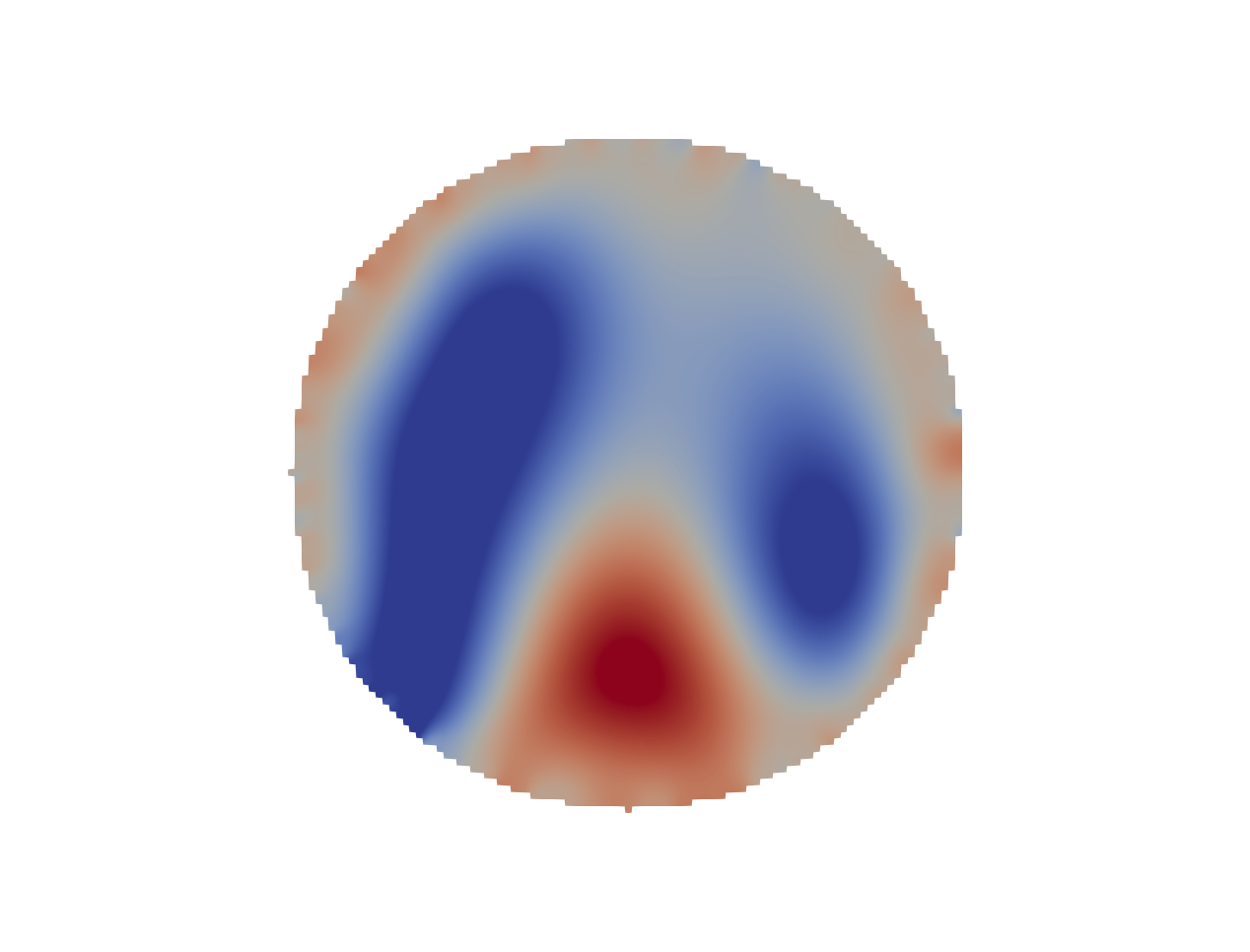}\includegraphics[scale = 0.1, trim={6cm 5cm 9cm 5cm},clip, valign=t]{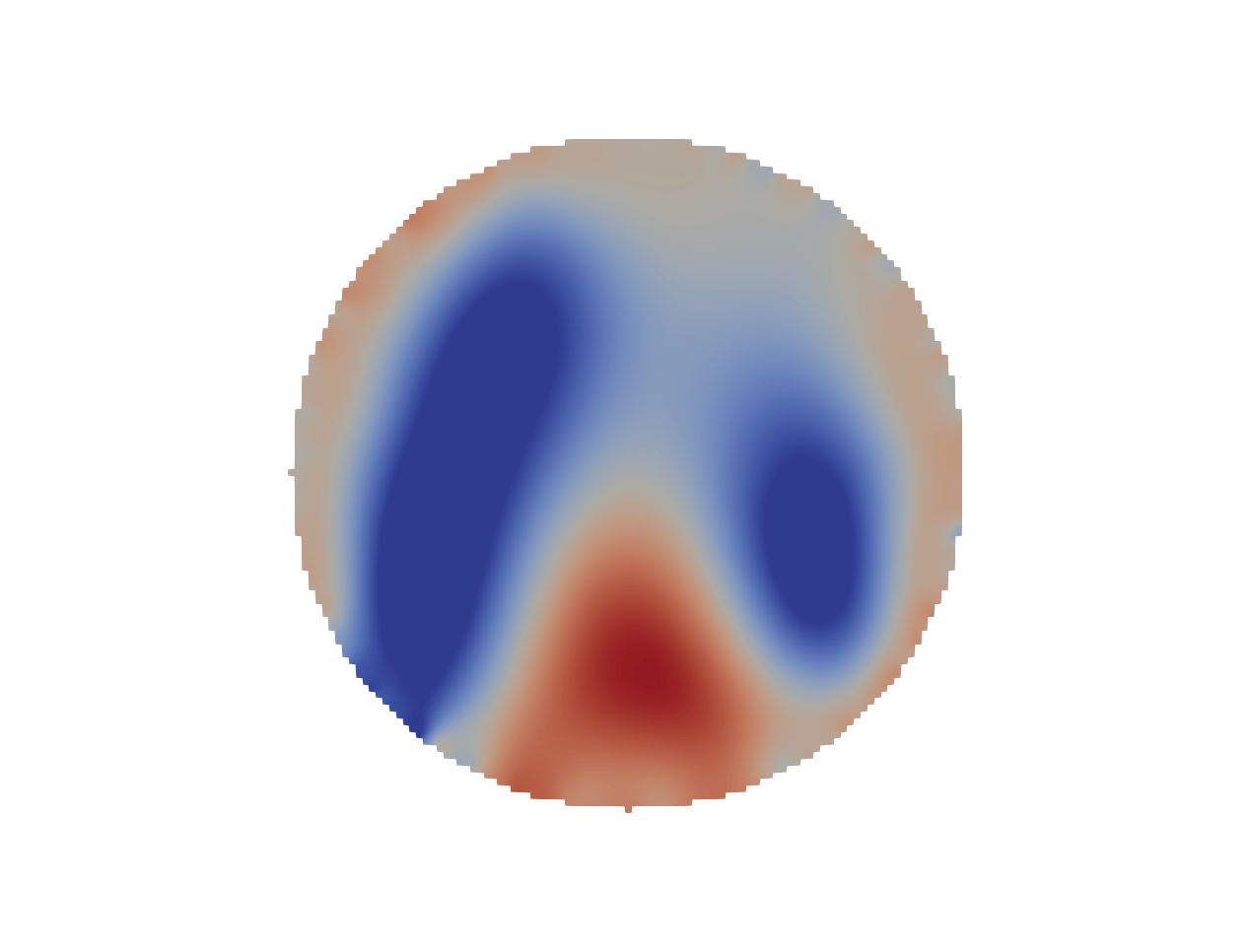}\includegraphics[scale = 0.1, trim={6cm 0 9cm 5cm},clip, valign=t]{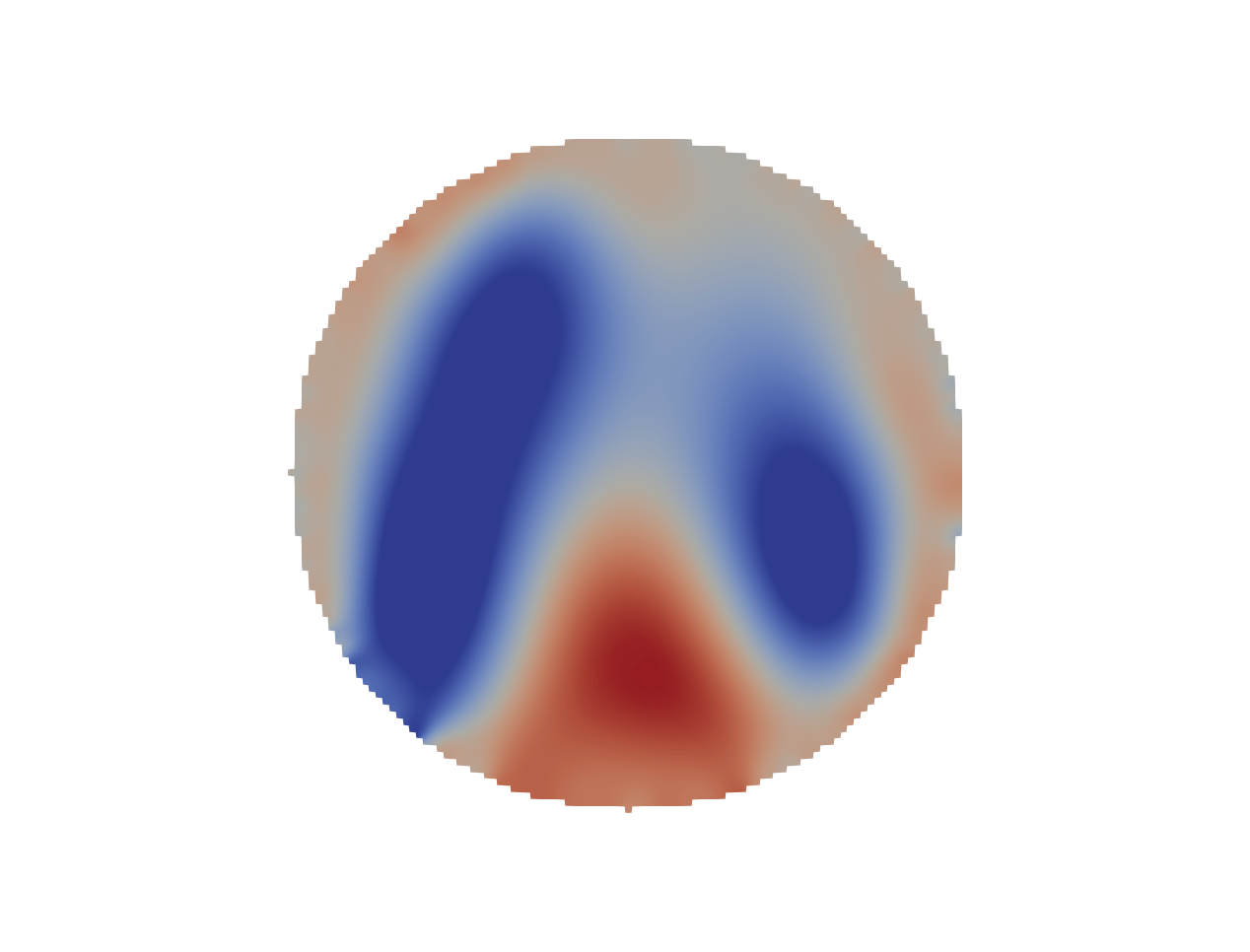}\includegraphics[scale = 0.20, trim={35cm 0 0 15cm},clip, valign=t]{eit_figures/ref4.png}}
    
    \resizebox{\textwidth}{!}{\includegraphics[scale = 0.1, trim={6cm 5cm 9cm 5cm},clip, valign=t]{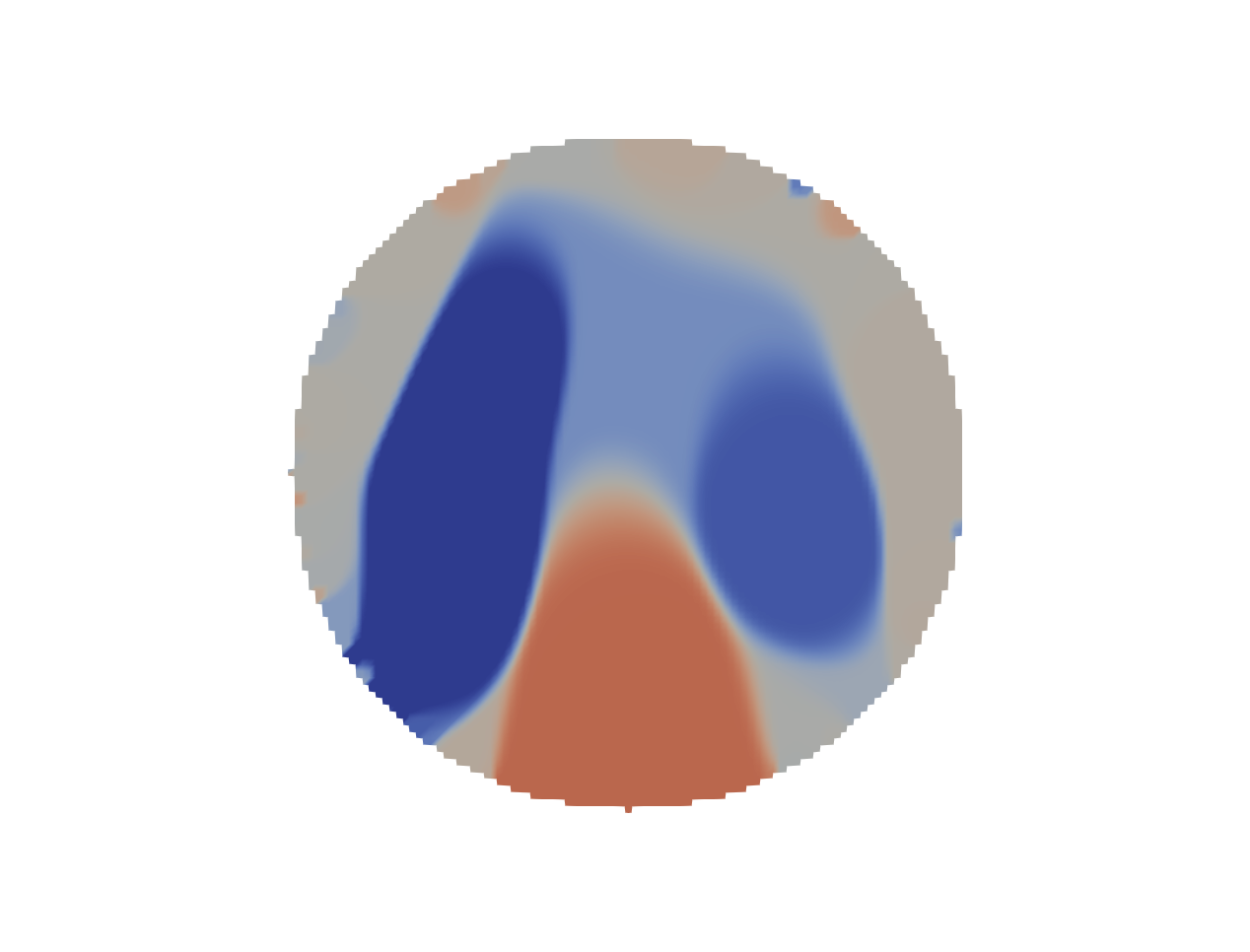}\includegraphics[scale = 0.1, trim={6cm 5cm 9cm 5cm},clip, valign=t]{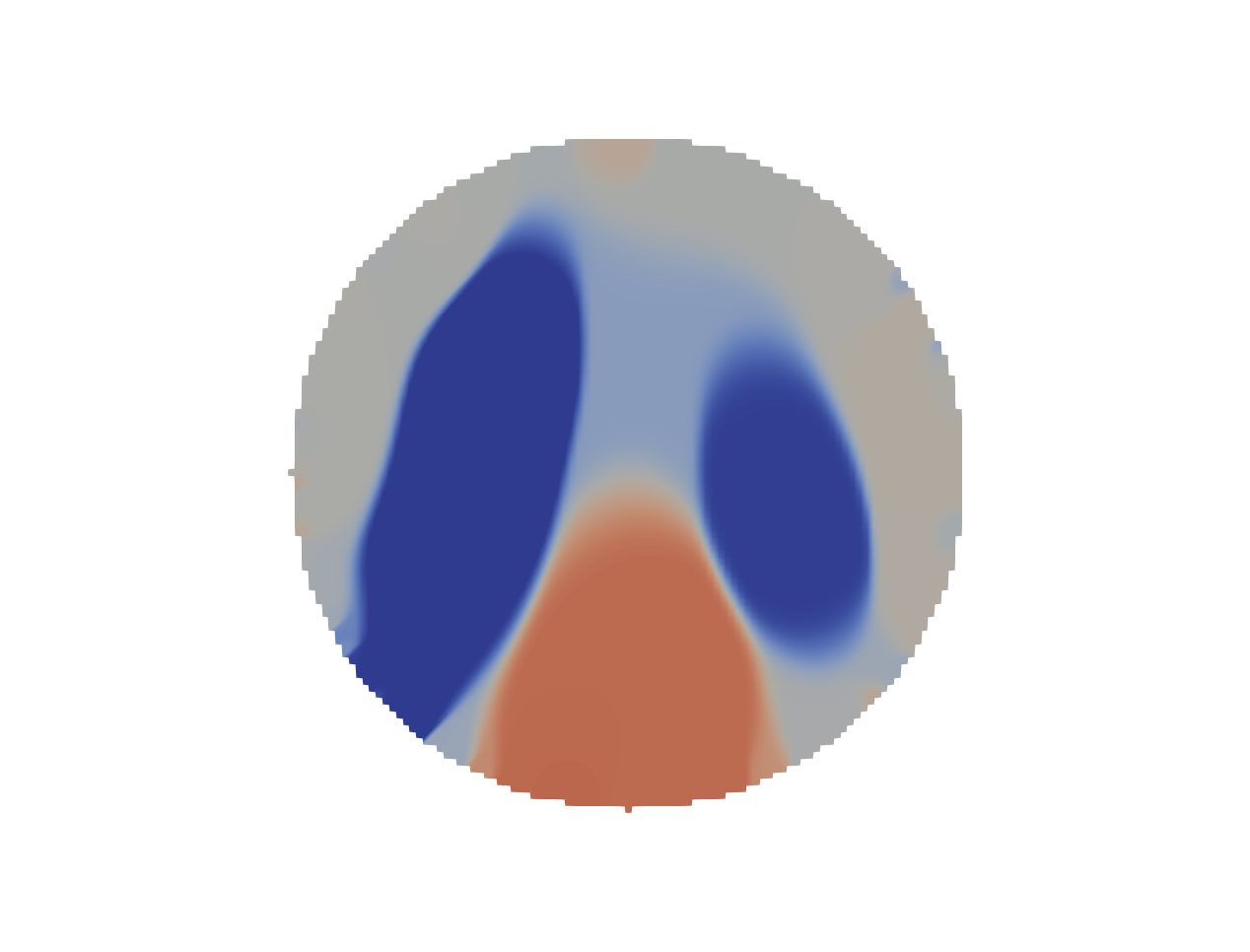}\includegraphics[scale = 0.1, trim={6cm 5cm 9cm 5cm},clip, valign=t]{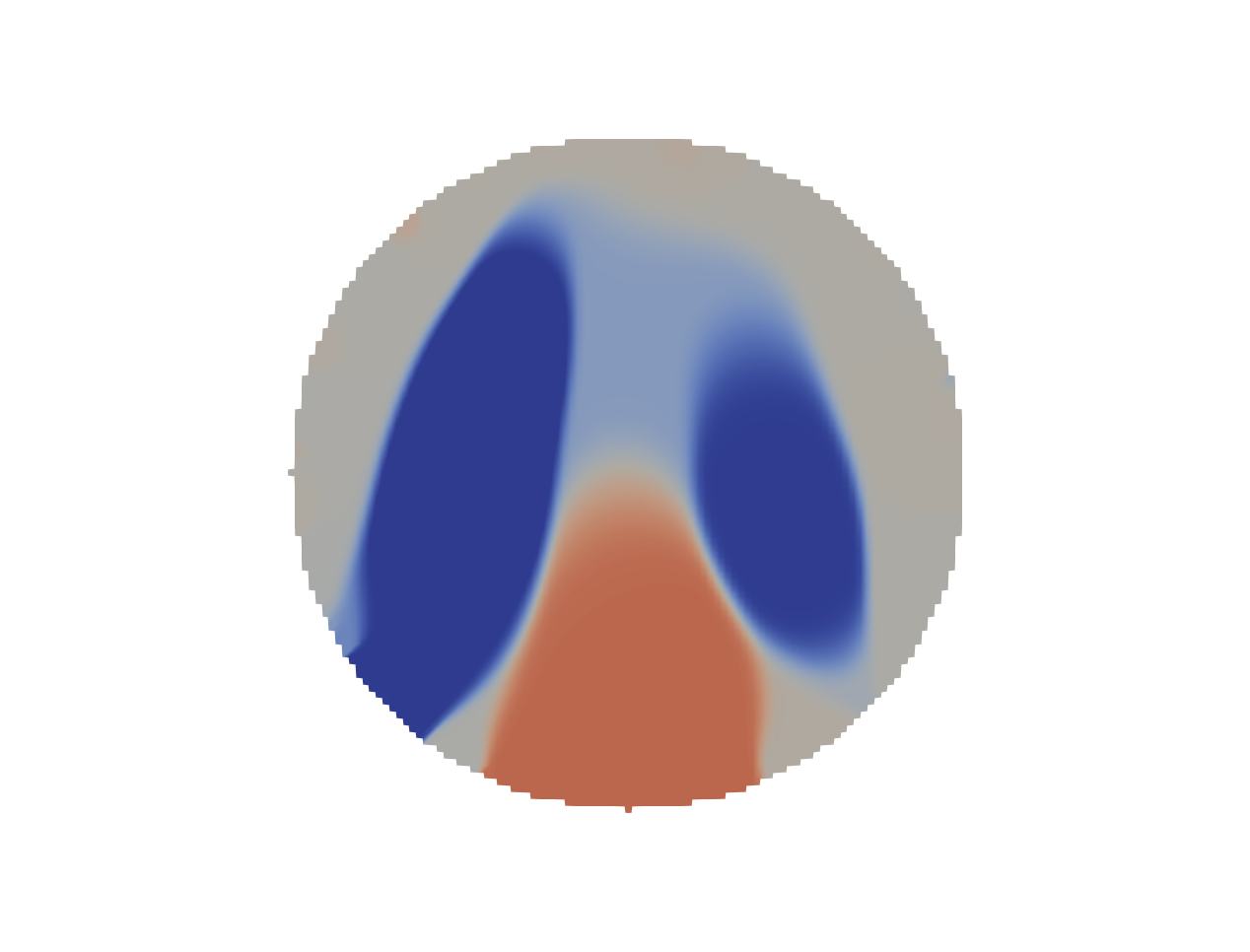}\includegraphics[scale = 0.1, trim={6cm 0 9.5cm 5cm},clip, valign=t]{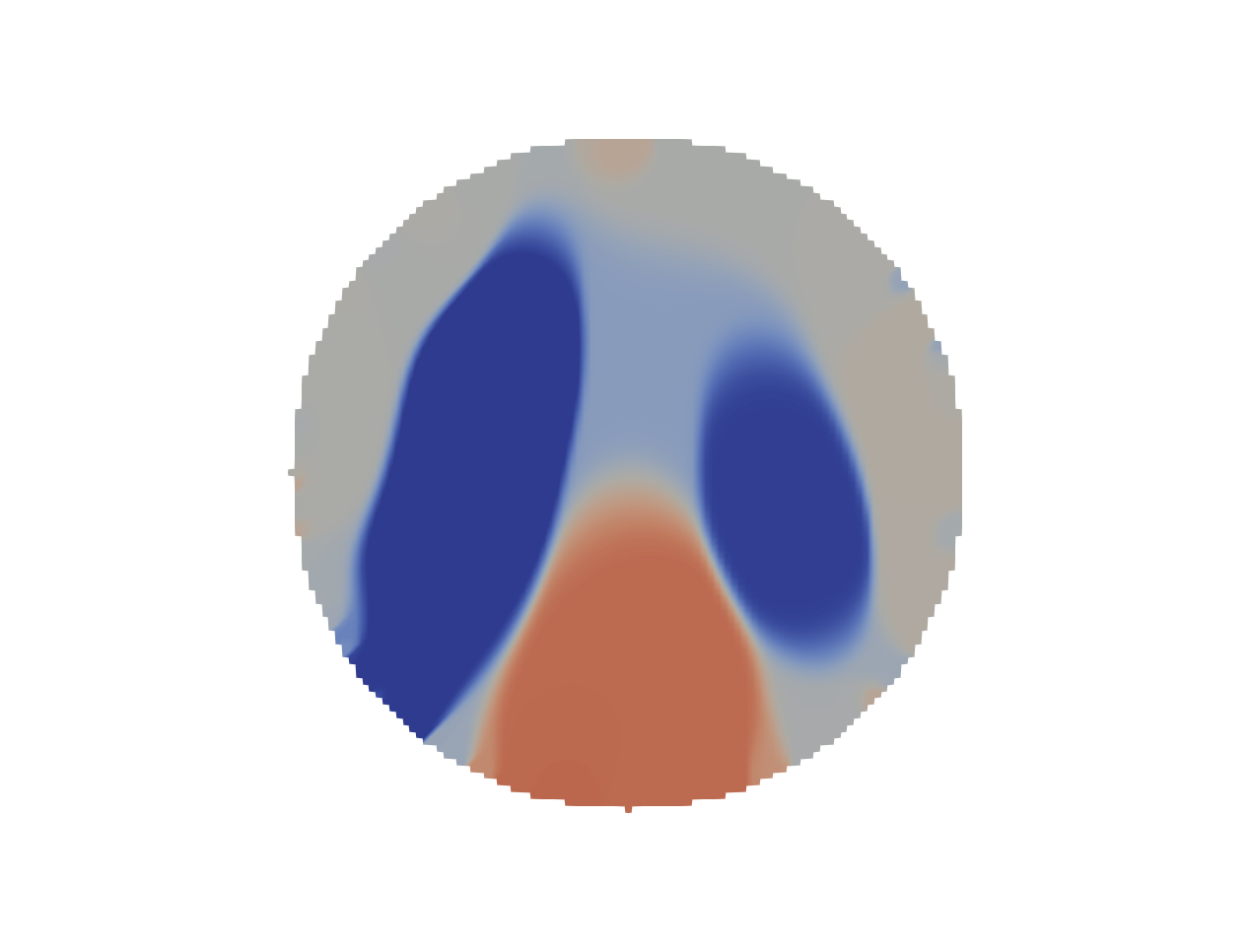}\includegraphics[scale = 0.20, trim={35cm 0 0 15cm},clip, valign=t]{eit_figures/mon_ref4.png}}
    \caption{ADMM and monolithic reconstructions over multiple mesh refinements}
    \label{fig:prob}
\end{figure}

\subsubsection{Computational scalability with respect number of PDE models}

\begin{figure}[tbh]
    \centering
    \includegraphics[scale = 0.3]{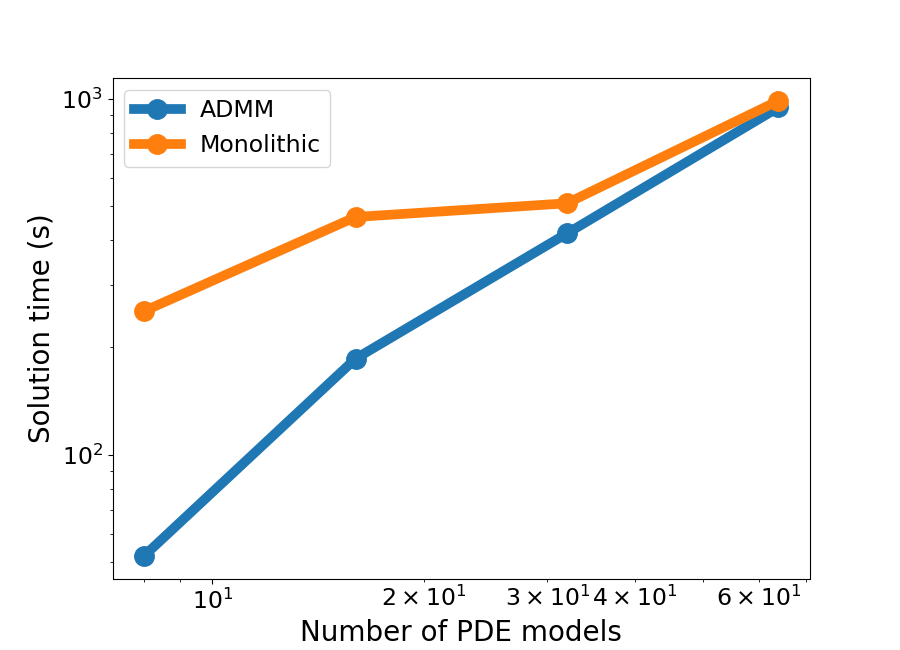}
    \caption{Solution time(s) vs number of PDE models}
    \label{fig:prob_time}
\end{figure}

\begin{figure}[tbh]
    \centering
    \resizebox{\textwidth}{!}{\includegraphics[scale = 0.4, trim={3cm 0 3.5cm 0},clip]{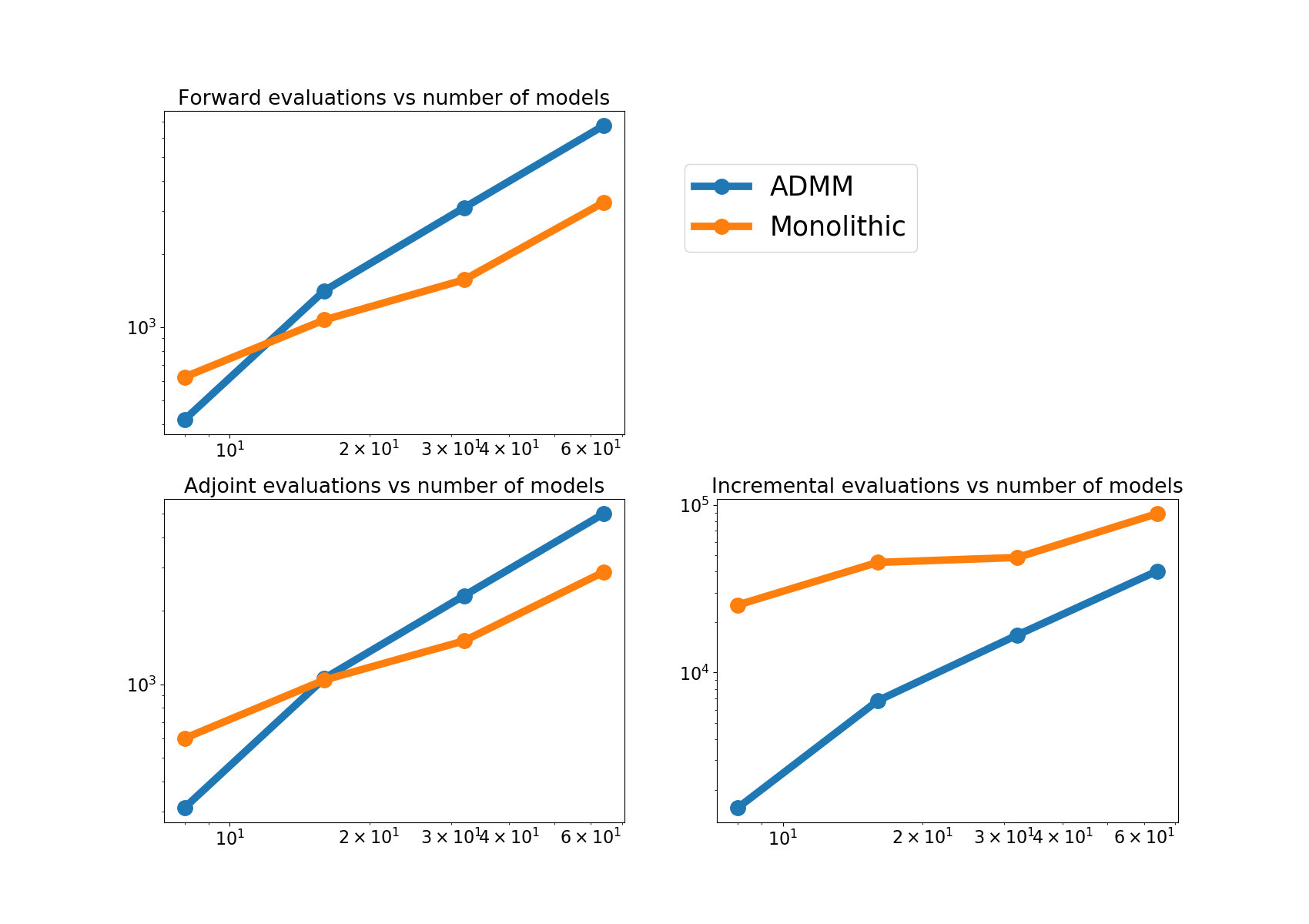}}
    \caption{Number of forward solves, adjoint solves, and incremental solves for ADMM and monolithic approaches vs number of PDE models}
    \label{fig:comp_prob}
\end{figure}

 The time required for each reconstruction is plotted against the number of PDE models in Figure \ref{fig:prob_time}. We continue to assess computational cost based on the number of forward, adjoint, and incremental solves. 
 The total number of forward solves, adjoint solves, and incremental solves is plotted against the number of PDE models in Figure \ref{fig:comp_prob}. 

The ADMM still achieves satisfactory accuracy in a shorter amount of time for each number of PDE models. Once again, this is is because it requires much fewer incremental evaluations for each solution. These results are summarized in Table \ref{Tab:prob_table}.

\begin{table}[tbh]
\caption{Runtime comparison between ADMM and monolithic approach as a function of the number of forward models} 
\label{Tab:prob_table}

\resizebox{\textwidth}{!}{
\begin{tabular}{ l|llll } 
 \textnormal{PDE models}& 8 & 16  & 32 &  64 \\
 \hline
 \textnormal{Solution time} & 52s/4m 13s & 3m 6s/7m 45s &6m 58s/8m 28s & 15m 44s /16m 26s\\ 
 
 \textnormal{Relative error} & 0.1261/ 0.0981 & 0.1116/0.0817 & 0.1059/0.0817& 0.1033/0.0707\\ 
 \textnormal{State misfit} & 1.284/1.127 & 1.221/1.242 &1.264/1.266 & 1.244/1.307\\ 
 \textnormal{Forward solves} & 417/624 & 1414/1072 &3104/1568 & 6758/3264\\
 \textnormal{Adjoint solves} & 312/600 &1056/1040 & 2304/1504& 4992/2880\\ 
 \textnormal{Incremental solves} & 1566/25312 & 6792/46360  & 16748/48544& 40278/88576\\

\end{tabular}
}
\end{table}

\section{Application to multi-wavelength quantitative photoacoustic tomography problem}\label{sec:qpact}

This section presents an application of the proposed ADMM method to quantitative photoacoustic tomography(qPACT).
qPACT) is an emerging medical imaging technique that holds great promise for early cancer diagnosis because it is non-invasive, radiation-free, and inexpensive. qPACT is a hybrid modality that combines endogenous contrast of optical imaging with the high-resolution of ultrasound detection technologies to provide maps of total hemoglobin content and oxygen saturation within the tissue\cite{WangAnastasio15, KloseLarsen06}.

This process requires solving a series of two inverse problems and thus involves two separate stages. The process begins with a fast laser pulse in the infrared range being sent into the object of interest. The underlying material then absorbs this optical energy and generates heat and a local increase pressure distribution. This pressure distribution then transitions into acoustic waves that can then be observed at the boundary of the domain. This process can be viewed as two inverse problems. The first involves reconstructing for the initial pressure distribution given measurements on the boundary of the domain. The second involves reconstructing the optical properties of the tissue based on the initial pressure distribution. Here we only worry about the second inverse problem of reconstructing tissue composition given the initial pressure distribution.
\subsection{Formulation of the qPACT problem}
Let $p_0$ denote this initial pressure distribution. $p_0$ is related to the absorption coefficient $\mu_a$ of the domain via \eqref{eqn:data_mua}.

\begin{equation}\label{eqn:data_mua}
    d= \frac{p_0}{\Gamma} =  \mu_a \phi
\end{equation}
 $\Gamma$ is the Grunesian parameter, which we will treat as a known value of constant 1.The fluence $\phi$ is then determined as the solution of the PDE relationship below, known as the diffusion approximation.

\begin{equation}\label{eqn:diffusion_approx}
    \begin{array}{cc}
        -\nabla\cdot \frac{1}{3(\mu_a + \mu_s')}\nabla\phi + \mu_a\phi = 0 & x \in \Omega \\
        \frac{1}{3(\mu_a + \mu_s')}\frac{\partial\phi}{\partial \eta} + \frac{1}{2}\phi = \frac{1}{2}\phi_0 & x \in \partial \Omega
    \end{array}
\end{equation}

where $\mu_s'$ is the reduced scattering coefficient and $\phi_0$ is the intensity of the incident 
illumination. Here $\mu_s' = (1-g)\mu_s$ where $\mu_s$ is scattering coefficient and $g$ is the anisotropy 
factor. We will fix $g = .9$. The absorption coefficient $\mu_a$ is a linear combination of basis materials 
called chromophores and reads
\begin{equation}\label{eqn:mua_chrom}
    \mu_a = \sum_i \varepsilon_i(\lambda)c_i
\end{equation}
 $\varepsilon_i$ is a known function of the incident wavelength and $c_i$ is the concentration of the $i$-th chromophore.

We will be attempting to reconstruct for a 2d maximum intensity projection of a human breast. Using human breast tissues as our domain means that we can safely limit our chromophores to deoxygenated hemoglobin $c_{hb}$ and oxygenated hemoglobin $c_{hb0_2}$. These two values are then related to the oxygen saturation $s$ and the total hemoglobin concentration $c_{thb}$ via \eqref{eqn:hb_calcs}.

\begin{equation}\label{eqn:hb_calcs}
    \begin{array}{c}
        c_{hb} = (1-s)\cdot c_{thb}   \\
        c_{hb0_2} = s\cdot c_{thb} 
    \end{array}
\end{equation}

We can then invert for $s,\ c_{thb},$ and $\mu_s$ with a PDE model associated with each incident wavelength. For each PDE model, there will be an associated fluence $\phi_i$, which acts as the state variable. Letting $d_i$ be the observations associated with the $i$-th incident wavelength, we can then form the data fidelity term

$$\frac{1}{q}\sum_{i=1}^d \loss_i(s,c_{thb}, \mu_s') = \frac{1}{q}\sum_{i=1}^d ||\ln(\mu_{a,i}\phi_i) - \ln(d_i)||^2$$
\subsection{Inversion results}
\begin{figure}[tbh]
    \centering
    \resizebox{\textwidth}{!}{\includegraphics[scale = 0.15, trim={10cm 0 5cm 0 },clip]{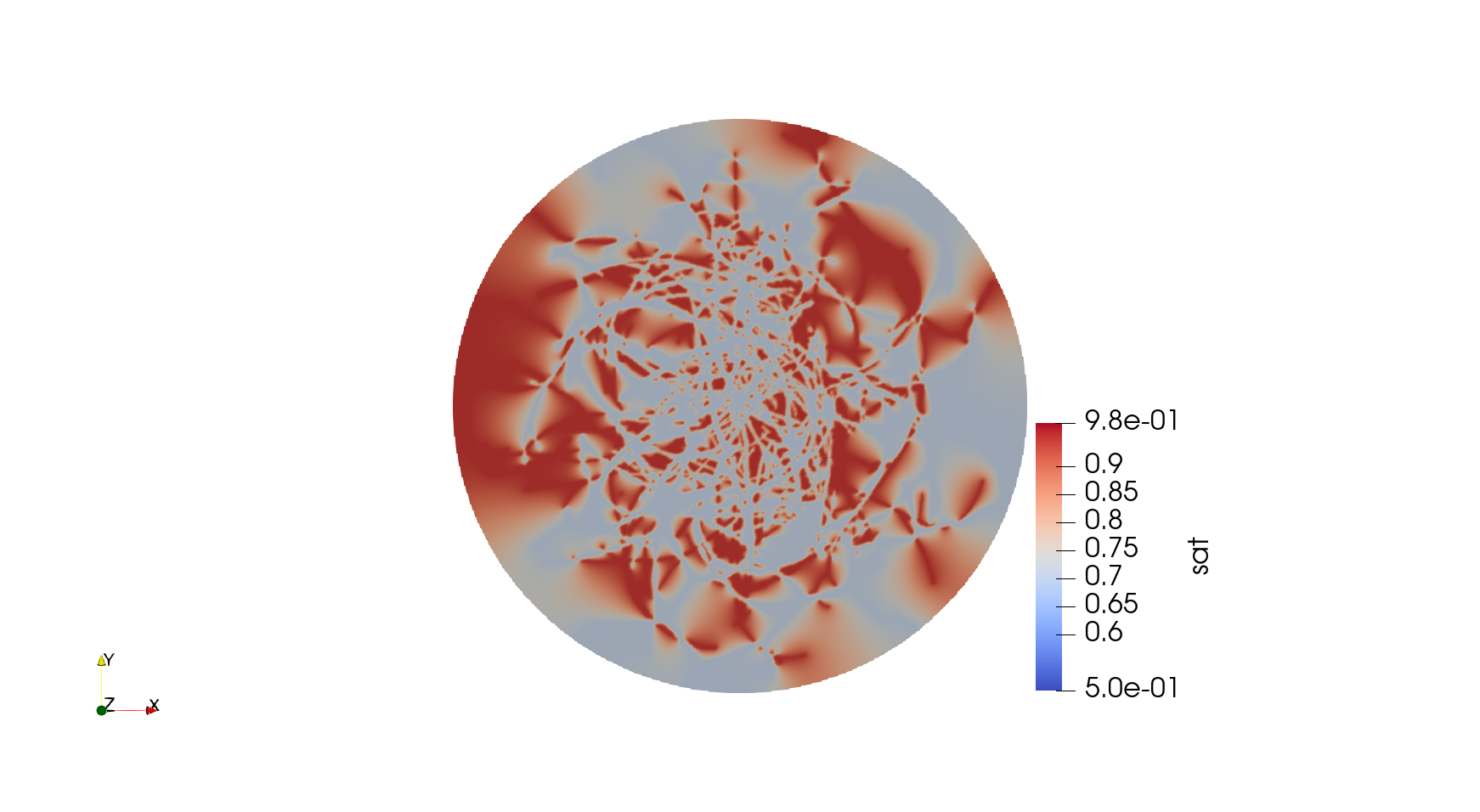}\includegraphics[scale = 0.15, trim={10cm 0 5cm 0 },clip]{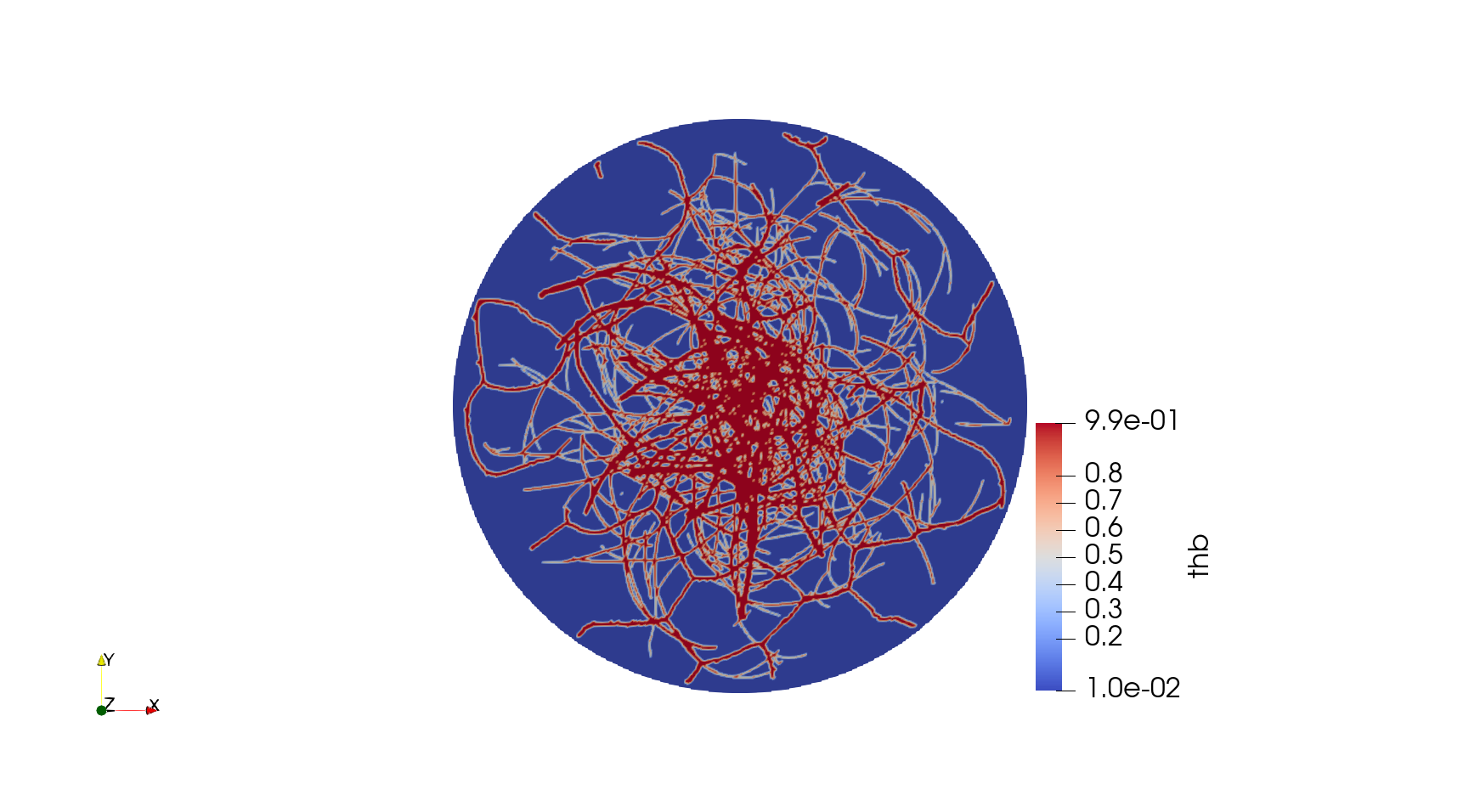}\includegraphics[scale = 0.15, trim={10cm 0 5cm 0 },clip]{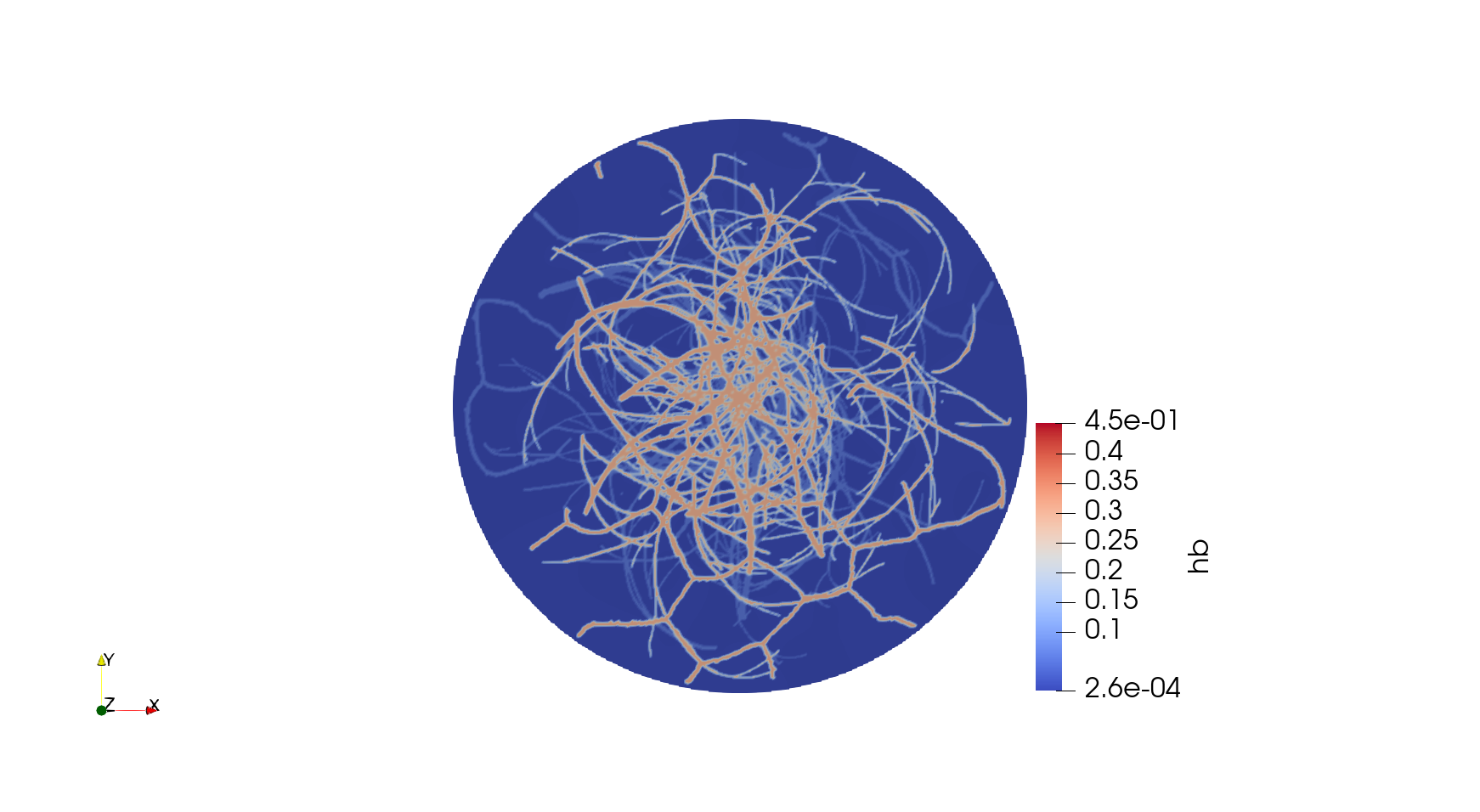}\includegraphics[scale = 0.15, trim={10cm 0 5cm 0 },clip]{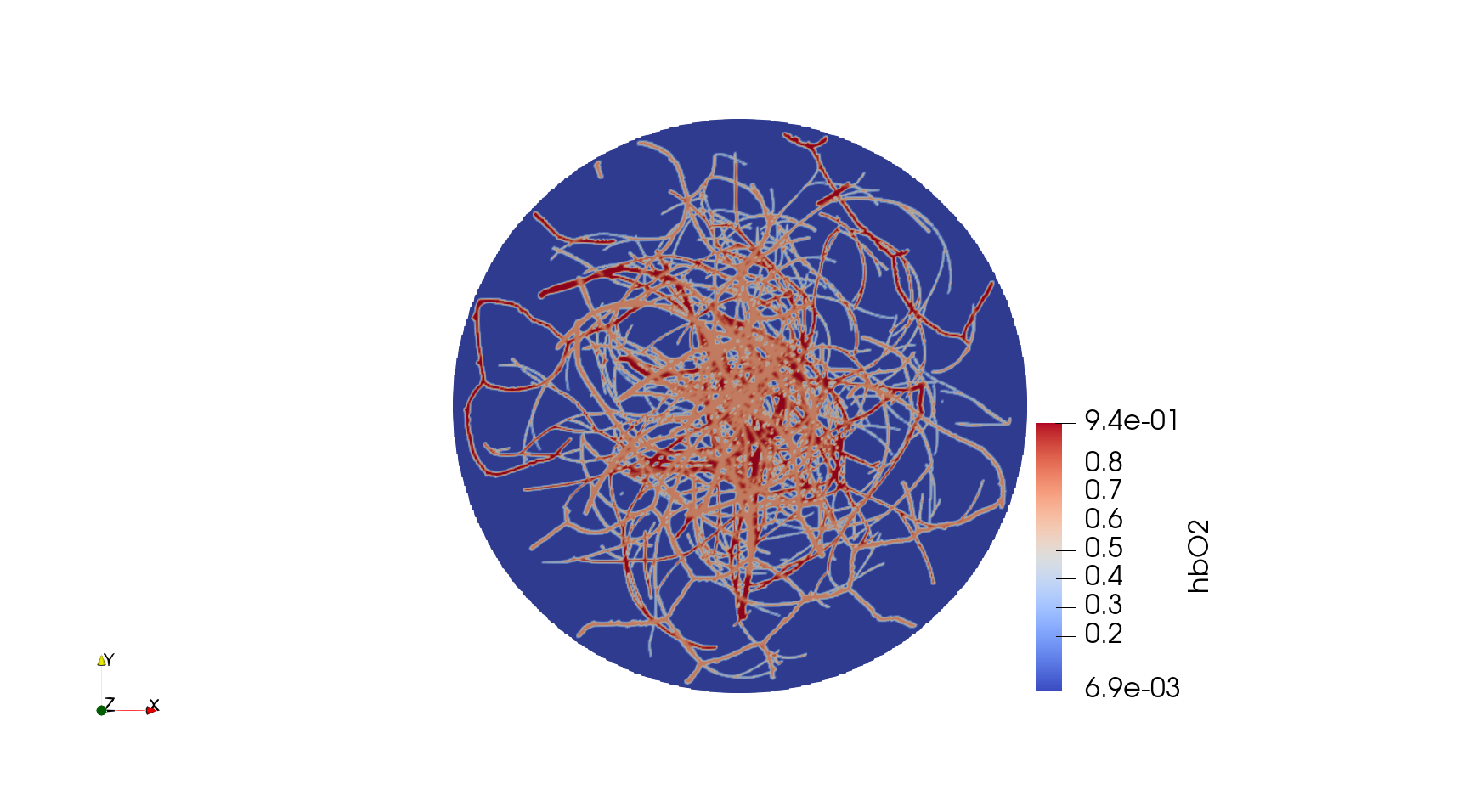}}
    \caption{True $s, \ c_{thb}, \ c_{hb},$ and $c_{hbO_2}$}

    \centering
    \resizebox{\textwidth}{!}{\includegraphics[scale = 0.15, trim={10cm 0 5cm 0 },clip]{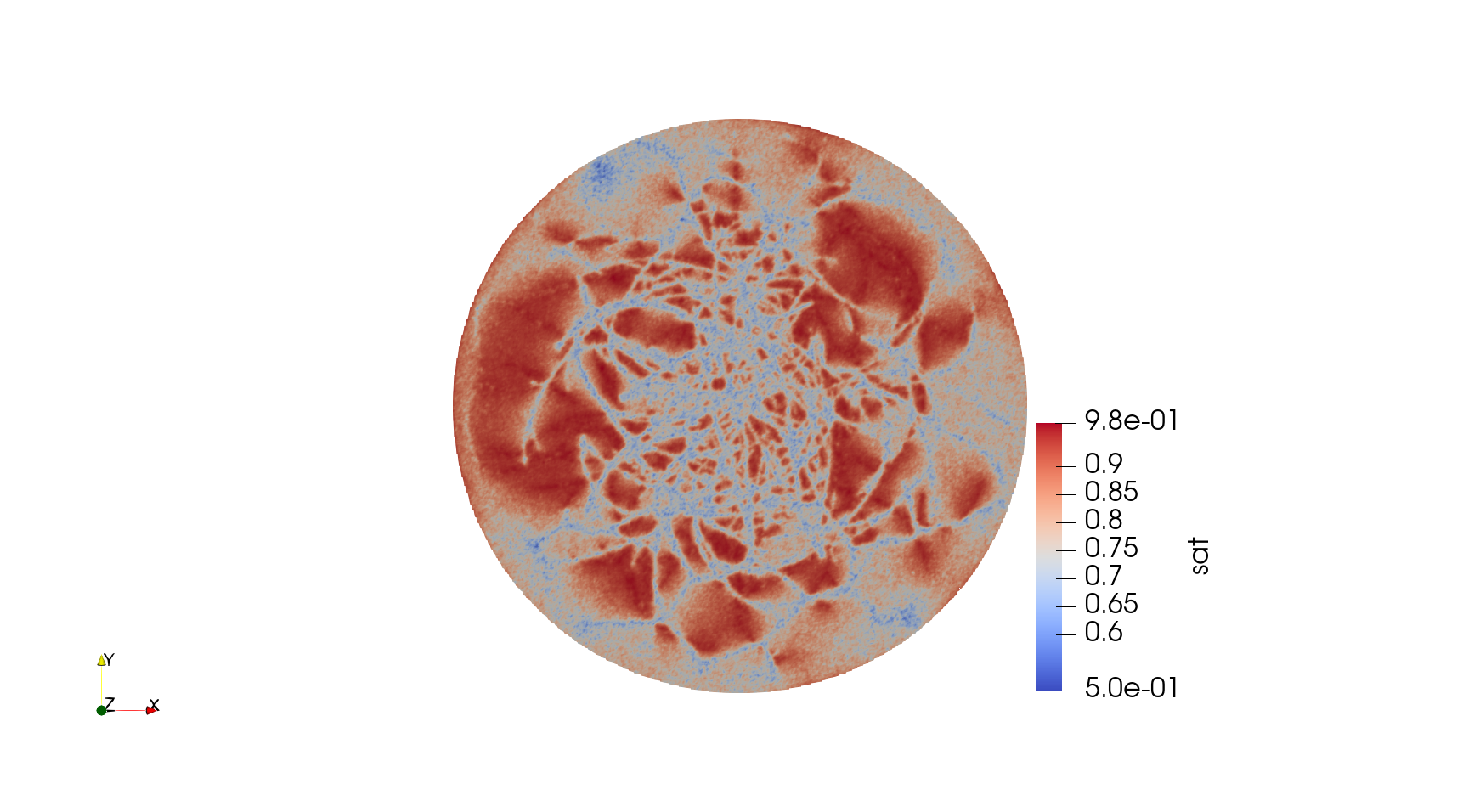}\includegraphics[scale = 0.15, trim={10cm 0 5cm 0 },clip]{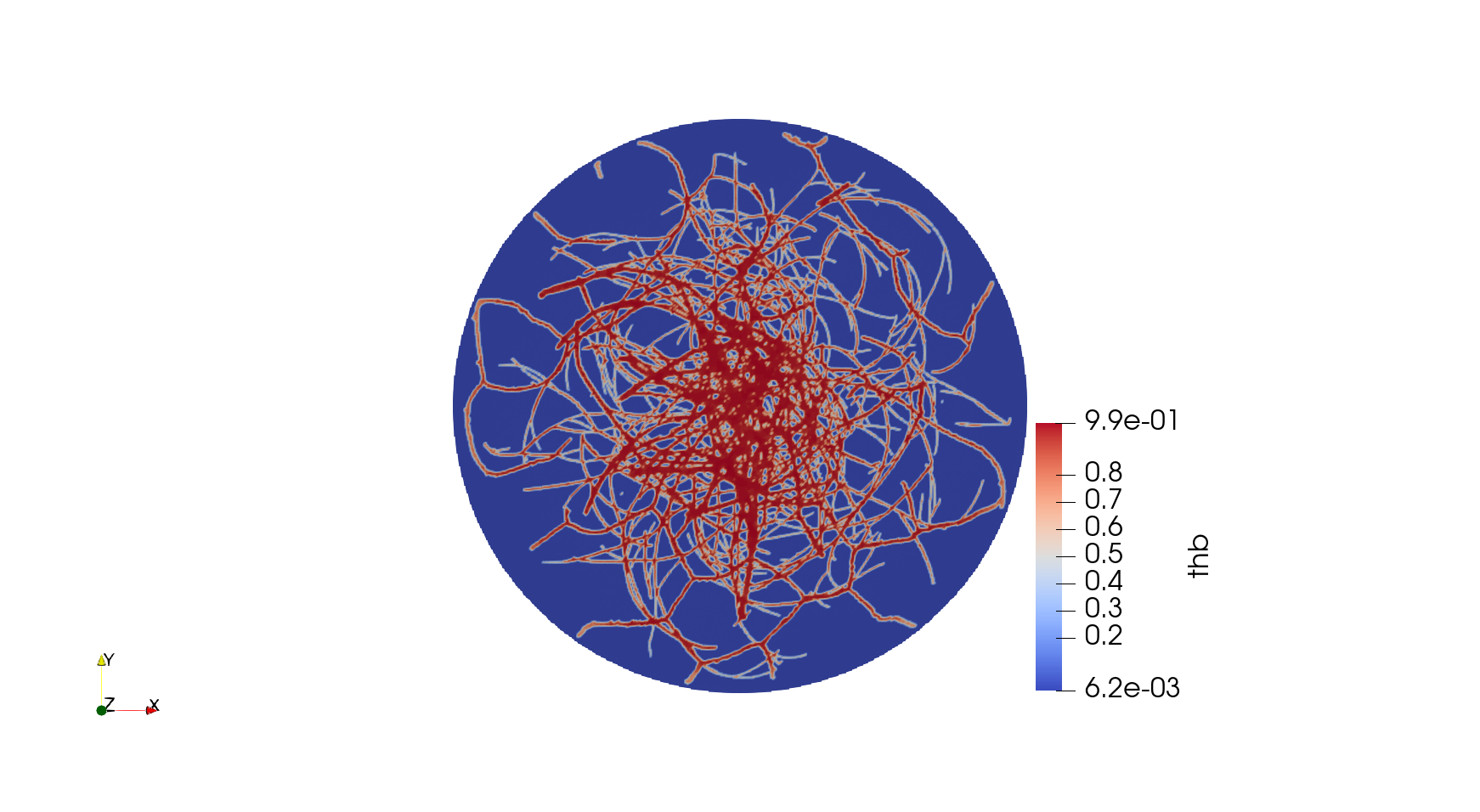}\includegraphics[scale = 0.15, trim={10cm 0 5cm 0 },clip]{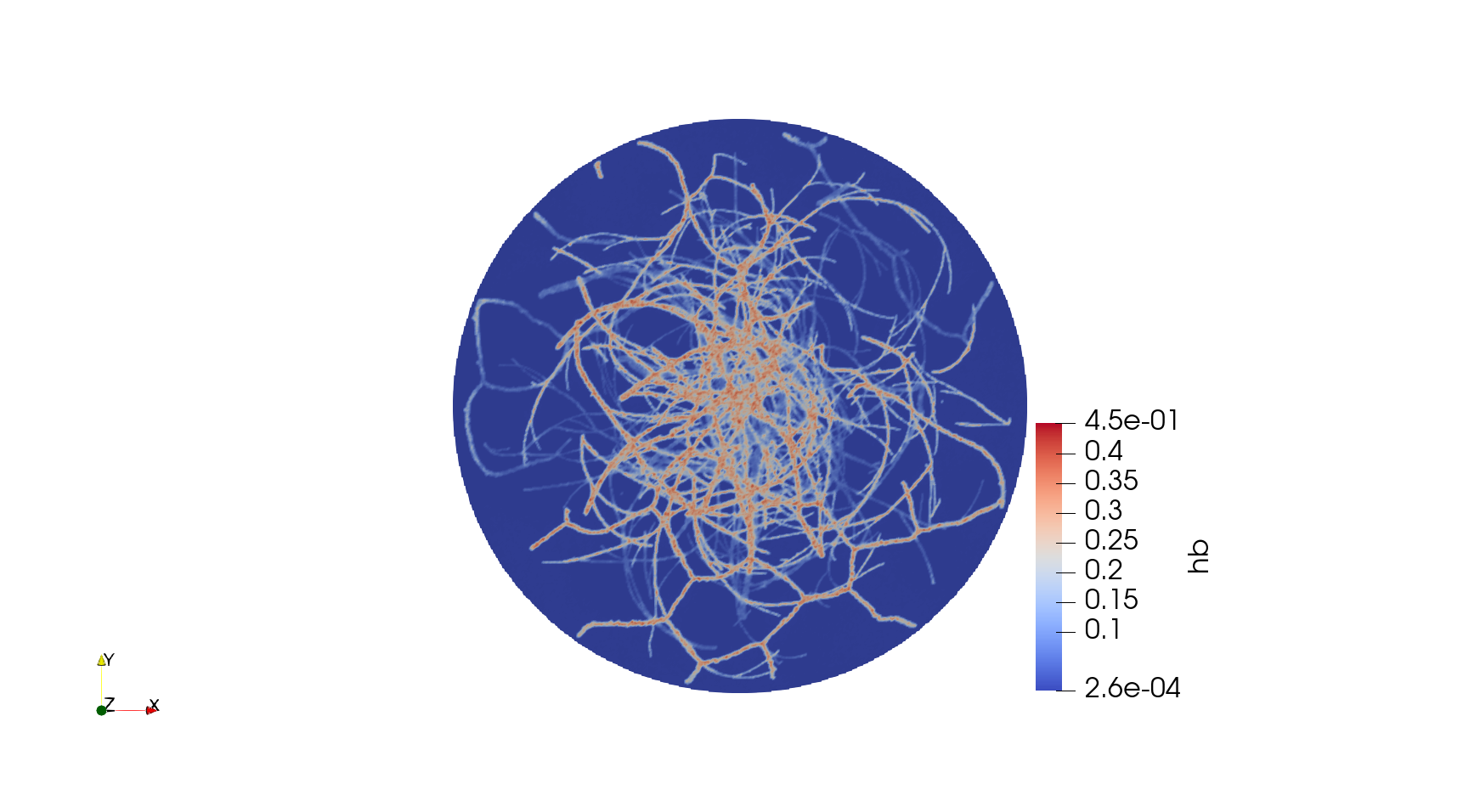}\includegraphics[scale = 0.15, trim={10cm 0 5cm 0 },clip]{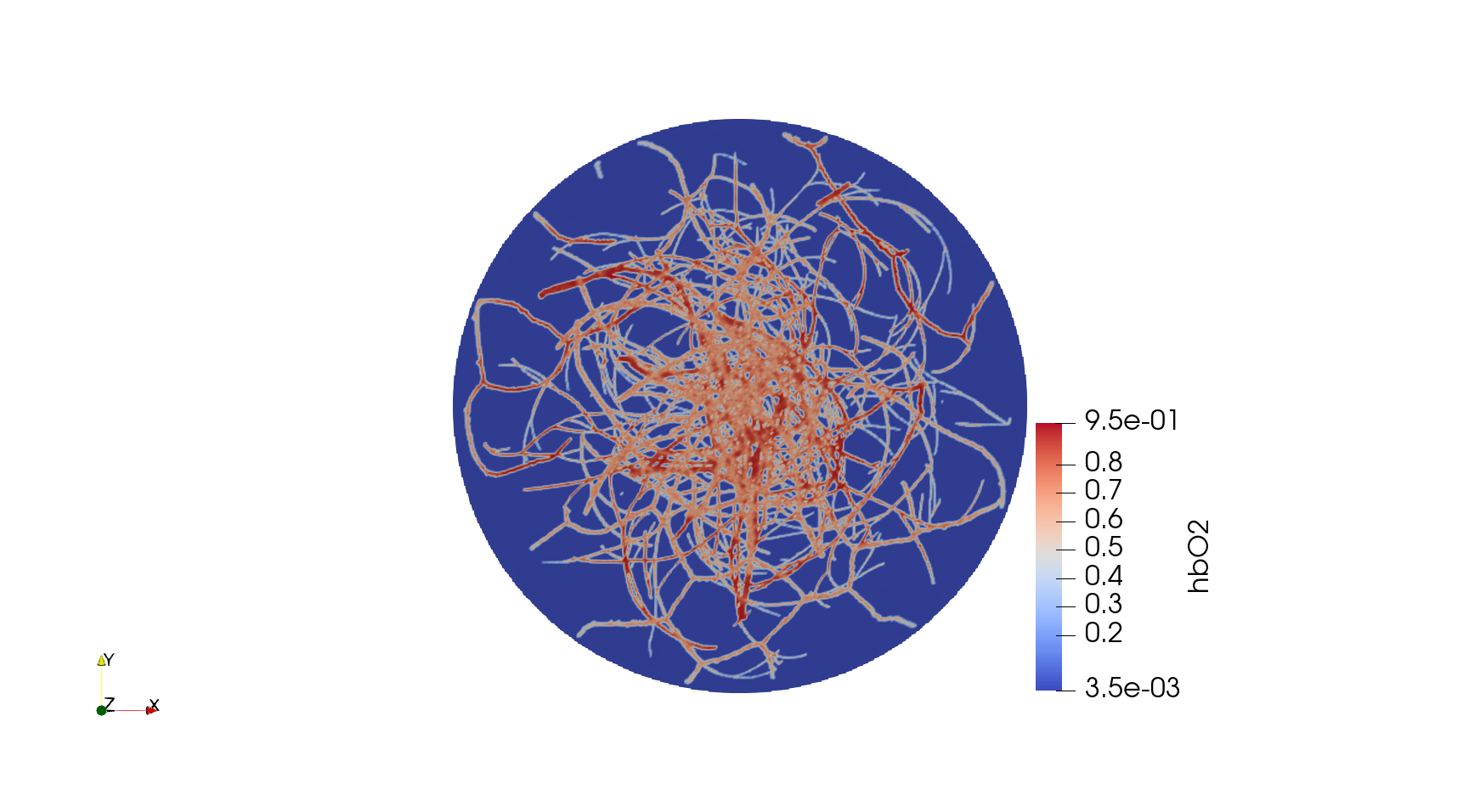}}
    \caption{Reconstructed $s, \ c_{thb}, \ c_{hb},$ and $c_{hbO_2}$}
    \label{fig:qpact}
\end{figure}

We then perform an inversion with measurements from $757,\ 800,$ and $850$ wavelengths with uniform intensity around the 2d boundary.The ground truth values for $s,\ c_{thb}, \ c_{hb},$ and $c_{hb0_2}$ are shown in the top row of  Figure \ref{fig:qpact}.We ran forward model for each wavelength and perturbed each of the data measurements by white Gaussian noise with a standard deviation equal to a hundredth of their maximum value. We then used the regularization shown in \eqref{eqn:qpact_reg}.
\begin{equation}\label{eqn:qpact_reg}
    \begin{array}{lrc}
        R(s,c_{thb},{\mu_s'}) =& \gamma_{s}\int_{\Omega}||\nabla s||^2 d \boldsymbol x&  + \delta_{s}\int_{\Omega}||s||^2 d \boldsymbol{x} \\ 
         & \gamma_{c_{thb}}\int_{\Omega}|\nabla c_{thb} |_\eps d \boldsymbol x
         & + \delta_{c_{thb}}\int_{\Omega}| c_{thb} |_\eps d \boldsymbol x\\
         
         & \gamma_{\mu_s'}\int_{\Omega}||\nabla \mu_s'||^2 d \boldsymbol x&  + \delta_{\mu_s'}\int_{\Omega}||\mu_s'||^2 d \boldsymbol{x}
    \end{array}
\end{equation}
Using ADMM is of much interest because it can iteratively solve stable subproblems to continuously advance to the true solution while appropriately handling the complex regularization function. We then applied the ADMM to solve this problem with the three different wavelength and a regularization determined by
$\gamma_s =0.05, \delta_s = 0.001,\gamma_{c_{thb}} = 0.005, \delta_{c_{thb}} = 10^{-6}, \gamma_{\mu_s'} =10,\delta_{\mu_s'}= 10,$ and$ \varepsilon = 10^{-6}$.
In doing this, we set an absolute global tolerance of $10^{-4}$ and a relative global tolerance of $10^{-3}$. The optimal values of $\ipar^k$ were calculated using an INCG solver with at most $50$ iterations, a relative tolerance of $1e{-6}$, and an absolute tolerance of $10^{-9}$. The number of iterations on this solver was more than needed and always converged within $30$ iterations.
The consensus variable was updated using a PETScTAOSolver solver with a relative tolerance of $1e{-9}$ and an absolute tolerance of $1e{-12}$. We implemented the adaptive $\rho$ scheme with $\mu = 4$ and $\tau = 2$. We were then able to reconstruct for the parameters shown in the bottom row of Figure \ref{fig:qpact}.

This reconstruction was very accurate and had minimal errors. The global relative error on the saturation was only $0.061$, and the global relative error on the total hemoglobin concentration was only $.097$. These errors were even lower on the arteries and veins, these places with high hemoglobin content and saturation of $0.95$ for arteries and $0.7$ for veins. The percentage errors on the arteries were $0.040$ for the saturation and $0.061$ for the total hemoglobin content. On the veins, the percentage errors were $0.050$ and $0.055$. These low errors mean that we can accurately identify arteries and veins using this reconstruction process. The error on this reconstruction is summarized in Table \ref{Tab:qpact_error}.

\begin{table}[tbh]
\caption{Reconstruction errors for the qPACT problem}
\label{Tab:qpact_error}
 
($s$/$c_{thb}$)\\
\resizebox{\textwidth}{!}{\begin{tabular}{ l|lll} 
Region &	$L^2$ norm & error & relative error \\
\hline 
Global &	82.73101184/48.23308111	& 5.041813126/4.669520028&	0.060942239/0.096811564\\
Artery &	30.3685615/27.48296848 &	1.207798494/1.685043185 &0.03977134/0.06131227 \\
Vein &	30.86117651/37.2034497 &	1.543544399/2.028329233 &	0.050015734/0.054519924 \\
Inner &	46.56081498/39.39550923	& 2.364520795/2.940342136	&0.050783492/0.07463648\\
Outer &	68.38501903/27.82847402	& 4.452967663/3.627506777	&0.065116128/0.130352342

\end{tabular}}
\end{table}

\begin{figure}[tbh]
    \centering
    \resizebox{\textwidth}{!}{\includegraphics[scale = 0.4]{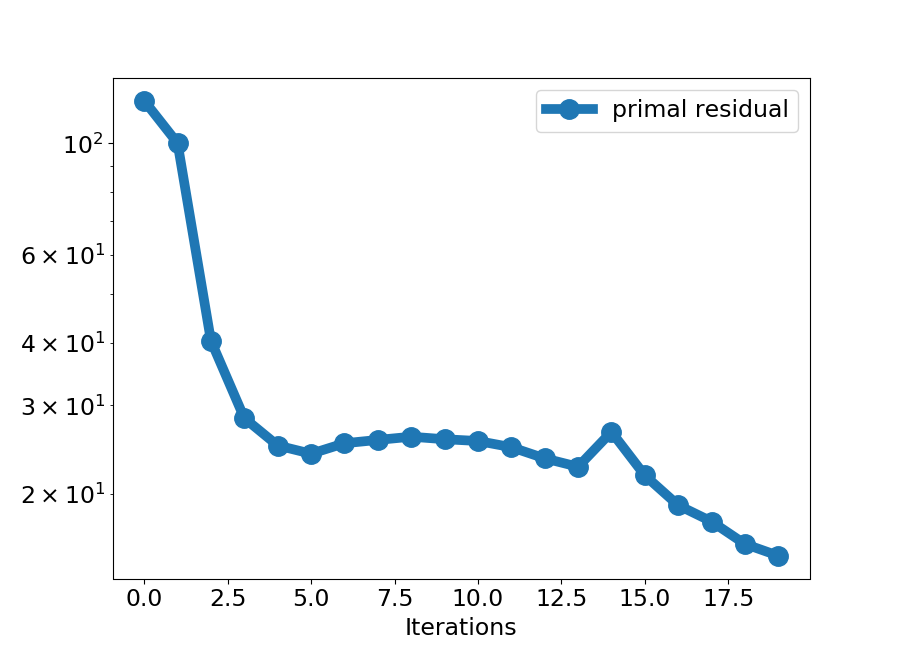}\includegraphics[scale = 0.4]{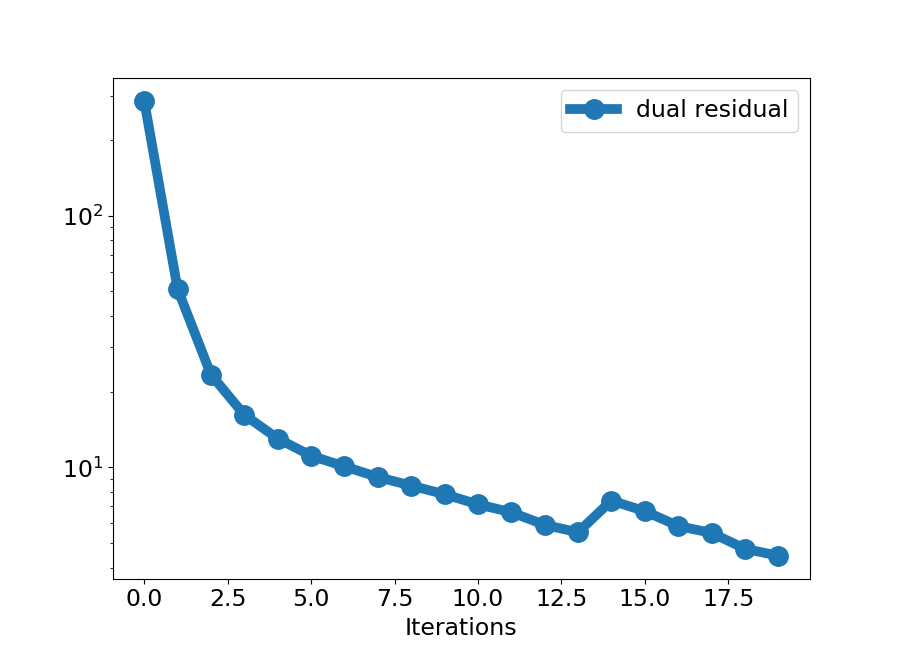}}
    \caption{Primal and dual residuals from the qPACT reconstruction}
    \label{fig:qpact_residuals}
\end{figure}

\begin{figure}[tbh]
    \centering
    \includegraphics[scale = 0.3]{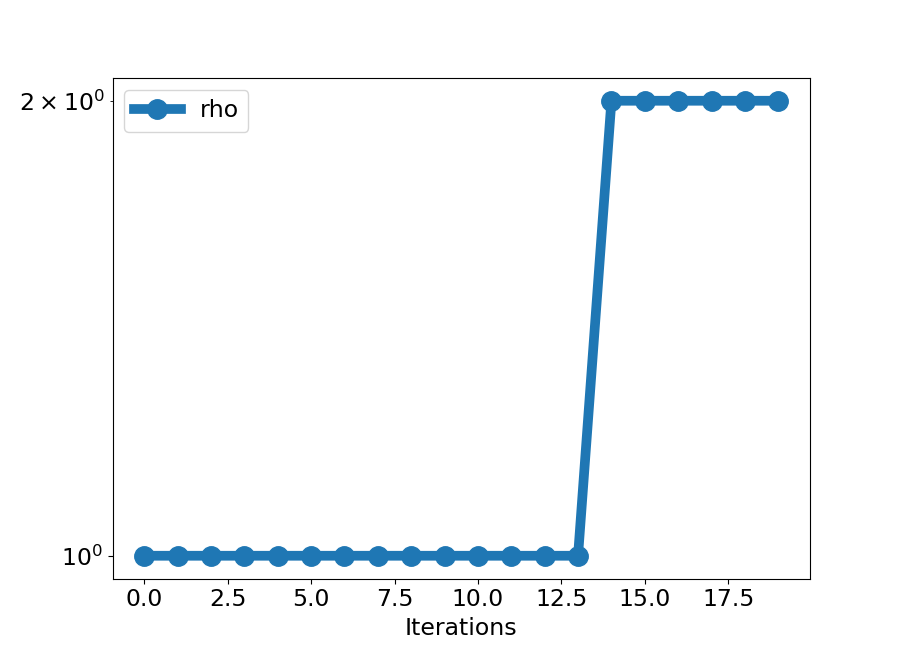}
    \caption{$\rho$ values of each iteration of qPACT reconstruction}
    \label{fig:qpact_rho}
\end{figure}

This reconstruction saw a continuous decrease in the residuals with each iteration. These residuals are shown in Figure \ref{fig:qpact_residuals}, and the $\rho$ for each iteration is shown in Figure \ref{fig:qpact_rho}.
\section{Conclusions}
In this paper, we presented a framework for solving \newline inverse problems governed by PDE forward models using ADMM. Through our numerical studies with electrical impedance tomography, we have shown the natural way in which ADMM can solve problems involving several large-scale PDE models with nonsmooth regularization. This ADMM solution method significantly reduced these solutions' computational costs while still achieving satisfactory accuracy at various scales. This framework preserves consistency with the infinite formulation of these inverse problems and utilizes the underlying Hilbert spaces' norm to enforce the consensus condition. The 
effectiveness of the ADMM framework was also demonstrated on a complex multiphysics problem related to 
photoacoustic tomography. Solving the photoacoustic tomography problem accurately and efficiently shows the 
viability of the ADMM framework outside the model problem related to electrical impedance tomography. In the 
future, we plan to improve upon this framework in two ways. First, we plan on implementing more advanced solvers for the total variation denoising problem, including the 
primal-dual method in \cite{ChanGolubMulet99} and the proximal splitting methods in \cite{HerrmannHerzogSchmidtEtAl19}. Second, we plan on implementing the ADMM process on 
several processors, with each PDE model being handled by its own set of processors. Splitting PDE models 
along processor, sets would allow every parameter associated with a PDE model to be updated simultaneously and accelerate the entire process. 

\bibliographystyle{siam}
\bibliography{local, references}

\end{document}